\documentclass[11pt]{article}
\usepackage[mathscr]{eucal}
\DeclareMathAlphabet{\mathpzc}{OT1}{pzc}{m}{it}
\usepackage{bbm}
\usepackage{dsfont}
\usepackage{color}
\usepackage{mathrsfs}
\usepackage{amssymb}
\usepackage{calligra}
\usepackage{fontenc}
\usepackage{amsbsy}
\newcommand{\bfscr}[1]{{\pmb{\mathscr{#1}}}}
\newcommand{\bfpzc}[1]{{\pmb{\mathpzc{#1}}}}
\newcommand{\bfsf}[1]{{\textbf{\textsf{#1}}}}
\newcommand{\sub}[1]{{\mbox{\footnotesize $#1$}}}
\usepackage{graphicx}
\graphicspath{%
    {converted_graphics/}
    {/}
    {C:/Users/GALDI/Desktop/Vibration/}
}
\begin{document}

\newcommand{\sfW}{{\sf W}}
\newcommand{\sfD}{{\sf D}}
\newcommand{\sfL}{{\sf L}}
\newcommand{\sfw}{{\sf w}}
\newcommand{\sfp}{{\sf p}}
\newcommand{\simge}{\ba{cc}\vspace*{-2.4mm}>\\ \sim\ea }
\newcommand{\simle}{\ba{cc}\vspace*{-2.4mm}<\\ \sim\ea }
\newcommand{\Cdot}{\!\cdot\!}
\newcommand{\sq}{{$\sqcap\!\!\!\!\sqcup$}}
\newcommand{\Eu}{{\rm I\,\!\! E}}
\newcommand{\Io}{\Int{\Omega}{}}
\newcommand{\Id}{\Int{\cald}{}}
\newcommand{\Div}{\mbox{\rm div}\,}
\newcommand{\tr}{\mbox{\rm tr}\,}
\newcommand{\grad}{\mbox{\rm grad}\,}
\newcommand{\supp}{\mbox{\rm supp}\,}
\newcommand{\curl}{\mbox{\rm curl}\,}
\newcommand{\Ido}{\Int{\partial\Omega}{}}
\newcommand{\IdS}{\Int{\Sigma}{}}
\newcommand{\Oint}[2]{{\displaystyle \oint_{#1}^{#2}}}
\newcommand{\Int}[2]{{\displaystyle \int_{ #1}^{ #2}}}
\newcommand{\Lim}[1]{{\displaystyle \lim_{ #1}}}
\newcommand{\Limsup}[1]{{\displaystyle \limsup_{\footnotesize #1}}}
\newcommand{\Liminf}[1]{{\displaystyle \liminf_{\footnotesize #1}}}
\newcommand{\Sup}[1]{{\displaystyle \sup_{#1}}}
\newcommand{\Inf}[1]{{\displaystyle \inf_{#1}}}
\newcommand{\Max}[1]{{\displaystyle \max_{#1}}}
\newcommand{\Min}[1]{{\displaystyle \min_{#1}}}
\newcommand{\Sum}[2]{{\displaystyle \sum_{#1}^{#2}}}
\newcommand{\Prod}[2]{{\displaystyle \prod_{#1}^{#2}}}
\newcommand{\BCup}[2]{{\displaystyle \bigcup_{#1}^{#2}}}
\newcommand{\BCap}[2]{{\displaystyle \bigcap_{#1}^{#2}}}
\newcommand{\Frac}[2]{\displaystyle{\frac{\displaystyle{#1}}{\displaystyle{#2}}}}
\newcommand{\norm}[1]{\left\|{#1}\right\|}
\newcommand{\Norm}[1]{\langle\langle{#1}\rangle\rangle_q}
\newcommand{\No}[1]{\langle\!\langle{#1}\rangle\!\rangle}
\newcommand{\NO}[1]{{\langle{#1}\rangle}_{\lambda,q}}
\newcommand{\beea}{\begin{eqnarray}}
\newcommand{\eeea}{\end{eqnarray}}
\newcommand{\ms}{\medskip\smallskip}
\newcommand{\bs}{\bigskip}
\newcommand{\ps}{\par\smallskip}
\newcommand{\bfe}{{\mbox{\boldmath $e$}} }
\newcommand{\pni}{{\par\noindent}}
\newcommand{\bfq}{{\mbox{\boldmath $q$}} }
\newcommand{\bfz}{{\mbox{\boldmath $z$}} }
\newcommand{\0}{{\mbox{\boldmath $0$}} }
\newcommand{\LE}{\!\!\!&\le&\!\!\!}
\newcommand{\BL}[1]{{\par\smallskip{\bf Lemma #1.}}}
\newcommand{\BT}[1]{{\par\smallskip{\bf Theorem #1.}}}
\newcommand{\Ln}{[\!|}
\newcommand{\Rn}{|\!]}
\newcommand{\n}[1]{{\Ln{#1}\Rn}} 
\newcommand{\nq}[1]{{\Ln{#1}\Rn}_{q}} 
\newcommand{\nqr}[1]{{\Ln{#1}\Rn}_{q,r}} 
\newcommand{\Nq}[1]{{\langle{#1}\rangle}_{q}} 
\newcommand{\Nql}[1]{{\langle{#1}\rangle}_{\lambda,q}} 
\newcommand{\Nqr}[1]{{\langle{#1}\rangle}_{q,r}}
\newcommand{\N}[1]{{|\!\!|\!\!|\,{#1}\,|\!|\!\!|_2}}
\newcommand{\EA}[2]{$$#1$$%
\vspace{-6.mm}
\begin{equation}
\end{equation}
\vspace{-6.mm}
$$
#2
\setlength{\belowdisplayskip}{3mm}
\setlength{\belowdisplayshortskip}{3mm}
$$
}
\newcommand{\A}[2]{$$#1$$%
\vspace{-4.mm}
$$
#2
\setlength{\belowdisplayskip}{3mm}
\setlength{\belowdisplayshortskip}{3mm}
$$
}
\newcommand{\BF}{\begin{footnotesize}}
\newcommand{\EF}{\end{footnotesize}}
\setlength{\jot}{.15in}
\newcommand{\pde}[2]{{\displaystyle \frac{\mbox{$\partial #1$}}{\mbox{$\partial #2$}}}}
\newcommand{\ode}[2]{{\displaystyle \frac{\mbox{$d #1$}}{\mbox{$d #2$}}}}
\newcommand{\f}[2]{\frac{\mbox{$#1$}}{\mbox{$ #2$}}}
\newcommand{\bi}{\begin{itemize}}
\newcommand{\ei}{\end{itemize}}
\newcommand{\ed}{\end{document}}
\newcommand{\be}{\begin{equation}}
\newcommand{\ba}{\begin{array}}
\newcommand{\ea}{\end{array}}
\newcommand{\ee}{\end{equation}}
\newcommand{\eeq}[1]{\label{eq:#1}\end{equation}}
\newcommand{\real}{{\mathbb R}}
\newcommand{\compl}{{\mathbb C}}
\def\Id{\mbox{\boldmath $1$}}
\def\zero{\mbox{\boldmath $0$}}
\newcommand{\PP}{{\rm I\!\!\,P}}
\newcommand{\nat}{{\mathbb N}}
\newcommand{\bfpsi}{\mbox{\boldmath $\psi$}}
\newcommand{\bfchi}{\mbox{\boldmath $\chi$}}
\newcommand{\bfomega}{\mbox{\boldmath $\omega$}}
\newcommand{\bfvaromega}{\mbox{\boldmath $\varpi$}}
\newcommand{\bfOmega}{\mbox{\boldmath $\Omega$}}
\newcommand{\bfTheta}{\mbox{\boldmath $\Theta$}}
\newcommand{\bfxi}{\mbox{\boldmath $\xi$}}
\newcommand{\bfmu}{\mbox{\boldmath $\mu$}}
\newcommand{\bfx}{\mbox{\boldmath $x$}}
\newcommand{\bfy}{\mbox{\boldmath $y$}}
\newcommand{\bfPsi}{\mbox{\boldmath $\Psi$}}
\newcommand{\bfphi}{\mbox{\boldmath $\varphi$}}
\newcommand{\bfhi}{\mbox{\boldmath $\phi$}}
\newcommand{\bfPhi}{\mbox{\boldmath $\Phi$}}
\newcommand{\bfv}{{\mbox{\boldmath $v$}} }
\newcommand{\bfu}{{\mbox{\boldmath $u$}} }
\newcommand{\bfuf}{{\mbox{\footnotesize\boldmath $u$}} }
\newcommand{\bfw}{{\mbox{\boldmath $w$}} }
\newcommand{\bff}{{\mbox{\boldmath $f$}} }
\newcommand{\bfa}{{\mbox{\boldmath $a$}} }
\newcommand{\bfi}{{\mbox{\boldmath $i$}} }
\newcommand{\bfj}{{\mbox{\boldmath $j$}} }
\newcommand{\bfc}{{\mbox{\boldmath $c$}} }
\newcommand{\bfo}{{\mbox{\boldmath $o$}} }
\newcommand{\bfp}{{\mbox{\boldmath $p$}} }
\newcommand{\bfkp}{{\mbox{\footnotesize{\boldmath $k$}}} }
\newcommand{\bfka}{{\mbox{\footnotesize{\boldmath $k^*$}}} }
\newcommand{\bft}{{\mbox{\boldmath $t$}} }
\newcommand{\bfd}{{\mbox{\boldmath $d$}} }
\newcommand{\bfl}{{\mbox{\boldmath $l$}} }
\newcommand{\bfr}{{\mbox{\boldmath $r$}} }
\newcommand{\bfk}{{\mbox{\boldmath $k$}} }
\newcommand{\bfA}{{\mbox{\boldmath $A$}} }
\newcommand{\bfS}{{\mbox{\boldmath $S$}} }
\newcommand{\bfO}{{\mbox{\boldmath $O$}} }
\newcommand{\bfM}{{\mbox{\boldmath $M$}} }
\newcommand{\bfP}{{\mbox{\boldmath $P$}} }
\newcommand{\bfB}{{\mbox{\boldmath $B$}} }
\newcommand{\bfR}{{\mbox{\boldmath $R$}} }
\newcommand{\bfC}{{\mbox{\boldmath $C$}} }
\newcommand{\bfD}{{\mbox{\boldmath $D$}} }
\newcommand{\bfQ}{{\mbox{\boldmath $Q$}} }
\newcommand{\bfZ}{{\mbox{\boldmath $Z$}} }
\newcommand{\bfG}{{\mbox{\boldmath $G$}} }
\newcommand{\bfE}{{\mbox{\boldmath $E$}} }
\newcommand{\bfX}{{\mbox{\boldmath $X$}} }
\newcommand{\bfY}{{\mbox{\boldmath $Y$}} }
\newcommand{\bfH}{{\mbox{\boldmath $H$}} }
\newcommand{\bfI}{{\mbox{\boldmath $I$}} }
\newcommand{\bfJ}{{\mbox{\boldmath $J$}} }
\newcommand{\bfN}{{\mbox{\boldmath $N$}} }
\newcommand{\bfh}{{\mbox{\boldmath $h$}} }
\newcommand{\bfm}{{\mbox{\boldmath $m$}} }
\newcommand{\bfone}{{\mbox{\boldmath $1$}} }
\newcommand{\hs}{{\rm I}\!\!\,{\rm R}^3_+}
\newcommand{\cala}{{\cal A}}
\newcommand{\calb}{{\cal B}}
\newcommand{\calc}{{\cal C}}
\newcommand{\cald}{{\cal D}}
\newcommand{\cale}{{\cal E}}
\newcommand{\calf}{{\cal F}}
\newcommand{\calg}{{\cal G}}
\newcommand{\calh}{{\cal H}}
\newcommand{\cali}{{\cal I}}
\newcommand{\calj}{{\cal J}}
\newcommand{\calk}{{\cal K}}
\newcommand{\call}{{\cal L}}
\newcommand{\calm}{{\cal M}}
\newcommand{\caln}{{\cal N}}
\newcommand{\calo}{{\cal O}}
\newcommand{\calp}{{\cal P}}
\newcommand{\calq}{{\cal Q}}
\newcommand{\calr}{{\cal R}}
\newcommand{\cals}{{\cal S}}
\newcommand{\calt}{{\cal T}}
\newcommand{\calu}{{\cal U}}
\newcommand{\calv}{{\cal V}}
\newcommand{\calx}{{\cal X}}
\newcommand{\caly}{{\cal Y}}
\newcommand{\calw}{{\cal W}}
\newcommand{\calz}{{\cal Z}}
\newcommand{\bfsigma}{\mbox{\boldmath $\sigma$}}
\newcommand{\bfSigma}{\mbox{\boldmath $\Sigma$}}
\newcommand{\bftau}{\mbox{\boldmath $\tau$}}
\newcommand{\bfeta}{\mbox{\boldmath $\eta$}}
\newcommand{\bfT}{{\mbox{\boldmath $T$}} }
\newcommand{\bfV}{{\mbox{\boldmath $V$}} }
\newcommand{\bfU}{{\mbox{\boldmath $U$}} }
\newcommand{\bfW}{{\mbox{\boldmath $W$}} }
\newcommand{\bfF}{{\mbox{\boldmath $F$}} }
\newcommand{\bfK}{{\mbox{\boldmath $K$}} }
\newcommand{\bfL}{{\mbox{\boldmath $L$}} }
\newcommand{\bfb}{{\mbox{\boldmath $b$}} }
\newcommand{\bfg}{{\mbox{\boldmath $g$}} }
\newcommand{\bfn}{{\mbox{\boldmath $n$}} }
\newcommand{\bfs}{{\mbox{\boldmath $s$}} }
\newcommand{\cf}{{\it cf.} }
\newcommand{\io}{\int_\Omega}
\newcommand{\1}{\item[({\it i})]}
\newcommand{\2}{\item[({\it ii})]}
\newcommand{\3}{\item[({\it iii})]}
\newcommand{\4}{\item[({\it iv})]}
\newcommand{\5}{\item[({\it v})]}
\newcommand{\6}{\item[({\it vi})]}
\newcommand{\7}{\item[({\it vii})]}
\newcommand{\8}{\item[({\it viii})]}
\newcommand{\9}{\item[({\it xi})]}
\newcommand{\ido}{\int_{\partial\Omega}}
\newcommand{\half}{\mbox{$\frac{1}{2}$}}
\def\parallel{\|}
\def\mid{|}
\def\Bbb R{\real}
\def\hat{\widehat}
\def\tilde{\widetilde}
\def\bar{\overline}
\newcommand{\threehalves}{3\over 2}
\newcommand{\bfPi}{\mbox{\boldmath $\Pi$}}
\newcommand{\bfXi}{\mbox{\boldmath $\Xi$}}
\newcommand{\bfalpha}{\mbox{\boldmath $\alpha$}}
\newcommand{\bfbeta}{\mbox{\boldmath $\beta$}}
\newcommand{\bfgamma}{\mbox{\boldmath $\gamma$}}
\newcommand{\bfdelta}{\mbox{\boldmath $\delta$}}
\newcommand{\bfzeta}{\mbox{\boldmath $\zeta$}}
\newcommand{\bfUpsilon}{\mbox{\boldmath $\Upsilon$}}
\newcommand{\bfGamma}{\mbox{\boldmath $\Gamma$}}
\newcommand{\bfcala}{\mbox{\boldmath ${\cal A}$}}
\newcommand{\bfcalm}{\mbox{\boldmath ${\cal M}$}}
\newcommand{\bfcaln}{\mbox{\boldmath ${\cal N}$}}
\newcommand{\bfcalq}{\mbox{\boldmath ${\cal Q}$}}
\newcommand{\bfcalb}{\mbox{\boldmath ${\cal B}$}}
\newcommand{\bfcalc}{\mbox{\boldmath ${\cal C}$}}
\newcommand{\bfcali}{\mbox{\boldmath ${\cal I}$}}
\newcommand{\bfcalg}{\mbox{\boldmath ${\cal G}$}}
\newcommand{\bfcalh}{\mbox{\boldmath ${\cal H}$}}
\newcommand{\bfcalk}{\mbox{\boldmath ${\cal K}$}}
\newcommand{\bfcalt}{\mbox{\boldmath ${\cal T}$}}
\newcommand{\bfcalx}{\mbox{\boldmath ${\cal X}$}}
\newcommand{\bfcall}{\mbox{\boldmath ${\cal L}$}}
\newcommand{\bfcalf}{\mbox{\boldmath ${\cal F}$}}
\newcommand{\bfcalr}{\mbox{\boldmath ${\cal R}$}}
\newcommand{\bfcals}{\mbox{\boldmath ${\cal S}$}}
\newcommand{\bfcalw}{\mbox{\boldmath ${\cal W}$}}
\newcommand{\bfcalu}{\mbox{\boldmath ${\cal U}$}}
\newcommand{\bfcalv}{\mbox{\boldmath ${\cal V}$}}
\newcommand{\bfcalz}{\mbox{\boldmath ${\cal Z}$}}
\pagenumbering{roman}
\newcommand{\art}[6]{{\I[{\sc #1,}] {#2}, {\it #3}, {\bf #4}, {#5} {[#6]}}}
\newcommand{\ED}{\end{description}}
\newcommand{\I}{\item }
\newcommand{\ra}{\rm a}
\newcommand{\rb}{\rm b}
\newcommand{\rc}{\rm c}
\newcommand{\Hsp}{{\rm I}\!\!\,{\rm R}^n_+}
\newcommand{\Hsn}{{\rm I}\!\!\,{\rm R}^n_-}
\newcommand{\po}[1]{\mbox{$\displaystyle \frac{\mbox{$\partial #1$}}
{\mbox{$\partial x_{1}$}}$}}
\newcommand{\PO}[1]{\mbox{$\displaystyle \frac{\mbox{$\partial #1$}}
{\mbox{$\partial y_{1}$}}$}}
\newcommand{\OP}{\left(\Delta+2\lambda\PO{}\right)}
\newcommand{\op}{\left(\Delta+2\lambda\po{}\right)}
\newcommand{\ft}[1]{
\Frac{1}{(2\pi)^{n/2}}\Int{{\Bbb R}^{n}}{}e^{i{\bf x}\cdot \bfxi}
#1(\xi)d\xi}
\newcommand{\Ft}[1]{
\Frac{1}{2\pi}\Int{{\Bbb R}^{2}}{}e^{i{x}\cdot \xi}
#1(\xi)d\xi}
\newcommand{\Z}{\item[({\it a})]}
\newcommand{\B}{\item[({\it b})]}
\newcommand{\C}{\item[({\it c})]}
\newcommand{\D}{\item[({\it d})]}
\newcommand{\E}{\item[({\it e})]}
\newcommand{\G}{\item[({\it g})]}
\newcommand{\Š}{\`e}
\newcommand{\…}{\`a}
\newcommand{\•}{\`o}
\newcommand{\—}{\`u}
\newcommand{\}{\`{\i}}
\def\tag{\renewcommand{\theequation}}
\newcommand{\Footnote}{~\footnote}
\newcommand{\ie}{{\it i.e.}}
\newcommand{\dist}{\mbox{\rm dist\,}}
\newcommand{\const}{\mbox{\rm const}}
\newcommand{\trace}{\mbox{\rm trace}}
\newcommand{\Bo}{\par\hfill{$\Box$}\par\noindent}
\newcommand{\Nor}[1]{\langle{#1}\rangle_q}
\newcommand{\vs}{\vspace*{.5cm}\par\noindent}
\newcommand{\Vs}{\vspace*{.6cm}\par\noindent}
\newcommand{\Vvs}{\vspace*{.7cm}\par\noindent}
\newcommand{\VVs}{\vspace*{.8cm}\par\noindent}
\newtheorem{definition}{Definition}[section]
\newcommand{\Bd}{\begin{definition}\begin{rm}}
\newcommand{\Ed}{\end{rm}\end{definition}}
\newtheorem{remark}{Remark}[section]
\newcommand{\Br}{\begin{remark}\begin{rm}}
\newcommand{\Er}{\end{rm}\end{remark}}
\newtheorem{proposition}{Proposition}[section]
\newcommand{\Bp}{\begin{proposition}\begin{sl}}
\newcommand{\EP}[1]{\end{sl}\label{proposition:#1}\end{proposition}}
\newcommand{\propref}[1]{{\rm Proposition \ref{proposition:#1}}}
\newcommand{\Bt}{\begin{theorem}\begin{sl}}
\newcommand{\Et}{\end{sl}\end{theorem}}
\newcommand{\Bl}{\begin{lemma}\begin{sl}}
\newcommand{\El}{\end{sl}\end{lemma}}
\newtheorem{theorem}{Theorem}[section]
\newtheorem{lemma}{Lemma}[section]
\newtheorem{corollary}{Corollary}[section]
\newcommand{\eqref}[1]{{\rm (\ref{eq:#1})}}
\newcommand{\Bc}{\begin{corollary}\begin{sl}}
\newcommand{\Ec}{\end{sl}\end{corollary}}
\newcommand{\ET}[1]{\end{sl}\label{theorem:#1}\end{theorem}}
\newcommand{\EDD}[1]{\end{rm}\label{definition:#1}\end{definition}}
\newcommand{\EL}[1]{\end{sl}\label{lemma:#1}\end{lemma}}
\newcommand{\theoref}[1]{{\rm Theorem \ref{theorem:#1}}}
\newcommand{\defref}[1]{{\rm Definition \ref{definition:#1}}}
\newcommand{\ER}[1]{\end{rm}\label{remark:#1}\end{remark}}
\newcommand{\EC}[1]{\end{sl}\label{corollary:#1}\end{corollary}}
\newcommand{\remref}[1]{{\rm Remark \ref{remark:#1}}}
\newcommand{\cororef}[1]{{\rm Corollary \ref{corollary:#1}}}
\newcommand{\lemmref}[1]{{\rm Lemma \ref{lemma:#1}}}
\newcommand{\essup}[1]{{\rm ess}\,{{\displaystyle \sup_{\hspace*{-5mm}{#1}}}}}

\renewcommand{\i}{{\rm i}}

\pagenumbering{arabic}
\newcommand{\QED}{{\par\par\hfill$\square$\par}}
\renewcommand{\thefootnote}{(\arabic{footnote})}
\title{On the Self-Propulsion of a Rigid Body in a Viscous Liquid\\ by Time-Periodic Boundary Data} 
\author{ Giovanni P. Galdi 
\thanks{Department of Mechanical Engineering and Materials Science, University of Pittsburgh, PA 15261. 
}}
\date{}
\maketitle
\begin{abstract} Consider a rigid body, $\mathscr B$,  constrained to move by translational motion in an unbounded viscous liquid. The driving mechanism is a given distribution of time-periodic velocity field, $\bfv_*$,  at the interface body-liquid, of magnitude $\delta$ (in appropriate function class). The main objective is to find conditions on $\bfv_*$ ensuring that $\mathscr B$ performs a non-zero net motion, namely, $\mathscr B$ can cover any given distance in a finite time.  The approach to the problem depends on whether the averaged value of $\bfv_*$ over a period of time is (case (b)) or is not (case (a)) identically zero. In case (a) we solve the problem in a relatively straightforward way, by showing that, for small $\delta$, it reduces to the study of a suitable amd well-investigated time-dependent Stokes (linear) problem. In case (b), however, the question is much more complicated, because we show that it {\em cannot} be brought to the study of a linear problem. Therefore, in case (b), self-propulsion is a genuinely nonlinear issue that we solve directly on the nonlinear system by a contradiction argument approach. In this way, we are able to give, also in case (b), sufficient conditions for self-propulsion (for small $\delta$). Finally, we demonstrate, by means of counterexamples, that such conditions are, in general, also necessary.  
 \end{abstract}

\renewcommand{\theequation}{\arabic{section}.\arabic{equation}}
\section*{Introduction} The {\em rigorous} mathematical analysis of self-propulsion of an either rigid or shape-changing ``bulky" body in a {\em viscous} fluid is a relatively new area of investigation. Actually, a systematic   and consistent study of this problem has began about  two decades ago; see, e.g., \cite{GaSP,Garev,Court,Hishida,MS,Macha,Necasova,Raymond,MT,Staro} and the references therein. In the time-dependent case, the main, and by no means trivial, contribution of these papers   consists  in establishing   well-posedness of the corresponding initial-boundary value problem in different functional settings --weak and strong solutions,  different situations --one or more bodies, and different flow regions --bounded and unbounded.\par
However, it should also be noted that the above works seem to leave out what, we believe, is a rather fundamental question. Precisely, to establish sufficient conditions on the shape-changes or (in the case of a rigid body) boundary velocity distributions, ensuring that the body performs a non-zero net motion, that is, the body is indeed able to self-propel. This question was only considered in \cite[Section 7]{MT}, and resolved by  numerical simulation of the relevant  system of equations in the two-dimensional case. As a consequence, to the best of our knowledge, there is no rigorous analytical condition on the driving mechanism securing  self-propulsion of the body. The main objective of   
 this article is to furnish a first contribution to  such a question.
\par
Specifically, we shall consider the motion of the coupled system constituted by a rigid body, $\mathscr B$, moving in an otherwise quiescent viscous liquid, $\mathscr L$, that fills the whole space outside $\mathscr B$. The driving mechanism is a distribution of velocity, $\bfv_*$, at the boundary of $\mathscr B$, that is, at the interface body-liquid. We also suppose that on $\mathscr B$ a suitable torque  is applied  preventing  $\mathscr B$ from spinning. Thus, the motion of $\mathscr B$ will be  merely translational, and we shall denote by  $\bfgamma$ the corresponding translational velocity. 
\par 
We assume that  $\bfv_*$ is a given periodic function of time with period $T>0$ (in short: $T$-periodic function), and  look for a {\em special class} (to be defined shortly) of corresponding $T$-periodic motions of the system body-liquid. 
Under the stated assumptions, this question leads us to investigate $T$-periodic solutions to the following set of equations \cite{Garev}
\be\ba{cc}\medskip\left.\ba{ll}\medskip
\partial_t\bfv+(\bfv-\bfgamma)\cdot\nabla\bfv=\nu\Delta\bfv-\nabla p\\
\Div\bfv=0\ea\right\}\ \ \mbox{in $\Omega\times\real$}\,,\\ \medskip
\bfv(x,t)=\bfv_*(x,t)+\bfgamma(t)\,,\ \ \mbox{on\ $\partial\Omega\equiv\partial\mathscr B\times\real$}\,,\\
M\dot{\bfgamma}+\Int{\partial\Omega}{}\left[\mathbb T(\bfv,p)-(\bfv_*+\bfgamma)\otimes\bfv_*\right]\cdot\bfn=\0 \ \ \mbox{in $\real$}\,.
\ea
\eeq{SE}
Here $\bfv=\bfv(x,t)$ and  $\rho\,p$,  are velocity and pressure fields of $\mathscr L$, respectively, $\rho$ is the (constant) density of $\mathscr L$, $\nu:=\mu/\rho$ with $\mu$ shear viscosity coefficient of $\mathscr L$ and, finally, $M:=m/\rho$, where $m$ is the mass of $\mathscr B$. We observe that equation \eqref{SE}$_4$ translates into mathematical terms the requirement that $\mathscr B$ self-propels, namely, $\mathscr B$ moves without the action of {\em external} forces. \cite[Section 6]{Garev}.
\par 
The {\em main focus}  of this paper is to address the following question: 
Suppose the boundary velocity distribution $\bfv_*$ is sufficiently smooth. 
Then, 
under which conditions on $\bfv_*$ the body $\mathscr B$ performs a non-zero net motion   or, equivalently, $\mathscr B$ self-propels? \par
From the mathematical viewpoint, this question is equivalent to look for a class of $T$-periodic solutions $(\bfv,p,\bfgamma)$ to \eqref{SE} such that the translational velocity $\bfgamma$ possesses a nonzero average:
\be 
\bfxi:=\frac1T\int_0^T\bfgamma(t)\,{\rm d}t\neq\0\,.
\eeq{0.2}
\par
The physical problem that motivated our study is the self-propulsion of a ``fish," or any mechanical device where the net motion is produced by the continuous periodic movement of parts of its body.  Even though  modeling a fish as a rigid body and its moving parts as a boundary velocity distribution might look a bit coarse at first sight, it should be observed that, as is well known, the motion of a  shape-changing object in a liquid can mathematically be reduced --by an appropriate
 transformation-- to that of an object of fixed shape with a suitable distribution of velocity at its boundary; see, for example, \cite{SC} and the reference therein.  
\par
Our analysis is divided into two parts, defined by the following mutually exclusive properties of the average of $\bfv_*$ over a period:\medskip\par
\par
(a) \,\ $\bar{\bfv_*}(x):=\Frac1T\Int0T\bfv_*(x,t)\,{\rm d}t\not\equiv \0$\,;\medskip\par
(b) \,\ $\bar{\bfv_*}(x)\equiv \0$
\medskip\par
Let us denote by $\delta$ the magnitude of $\bfv_*$ (in a suitable functional class) and set 
\be 
\bfv_*=\delta\,\bfV_*.\eeq{0.2_0} 
In case (a) we  show that a sufficient condition for the validity of \eqref{0.2_0} is that, for $\delta\le \delta_0$, the $L^2$-projection, $\bfpzc G$, of $\bar{\bfV_*}$ onto a three-dimensional subspace, $\mathcal S$, of $L^2(\partial\Omega)$ is not zero.  The space $\mathcal S$  depends only on the geometric properties of
 $\mathscr B$, such as size or shape, but it is otherwise independent of its mass
and  the physical properties of $\mathscr L$. Thus, $\bfpzc G$ is the {\em thrust} for self-propulsion in case (a). Such a result, proved in the class of weak solutions to \eqref{SE} (see \theoref{Exi}), is obtained by showing that, in the limit of ``small" $\delta$ the averaged (over a period) solution must tend to the uniquely determined solution of the corresponding linear (time-independent) Stokes problem (see \lemmref{3.2}). The latter is obtained from \eqref{SE} by formally  setting  equal to 0 all time derivatives, disregarding all nonlinear terms and taking the time average of the resulting equations (see \eqref{Stsp}). We also show that the average velocity $\bfxi$ in \eqref{0.2}  in the limit of small $\delta$ has the following expression:
\be
\bfxi=\delta\, \mathbb A\cdot\bfpzc G +o(\delta)\,,
\eeq{0.3_0}
where $\mathbb A$ is an invertible matrix depending only on $\mathscr B$. 
\par 
We next consider case (b). In this situation, the boundary velocity distribution defines a purely oscillatory regime.   
From what we have already proved and, in particular, from \eqref{0.3_0} we infer that, being now $\bfpzc G=\0$,  self-propulsion must be searched at an order in $\delta$ higher than 1. Moreover, since the linearized approximation (Stokes problem) possesses in this case only the identically vanishing solution (see \lemmref{stoSP}), the solution  to the nonlinear problem that would ensure self-propulsion {\em cannot} be obtained by a perturbation argument around its linear counterpart as in case (a). In other words, in case (b) self-propulsion  becomes a {\em strictly nonlinear} phenomenon. We then use a completely different strategy that will be described next. In the first place, we split the solution $(\bfv,\bfgamma)$ to \eqref{SE} into its averaged and purely oscillatory components:
$$
\bfv=(\bfv-\bar{\bfv})+\bar{\bfv}:=\bfw+\bfu\,;\ \ \bfgamma=(\bfgamma-\bar{\bfgamma})+\bar{\bfgamma}:=\bfchi+\bfxi
$$
where the bar indicates average over a period. As a result, the problem \eqref{SE} splits into a {\em coupled} nonlinear problem of ``elliptic-parabolic" type in the unknowns $(\bfu,\bfxi)$ for the ``elliptic" part, and $(\bfw,\bfchi)$ for the ``parabolic" part (see \eqref{6.10}--\eqref{6.13}). If $\delta\le\delta_0$, we then show that the problem, in this new form, has a unique ``strong" solution (see \theoref{6.1}). This result is obtained by a successive   approximation scheme around an appropriate {\em nonlinear problem} (see \eqref{6.10_0}--\eqref{6.11_0}). The fundamental aspect of this finding is that, thanks to the self-propelling condition \eqref{SE}$_4$, we can show that the averaged component $\bfu$ belongs to a certain function class where, in general, {\em classical} boundary-value problems associated to Stokes and Navier-Stokes equations (namely, {\em without} the self-propelling condition) do {\em not} have a solution {\em if} $\bfxi=\0$ \cite{KoSo,GaNWP}. This space is the homogeneous Sobolev space $D^{1,\frac32}(\Omega)$ of locally integrable functions with spatial derivatives belonging to the Lebesgue space $L^{\frac32}(\Omega)$. In the next step, we solve the ``parabolic" problem for $(\bfw,\bfchi)$, considering  $\bfu$ as assigned in its function class, thus implicitly getting $(\bfw=\bfw(\bfu,\bfv_*),\bfchi=\bfchi(\bfu,\bfv_*))$, and replace the latter into the ``elliptic'' problem. This can be done, provided $\delta$ is sufficiently ``small." The final step is then a contradiction argument on the nonlinear problem satisfied by $\bfu$. Precisely, with the help of the previous step, and assuming  ad absurdum $\bfxi=\0$, we show that $\bfu$ must solve the following problem
\be  
\ba{cc}\medskip\left.\ba{ll}\medskip
\nu\Delta\bfu-\nabla p=\bfscr F(\bfu,\bfv_*)\\ \Div\bfu=0\ea\right\}\ \ \mbox{in $\Omega$}\,,
\\ 
\bfu=\0\ \ \mbox{on $\partial\Omega$\,,}
\ea
\eeq{0.4}
with $\bfscr F$  a sufficiently smooth, nonlinear function of its arguments. However,  $\bfu\in D^{1,\frac32}(\Omega)$ and so, as previously noticed, the problem \eqref{0.4} does not, in general,  admit solutions in that class. Actually, a solution in $D^{1,\frac32}(\Omega)$ may exists  {\em if and only if} $\bfscr F$ satisfies a nonlocal compatibility condition \cite[Section V.5]{Gab}. Our self-propulsion problem reduces then to find requirements on $\bfv_*$ that {\em violate} this condition, thus implying $\bfxi\neq\0$.
We then prove that, for $\delta$ sufficiently ``small", there is 
a vector, $\bfG$,  in $\real^3$ depending only on the ``shape" of $\mathcal B$, its mass  and $\bfV_*$ (see \eqref{G}) such that if $\bfG\neq\0$, then
the  compatibility condition mentioned above is violated, and hence $\bfxi\neq\0$ (see \theoref{6.2}). Thus, $\bfG$ is the {\em thrust}, in case (b). Furthermore,  we prove (see \theoref{7.1_0})
\be
\bfxi=\delta^2\,\mathbb A\cdot\bfG+O(\delta^{\frac{11}4})\,.
\eeq{0.7}
This formula shows, as expected, that if $\bfG\neq\0$ self-propulsion does occur occur at an order in $\delta$ higher than 1.
The natural question to ask is then  whether, if $\bfG=\0$, could self-propulsion  occur at an order in $\delta$  higher than 2. We show that, in general, the answer is {\em negative} (see Section 9). In fact, we give an example  showing that given a body $\mathscr B$ of {\em any shape and mass}, and an arbitrary period $T>0$, for any any $\delta>0$ there is always a $T$-periodic boundary velocity $\bfv_*$ such that, if $\bfG=\0$,  the averaged velocity field is {\em identically vanishing}, thus implying that  $\mathscr B$ can only ``oscillate," with  zero net motion.
\par
The outline of the paper is as follows. After some preliminary results collected in Section 1, in Section 2 we give a weak formulation of problem \eqref{SE} and prove in \theoref{Exi}  existence of corresponding weak solutions under suitable assumptions on $\bfv_*$. In the following Section 3, we show in \theoref{NZA} sufficient conditions for self-propulsion when $\bar{\bfv_*}\not\equiv\0$ in the class of weak solutions. The next two sections are dedicated to preparatory results  necessary for the investigation of the case $\bar{\bfv_*}\equiv 0$. Specifically, in Section 4, we prove existence and uniqueness of a steady-state {\em nonlinear} problem in a function class contained in the homogeneous Sobolev space $D^{1,\frac32}(\Omega)$ mentioned earlier on; see \lemmref{Fm}. This lemma plays a {\em key role in all our subsequent analysis}, and   the self-propelled condition becomes a necessary and sufficient requirement for its general validity. The main result of the following Section 5 is the proof of existence and uniqueness in the maximal regularity class of solutions to the linearized version of \eqref{SE}, in the case when the data have zero average; see \lemmref{1.6_0}. With the help of these findings we show, in Section 6, existence and uniqueness of  $T$-periodic solutions to \eqref{SE} in a rather regular function class, provided $\delta$ is sufficiently restricted; see \theoref{6.1}. By using the  contradiction argument mentioned earlier on, we provide in Section 7 sufficient conditions for self-propulsion when $\bar{\bfv_*}\equiv\0$ (see \theoref{6.2}), whereas in Section 8 we furnish the expression  of the velocity of propulsion; see \theoref{7.1_0}. In the final Section 9, we prove by means of counterexamples that the self-propelling condition determined in Section 7 is also necessary.

\setcounter{section}{0}
\section{Preliminary Results} In this section we shall recall and/or  introduce the main notation, and collect some basic results that will be frequently used later on in the paper. \par By $\Omega$ we indicate a domain of $\real^3$,  complement of the closure of a  bounded domain $\Omega_0$ ($\equiv$ the body $\mathscr B$). We assume $\Omega$ of class $C^2$. Moreover, with the origin in the interior of $\Omega_0$, we set $\Omega_R:=\Omega\cap\{|x|<R\}$ and $\Omega^R:=\Omega\cap\{|x|>R\}$, for $R>R_*:={\rm diam}\,\Omega_0$.  As customary, for a domain $A\subseteq\real^3$, $L^q=L^q(A)$ is the Lebesgue space with norm $\|\cdot\|_{q,A}$, and  $W^{m,q}=W^{m,q}(A)$ denotes Sobolev space, $m\in\nat$, $q\in[1,\infty]$, with norm $\|\cdot\|_{m,q,A}$. Corresponding trace norms at $\partial A$ are denoted by $\|\cdot\|_{m-\frac1q,q(\partial A)}$. Furthermore, $D^{m,q}=D^{m,q}(A)$ are homogeneous Sobolev spaces with semi-norm $|u|_{m,q,A}:=\sum_{|l|=m}\|D^lu\|_q$, whereas $D_0^{1,q}=D_0^{1,q}(A)$ is the completion of $C_0^\infty(A)$ in the norm $|\cdot|_{1,q,A}$. 
\footnote{A detailed analysis of homogeneous Sobolev spaces including  their main properties can be found in \cite[Section II.6]{Gab}.} In all the above  norms, the subscript $"A"$ will be omitted, unless confusion arises. An important embedding property of the spaces $D^{1,q}$ is recalled in the following lemma,  for whose proof we refer to \cite[Theorem II.6.1(i) and Theorem II.7.3]{Gab}.
\Bl Let  $u\in D^{1,q}(\Omega)\cap L^r(\Omega)$, for some  $q\in [1,3)$, $r\in [1,\infty)$. Then $u\in L^{\frac{3q}{3-q}}(\Omega)$ and there is $c=c(\Omega,q)$ such that 
$$
\|u\|_{\frac{3q}{3-q}}\le c\,|u|_{1,q}\,.
$$
Suppose, in addition, $u\in D^{2,2}(\Omega)$. Then 
$$
\|u\|_{s}+|u|_{1,\sigma}\le c\,\left(|u|_{1,q}+|u|_{2,2}\right)\,,\ \ \mbox{for all $s\in [\frac{3q}{3-q},\infty]$\, and $\sigma\in[q,6]$}\,.
$$
\EL{HSS}
\par
Let
$$
\calc\!=\!\mathcal C(\Omega)\!:=\!\big\{\bfphi\in C_0^\infty(\bar{\Omega})\!: \Div\bfphi=0 \ \mbox{in $\Omega$};
\bfphi(x)=\bfxi_{\mbox{\footnotesize $\bfphi$}},\, \mbox{some}\ \bfxi_{\sub{\bfphi}}\in\real^3, \mbox{in a neighborhood of $\partial\Omega$}\big\}\,,
$$
and  define
$$
\mathcal H=\calh(\Omega)\equiv \,\big\{\mbox{completion of $\calc(\Omega)$ in the norm $\|\mathbb D(\cdot)\|_2$}\big\}\,.
$$
The essential properties of the space $\calh$ are collected in the  next lemma, whose proof is given in \cite[Lemmas 9--11]{Garev}.
\Bl  $\calh$ is a Hilbert space endowed with the scalar product
\be
\big[\bfu,\bfw\big]:=\int_\Omega\mathbb D(\bfu):\mathbb D(\bfw)\,,\ \ \bfu,\bfw\in\calh\,,
\eeq{sp}
and the following characterizations hold
\be
\calh=\calh(\Omega):=\big\{\bfu\in W^{1,2}_{\rm loc}(\bar{\Omega}): \bfu\in L^6(\Omega),\,\mathbb D(\bfu)\in L^2(\Omega)\,;\ \Div\bfu=0\ \mbox{in\, $\Omega$}\,;\ \bfu(y)=\bfxi_\sub{\bfu}\,,\ y\in\partial\Omega\big\}\,.
\eeq{0_0}
Moreover, we have
\be
\|\nabla\bfu\|_2\le\sqrt{2}\|\mathbb D(\bfu)\|_2\le 2\|\nabla\bfu\|_2\,,
\eeq{WH1}
and
\be
\|\bfu\|_6\le c_0\,\|\mathbb D(\bfu)\|_2\,,\ \ \bfu\in\calh\,,
\eeq{SoB}
for some constant $c_0>0$.  Finally, 
there is another positive constant $c_1$ such that
\be
|\bfxi_\sub{\bfu}|\le c_1\,\|\mathbb D(\bfu)\|_2\,.
\eeq{WH2}
\EL{1}
We shall need also the following ``local" version of the above spaces:
$$
\mathcal C(\Omega_R):=\!\big\{\bfphi\in C_0^\infty(\bar{\Omega_R})\!: \Div\bfphi=0 \ \mbox{in $\Omega_R$};
\bfphi(x)=\bfxi_{\mbox{\footnotesize $\bfphi$}},\ \mbox{around $\partial\Omega$},\ \bfphi=\0 \ \mbox{around $\partial B_R$}\big\}\,,
$$
and 
$$
\calh(\Omega_R):=\big\{\bfu\in W^{1,2}(\Omega_R): \ \Div\bfu=0\ \mbox{in\, $\Omega_R$}\,;\ \bfu(y)=\bfxi_\sub{\bfu}\,,\ y\in\partial\Omega\,,\ \bfu=\0\ \mbox{at $\partial B_R$}\big\}\,.
$$
For each $R>R_*$, $\calh(\Omega_R)$ becomes a Hilbert space when endowed with  the scalar product defined in \eqref{sp}, by setting $\Omega\equiv\Omega_R$.
We denote by $\calh^{-1}(\Omega_R)$ its dual space. Also, $\calc(\Omega_R)$ is dense in $\calh(\Omega_R)$.\smallskip
\par
For $q\in (1,\infty)$ we introduce the following Banach spaces
$$
\ba{ll}\medskip
\sfD^{2,q}:=D^{2,2}(\Omega)\cap D^{1,q}(\Omega)\cap \calh(\Omega)\,;\ \ 

{\sf D}^{1,q}
:=D^{1,2}(\Omega)\cap 
L^{q}(\Omega)\,,\\
\sfL^q:=\{\bff=\Div\bfscr F\in L^{2}(\Omega):\ \bfscr F\in L^{2}(\Omega)\cap L^{q}(\Omega)\},,
\\
\ea
$$
endowed with the  norms
$$\ba{ll}\medskip
\|\bfu\|_{{\sfD}^{2,q}}:=|\bfu|_{2,2}+|\bfu|_{1,q}+\|\bfu\|_6\,,\ \ 
 \|u\|_{{\sf D}^{1,q}}:=|u|_{1,2}+\|u\|_{q}\,\\
\ \ \|\bff\|_{{\sf L}^{q}}:=\|\bff\|_{2}+\|\bfscr F\|_{2}+\|\bfscr F\|_{q}\,.
\ea
$$
From \lemmref{1},  \lemmref{HSS}, and  \cite[Theorem II.9.1]{Gab} we deduce the following  result.    
\Bl Let $q\in(1,3)$, and set ${q_1}:=\min\{\frac{3q}{3-q},6\}$, ${q_2}:=\min\{q,2\}$. Then, the  following continuous embedding holds
$$
{\sf D}^{2,q}\subset\left\{\ba{ll}\medskip
L^s(\Omega)\,,\ \ s\in [q_1,\infty]\\
D^{1,\sigma}(\Omega)\,,\ \ \sigma\in [q_2,6]\,.
\ea\right.
$$
Moreover, if $\bfu\in \sfD^{2,q}$, then
$$
\lim_{|x|\to \infty}\bfu(x)=\0.
$$
\EL{Emb0}

If $A\subset\real^3$, a function $u:A\times \real\mapsto \real^3$ is 
{\em $T$-periodic}, $T>0$, if $u(\cdot,t+T)=u(\cdot\,t)$, for a.a. $t\in \real$,
 and we set
$
{\bar u}:=\frac{1}{T}\int_{0}^{T}u(t){\rm d}t\,.
$
Let $B$ be a function space endowed with seminorm $\|\cdot\|_B$, $r=[1,\infty]$, and $T>0$. Then, $L^r(0,T;B)$ is the class of functions
$u:(0,T)\rightarrow B$ such that 
$$
\|u\|_{L^r(B)}\equiv\left\{\ba{ll}\smallskip\big( \Int{0}{T}\|u(t)\|_B^r \big)^{\frac 1r}<\infty, \ \ \mbox{if 
$r\in [1,\infty)\,;$}\\   
\essup{t\in[0,T]}\,\|u(t)\|_B <\infty, \ \ \mbox{if $r=\infty.$}
\ea\right.
$$
Likewise, we put
$$
W^{m,r}(0,T;B)=\Big\{u\in L^{r}(0,T;B): \textcolor{black}{\partial_t^ku\in L^{r}(0,T;B), \, k=1,\ldots,m}\Big\}\,.
$$
Unless confusion arises, we shall simply write $L^r(B)$ for $L^r(0,T;B)$, etc. Moreover, for $q\in(1,\infty)$, we introduce the following Banach spaces
$$\ba{ll}\medskip
L^{q}_\sharp(0,T)=\{\chi\in L^{q}(0,T), \ \mbox{$\chi$ is $T$-periodic with }\ \bar{\chi}=0\}\\\medskip
W^{1,q}_\sharp(0,T)=\{\chi\in L^{q}_\sharp(0,T), \ \dot\chi\in L^q(0,T)\}\\\medskip
\mathcal L_\sharp^{q}:=\{\bfu\in L^{q}(L^q); \ \mbox{$\bfu$ is $T$-periodic, with $\bar{\bfu}=\0$}\}
\\ \medskip
\mathcal W_\sharp^{2}:=\{\bfu\in W^{1,2}(L^2)\cap L^2(W^{2,2}\cap\calh); \ \mbox{$\bfu$ is $T$-periodic, with $\bar{\bfu}=\0$}\}\\
\hat{\mathcal W}_\sharp^{q}:=\{\bfu\in W^{1,q}(L^q)\cap L^q(W^{2,q}); \ \mbox{$\bfu$ is $T$-periodic, with $\bar{\bfu}=\0$}\}
\ea
$$
endowed with  natural norms, and define  
$$
\mathcal L^{2,q}_\sharp:= \mathcal L_\sharp^{2}\cap L^q(L^q)\,;\ \  \mathcal W_\sharp^{2,q}:=\mathcal W_\sharp^{2}\cap \hat{\mathcal W}_\sharp^{q}\,;\ \ \hat{\mathcal W}_\sharp^{2,q}:=\hat{\mathcal W}_\sharp^{2}\cap \hat{\mathcal W}_\sharp^{q}
$$
with associated norms
$$\ba{ll}\medskip
\|\bfu\|_{\mathcal L^{2,q}_\sharp}:=\|\bfu\|_{L^2(L^2)}+\|u\|_{L^q(L^q)}\\
\|\bfu\|_{\hat{\mathcal W}_\sharp^{2,q}}=\|\bfu\|_{\mathcal W_\sharp^{2,q}}:=\|\bfu\|_{W^{1,2}(L^2)\cap W^{1,q}(L^q)}+\|\bfu\|_{L^2(W^{2,2})\cap L^q(W^{2,q})}\,.
\ea
$$
Finally, we define
$$
\calp^{1,q}:= L^2(D^{1,2})\cap L^q(D^{1,q})\,.
$$
We recall some embedding properties of the spaces $\calw_\sharp^2$ which are a particular case of
 \cite[Theorem 2.1]{Mallo}.
\Bl The following continuous embedding holds, for all $r,s\in[q,\infty]$:
$$
\hat{\calw}_\sharp^q\subset\left\{\ba{ll}\medskip L^r(L^s)\,,\ \ \frac3s+\frac2r>\frac5q-2\,,\ \,\,,\\
L^r(D^{1,s})\,,\ \ \frac3s+\frac2r>\frac5q-1\,.\ea\right.
$$
\EL{Mallo}\par
We conclude this section with a couple of further {\bf notational remarks}. The first one regards the standard Landau notation. Precisely,  by $O(\delta^\alpha)$, [respectively, $o(\delta^\alpha)$],  $\alpha\ge0$, we indicate a generic function $f$ (say) depending on $\delta$ and such that $|f(\delta)|\le c\,\delta^\alpha$,  $\delta\le\delta_0$, for some positive constants $c,\delta_0$ [respectively, $\lim_{\delta\to 0}f(\delta)\delta^{-\alpha}=0$]. Finally, by  $c$, $c_0$, $c_1$, etc.,  we denote positive constants, whose particular value is unessential to the context. When we wish to emphasize
the dependence of $c$ on some parameter $\zeta$, we shall write  $c(\zeta)$ or $c_\zeta$.
\setcounter{equation}{0}
\section{Weak Formulation and Existence of of Weak Solutions}
In this section we shall prove existence to problem \eqref{SE} in the very general class of weak solutions, under the assumption  that the boundary velocity distribution has zero net flux through the boundary of $\mathscr B$ (see \eqref{flux}) and it is sufficiently 	``small" in appropriate norm. Probably, both conditions could be weakened, but this would not be relevant to our main objective of finding sufficient conditions for self-propulsion that will be discussed in the next section.
\par
To reach the goal above, we begin to put \eqref{SE} in a weak form. For  $A\equiv\Omega_R,\Omega$, we denote by $\calc_{\sharp}(A):=\calc_{\sharp}(A\times [0,T])$  the space of vector functions 
obtained by restriction to $[0,T]$ of functions $\bfpsi \in C^{1}(A\times \mathbb{R})$, satisfying: 
\begin{enumerate}
\item $\Div \bfpsi = 0$ in $A \times \mathbb{R}$\,; 
\item There exists $\bfgamma_\sub{\bfpsi} \in C^{1}(\mathbb{R})$ such that $\bfpsi (x,t) = \bfgamma_\sub{\bfpsi}(t) $ for $x$ in a neighborhood 
of $\partial\Omega$, and $t \in \mathbb{R}$\,;
\item For each $\bfpsi,$ there exists $\rho=\rho(\bfpsi) > R_*$ such that $\bfpsi (x,t) = 0$ for $|x| \geq \rho$ and $t \in \mathbb{R}$\,, with $\rho<R$ if $A\equiv\Omega_R$\,;
\item $\bfpsi$ is $T$-periodic.
\end{enumerate}
Multiplying formally \eqref{SE}$_1$ by the test function $\bfpsi\in\calc_\sharp(\Omega)$,  integrating by parts over $\Omega\times [0,T]$, and taking into account $T$-periodicity, we find
$$
- \Int{0}{T}\Int{{\Omega}}{}\partial_t{ \bfpsi} \cdot \bfv 
=
\Int{0}{T}\Int{\partial\Omega}{}\bfgamma_\sub{\bfpsi}\cdot \left[\mathbb T(\bfv,p)\cdot \bfn - (\bfv_{*}+\bfgamma)\bfv_{*}\cdot \bfn\right]
+ \Int{0}{T}\Int{{\Omega}}{}[(\bfv-\bfgamma) \cdot \nabla \bfpsi \cdot \bfv
-2\,\nu\, \mathbb D(\bfpsi):\mathbb D(\bfv)]\,.
$$
Then, imposing in the latter the self-propelling condition \eqref{SE}$_4$, we get
\be
- \Int{0}{T}\big(\Int{{\Omega}}{}\partial_t{ \bfpsi} \cdot \bfv +
M\dot{\bfgamma}_\sub{\bfpsi}\cdot\bfgamma\big)
= \Int{0}{T}\Int{{\Omega}}{}[(\bfv-\bfgamma) \cdot \nabla \bfpsi \cdot \bfv
-
2  \,\nu\,\mathbb D(\bfpsi):\mathbb D(\bfv)]
\eeq{ana1}
Following \cite{GaSi}, we give the following definition.
\Bd
The pair
$\{\bfv,\bfgamma\}$ is a \emph{$T$-periodic weak solution} to \eqref{SE} if the following conditions hold:
\begin{itemize}
\item[(i)] $\bfv$ and $\bfgamma$ are both $T$-periodic with $\bfv \in L^{2}(0,T;D^{1,2}({\Omega}))$, $\bfgamma\in L^2(0,T)$\,; 
\item[(ii)] $\Div\bfv(\cdot,t)=0$ in $\Omega$, for a.a. $t\in[0,T]$\,;  
\item[(iii)] $\bfv = \bfv_{*}+\gamma$ at  $\partial\Omega \times (0,T)$ (in the trace sense)\,;
\item[(iii)] $\{\bfv,\bfgamma\}$ verifies \eqref{ana1}, for all $\bfpsi \in {\cal C}_\sharp(\Omega).$
\end{itemize}
\EDD{1}
\par
We need a preparatory result concerning suitable extension of the boundary data. 
\Bl Let $\bfV_*\in W^{1,2}(W^{\frac12,2}(\partial\Omega))$ be $T$-periodic with 
\be
\int_{\partial\Omega}\bfV_*(x,t)\cdot\bfn=0\,,\ \ \mbox{for all $t\in [0,T]$.}
\eeq{flux}
Then, there exists a $T$-periodic field $\bfV\in W^{1,2}(W^{1,2}(\Omega))$ such that
\begin{itemize}
\item[{\rm (a)}] $\bfV(x,t)=\bfV_*(x,t)$\,,\ \ \mbox{for all $(x,t)\in \partial\Omega\times [0,T]$\,;}
\item[{\rm (b)}] $\Div\bfV=0$\ \ \mbox{in $\Omega\times [0,T]$}\,;
\item[{\rm (c)}] $\bfV(x,t)=0$\,,\ \ \mbox{for all $(x,t)\in \Omega^{2R_*}\times[0,T]$\,.}
\item[{\rm (d)}] $\|\bfV(t)\|_{1,2}\le c\, \|\bfV_*(t)\|_{\frac12,2(\partial\Omega))}$\,,\ \mbox{for a.a. $t\in[0,T]$}\,; 
\item[{\rm (e)}] $\|\bfV\|_{W^{1,2}(W^{1,2})}\le c\, \|\bfV_*\|_{W^{1,2}(W^{\frac12,2}(\partial\Omega))}$\,.  
\end{itemize}
\EL{xten}
{\em Proof.} For each $t\in [0,T]$, let $\bfU=\bfU(x,t)$ be the solution to the following Stokes problem
$$\ba{cc}\medskip\left.\ba{ll}\medskip
\Delta\bfU=\nabla P\\
\Div\bfU=0
\ea\right\}\ \ \mbox{in $\Omega$}\\
\bfU=\bfV_*(t)\ \ \mbox{at $\partial\Omega$\,,}\ \ \bfU=\0\ \ \mbox{at $\partial B_{3R_*}$\,.}
\ea
$$
From classical results \cite[Theorem IV.6.1]{Gab}, we know that there exists a unique solution
$(\bfU(t),P(t))\in W^{1,2}(\Omega)\times L^2(\Omega)$ satisfying, in addition, 
\be|\bfU (t)|_{1,2}\le c\,\|\bfV_*(t)\|_{\frac12,2(\partial\Omega)}\,.
\eeq{STOE} 
Clearly, $\bfU(t)$ is $T$-periodic. Moreover, $\bfU$ is time-differentiable and since $\partial_t\bfU$ satisfies the same problem as $\bfU$ with $\bfV_*$ replaced by $\partial_t\bfV_*$, we recover that $\partial_t\bfU(t)\in W^{1,2}(\Omega)$ and
\be|\partial_t\bfU (t)|_{1,2}\le c\,\|\partial_t\bfV_*(t)\|_{\frac12,2(\partial\Omega)}\,.
\eeq{STOE1}
Let $\varphi=\varphi(|x|)$ be a smooth function such that $\varphi=1$, for $|x|\le R_*$, and $\varphi=0$, for $|x|\ge 2R_*$, and consider the problem
\be\ba{ll}\medskip
\Div\bfw(t)=-\nabla\varphi\cdot\bfU:=f(t)\ \mbox{in $\Omega$}\,;\,\ \ \bfw,\partial_t\bfw\in W^{1,2}_0(\Omega_{2R_*})\,;\\ 
\|\bfw(t)\|_{1,2}\le c\,\|f(t)\|_{2}\,,\ \ \|\partial_t\bfw(t)\|_{1,2}\le c\,\|\partial_tf(t)\|_{2}\,.
\ea
\eeq{STOE2}
Since, by \eqref{flux}, $\int_{\Omega}f(t)=0$ for all $t\in[0,T]$, in view of \cite[Exercise III.3.6]{Gab}, we may assert that such function $\bfw$ exists and is also $T$-periodic. Thus, if we set $\bfV:=\varphi\,\bfU+\bfw$, extend $\bfw$ to 0 outside $\Omega_{2R_{*}}$, and employ the regularity assumption on $\bfV_*$ along with  \eqref{STOE}--\eqref{STOE2}, we recognize at once that $\bfV$ satisfies all the stated properties.
\par\hfill$\square$\medskip\par
The following lemma shows, in particular, a further property of a weak solution that also  furnishes a more specific way in which it satisfies the periodicity property. 
\Bl Let $(\bfv,\bfgamma)\in L^2(L^2(\Omega))\times L^2(0,T)$ satisfy \eqref{ana1} for all $\bfpsi\in \calc_\sharp(\Omega_R)$. Moreover, suppose that $\bfv-\bfgamma=\bfv_*$ at $\partial\Omega\times[0,T]$, with $\bfv_*\in W^{1,2}(W^{\frac12,2}(\partial\Omega))$, and that $\bfv-\bfcalv\in L^2(\calh)$, where $\bfcalv$ is the extension of $\bfv_*$ given in \lemmref{xten}. Then, 
$\partial_t\bfv\in L^1(0,T;\calh^{-1}(\Omega_R))$, and  there is a constant $c=c(R,\nu)$ such that
$$
\|\partial_t\bfv\|_{L^1(\calh^{-1}(\Omega_R))}\le c\,\left(\|\bfv\|_{L^2(\calh)}^2+
\|\bfv_*\|_{W^{1,2}(W^{\frac12,2}(\partial\Omega))}^2\right)\,.
$$
So, in particular,
$$
{\bfv}\in C([0,T];\calh^{-1}(\Omega_R))\,.
$$
\EL{Hprime}
{\em Proof.}
Let us choose in \eqref{ana1} $\bfpsi=\chi\bfphi$, where $\bfphi\in\calc(\Omega_R)$ and $\chi\in C_0^\infty(0,T)$.
We then obtain
$$
 \Int{0}{T}(\bfv,\bfphi)\dot{\chi}
= -\Int{0}{T}G_\sub{\bfphi}(t)\chi(t)\,,\ \ \mbox{for all $\chi\in C_0^\infty(0,T)$.}
$$
where
$$ G_\sub{\bfphi}(t):=
\Int{{\Omega}}{}[(\bfv-\bfgamma) \cdot \nabla \bfphi \cdot \bfv
-
2  \,\nu\,\mathbb D(\bfphi):\mathbb D(\bfv)]
$$
We now write $\bfv=\bfu+\bfcalv$. Thus, clearly, for a.a. $t\in[0,T]$, $\bfu\in\calh(\Omega)$ and $\bfu=\bfgamma$ at $\partial\Omega$.
Observing that, by \eqref{SoB} and \eqref{WH2}, 
$$
\|\bfu\|_{4,\Omega_R}\le C(R)\,\|\bfu\|_{6}\le c(R)\,\|\mathbb D(\bfu)\|_{2}\,;\ \ |\bfgamma|\le c\,\|\mathbb D(\bfu)\|_{2}
$$
and that, by \lemmref{xten} and embedding theorems,
$$
\|\bfcalv\|_4\le c\,\|\bfcalv\|_{1,2}\le c\,\|\bfv_*\|_{\frac12,2(\partial\Omega)}\,,
$$
with the help of Schwarz inequality and \eqref{WH1}, we obtain for a.a. $t\in[0,T]$
$$
|G_\sub{\bfphi}(t)|\le c(R,\nu)\,\left(\|\mathbb D(\bfu)\|_{2}^2+\|\mathbb D(\bfu)\|_{2}+\|\bfv_*\|_{\frac12,2(\partial\Omega)}^2+\|\bfv_*\|_{\frac12,2(\partial\Omega)}
\right)\|\mathbb D(\bfphi)\|_{2,\Omega_R} 
\,.
$$
Thus, since $\bfu\in L^2(\calh)$, we find $G_\sub{\bfphi}(t)=\langle\bfg(t),\bfphi\rangle$ with  $\bfg\in L^1(0,T; \calh^{-1}(\Omega_R))$ and $\langle\cdot,\cdot\rangle$ duality pairing, and, moreover,
$$
\ode{}{t}(\bfv,\bfphi)= \langle \bfg ,\bfphi\rangle\,,
$$
in the sense of distributions on $[0,T]$. The desired property is then proved.
\par\hfill$\square$\ms\par
We are now in a position to prove the following existence result.
\Bt Let $\bfv_*=\delta\,\bfV_*$, with $\bfV_*$ as in \lemmref{xten}. Then, there is $\delta_0>0$ such that for each $\delta\in (0,\delta_0)$, the problem \eqref{SE} has at least one weak solution. Moreover, 
the following estimate holds
\be
\|\bfgamma\|_{L^2(0.T)}+\|\nabla\bfv\|_{L^2(L^2)}\le c(\Omega)\,\delta\,.
\eeq{stima}
\ET{Exi}
{\em Proof.} The method we use is close to the one employed in \cite[Section 3]{GS1},  \cite[Section 3]{GaSi} in a similar context. We shall therefore restrict ourselves to provide the main steps, referring the reader to the cited work for details. The basic idea is to couple the classical Galerkin method with the ``invading domains" procedure. Thus, let $\{\Omega_{R_k}\}$ be a sequence  with $R_k\to\infty$ as $k\to\infty$, and $\Omega_{R_1}\supset {\rm supp}\,(\bfV)$, where $\bfV$ is the extension field constructed in \lemmref{xten}. For each fixed $k$,  we look for a weak solution  to \eqref{SE} in $\Omega_{R_k}$, namely, a pair $(\bfv^{(k)},\bfgamma^{(k)})$ such that
\be\ba{ll}\medskip
- \Int{0}{T}\Int{{\Omega_{R_k}}}{}\left(\partial_t{ \bfpsi} \cdot \bfv^{(k)} +
M\dot{\bfgamma}_\sub{\bfpsi}\cdot\bfgamma^{(k)}\right)
= \Int{0}{T}\Int{{\Omega}}{}[(\bfv^{(k)}-\bfgamma^{(k)}) \cdot \nabla \bfpsi \cdot \bfv^{(k)}
-
2  \,\nu\,\mathbb D(\bfpsi):\mathbb D(\bfv^{(k)})]\,,\\ 
\hspace*{9cm}\mbox{for all $\bfpsi\in C_\sharp(\Omega_{R_k})$}\,,
\ea
\eeq{ana1k}
and satisfying property (i) and (ii) of \defref{1} with $\Omega\equiv\Omega_{R_k}$, property (iii) and, moreover, $\bfv^{(k)}=0$ at $\partial B_{R_k}\times [0,T]$. 
Such a solution 
is then searched as suitable limit of a sequence of functions $\bfv_m^{(k)}:=\bfu_m^{(k)}+\bfcalv$, $\bfcalv:=\delta\bfV$,  solving the following sequence of ``approximating" problems, $\ell=1,\ldots,m$, $m\in\nat$: 
\be\ba{rl}\medskip
(\partial_t\bfu_m^{(k)},\bfphi_\ell)+M\dot{\bfgamma}_m^{(k)}\cdot\bfgamma_\sub{\bfphi_\ell}+2\,\nu\,[\bfu_m^{(k)},\bfphi_\ell]-\big((\bfu_m^{(k)}+&\!\!\!\!\bfcalv-\bfgamma_m^{(k)})\cdot\nabla\bfphi_\ell,(\bfu_m^{(k)}+\bfcalv)\big)\\
&\!\!\!+2\,\nu\,[\bfcalv,\bfphi_\ell]+(\partial_t\bfcalv,\bfphi_\ell)=0\,,
\ea
\eeq{x1}
where $\{\bfphi_\ell\}\subset \calc(\Omega_{R_k})$ is base of $\calh(\Omega_{R_k})$ normalized by the condition
\be
(\bfphi_\ell,\bfphi_{\ell'})+M\,\bfgamma_\sub{\bfphi_\ell}\cdot\bfgamma_\sub{\bfphi_{\ell'}}=\delta_{\ell,\ell'}\,,
\eeq{norm}
whereas
$$
\bfu_m^{(k)}:=\sum_{i=1}^mc_{mi}^{(k)}\bfphi_i\,,\ \ \bfgamma_m^{(k)}:=\sum_{i=1}^mc_{mi}^{(k)}\bfgamma_\sub{\bfphi_i}\,,
$$
where $c_{mi}^{(k)}(t)$ are functions of time only, requested to solve  \eqref{x1}. If we multiply both sides of \eqref{x1} by $c_{m\ell}^{(k)}(t)$,  sum over $\ell$ from 1 to $m$, and integrate by parts as necessary we get the following relation where, in order to alleviate the notation, we suppress the superscript ${}^{(k)}$:
\be\ba{rl}\medskip
\half\ode{}t(\|\bfu_m\|^2_2+M\,|\bfgamma_m|^2)+2\nu\|\mathbb D(\bfu_m)\|_2^2=-((\bfu_m&\!\!\!-\bfgamma_m)\cdot\nabla\bfcalv,\bfu_m)\\
&\!\!\!+(\bfcalv\cdot\nabla\bfu_m,\bfcalv)-2\nu[\bfcalv,\bfu_m]+(\partial_t\bfcalv,\bfu_m)\,.
\ea
\eeq{x2}
Denote by $\calr_m$ the right-hand side of \eqref{x2}. Since $\Omega_{R_1}$ strictly contains the support of $\bfcalv$, with the help of \lemmref{1} and classical embedding theorems we show
$$
|\calr_m|\le c\,\big[\|\bfcalv\|_{L^\infty(W^{1,2})}\|\mathbb D(\bfu_m)\|_2^2+\big(\|\bfcalv\|_{L^\infty(W^{1,2})}^2+\|\bfcalv\|_{L^\infty(W^{1,2})}\big)\|\mathbb D(\bfu_m)\|_2\big]\,,
$$
with the constant $c$ independent of $m$.
We now employ in this relation Cauchy-Schwarz inequality along with \lemmref{xten}(e), take $\delta$ sufficiently small and replace the resulting expression back into \eqref{x2} to infer 
\be
\ode{}t(\|\bfu_m\|^2_2+M\,|\bfgamma_m|^2)+\nu\,\|\mathbb D(\bfu_m)\|_2^2\le c\,\|\bfcalv\|_{L^\infty(W^{1,2})}^2\,.
\eeq{x3}
By \lemmref{1},
$$
\|\mathbb D(\bfu_m)\|_2^2\ge c(R)\,(\|\bfu_m\|_2^2+M\,|\bfgamma_m|^2)\,,
$$
which, once replaced in \eqref{x3} in conjunction with Gronwall lemma, allows us to deduce
$$
\|\bfu_m(t)\|_2^2+M\,|\bfgamma_m(t)|^2\le \left(\|\bfu_m(0)\|_2^2+M\,|\bfgamma_m(0)|^2\right)\,{\rm e}^{-c_0\,t}+ c_1\,\|\bfcalv\|_{L^\infty(W^{1,2})}^2\,,
$$
for some positive constant $c_0$. Employing this inequality we deduce, on the one hand, by \eqref{norm}, that $|c_{mi}(t)|$ are uniformly bounded, which implies the existence of a global, unique solution to the approximating problem \eqref{x1}; and, on the other hand, that the map  $\{c_{mi}(0)\}\mapsto \{c_{mi}(T)\}$ must have a fixed point, which implies that \eqref{x1} admits a $T$-periodic solution for each $m\in\nat$; see \cite[Lemma 3.1]{GS1} for details. Furthermore, from \eqref{x3}, \eqref{WH2} and \lemmref{xten} we deduce the {\em uniform estimate
} (in both $m$ and $k$)
\be
\|\bfgamma_m\|_{L^2(0,T)}+\|\mathbb D(\bfu_m)\|_{L^2(L^2)}\le c\,\delta \,.
\eeq{x4}
The latter together with \lemmref{1} implies, in particular, the existence of $t_m\in (0,T)$ such that
$$
\|\bfu_m(t_m)\|_2+M\,|\bfgamma(t_m)|\le c(R)\,\delta\,.
$$
as a result, integrating both sides of \eqref{x3} from $t_m$ to $t_m+T$ and using the periodicity of $\bfu_m$, we get, in particular,
\be
\|\bfu_m\|_{L^\infty(L^2)}\le c(R)\,\delta\,.
\eeq{x5}
Combining \eqref{x4}--\eqref{x5} with well-known procedures (see {\em e.g.} \cite[Section 3]{GaRev}), we prove the existence of a field $\bfu^{(k)}$ and of a subsequence, again denoted by $\{\bfu_{m}\}$, such that
\be\ba{ll}\medskip
\bfu^{(k)}\in L^2(0,T;\calh(\Omega_R))\cap L^\infty(0,T;L^2(\Omega_R))\,;\\ \medskip
\bfu_{m} \to \bfu^{(k)} \ \mbox{weakly in $L^2(0,T;\calh(\Omega_R))$}\,,\\   
\medskip
\bfu_{m} \to \bfu^{(k)} \ \mbox{strongly in $L^2(0,T;L^2(\Omega_R))$}\,,\\
\medskip
\bfu_{m}(t) \to \bfu^{(k)}(t) \ \mbox{weakly in $L^2(\Omega_R)$\,, for all $t\in[0,T]$}\,,\\
\bfgamma_m\to \bfgamma_\sub{\bfu^{(k)}}\ \mbox{strongly in $L^2(0,T)$}\,.
\ea
\eeq{pr}
Recalling that $\bfu_{m}(0)=\bfu_{m}(T)$, $m\in\nat$, the last condition in \eqref{pr} furnishes that $\bfu^{(k)}(0)=\bfu^{(k)}(T)$, namely, that $\bfu^{(k)}$ is $T-$periodic. Moreover, in view of \eqref{x4}, we find
\be
\|\bfgamma_\sub{\bfu^{(k)}}\|_{L^2(0,T)}+\|\mathbb D(\bfu^{(k)})\|_{L^2(L^2)}\le c\,\delta\,,
\eeq{perche}
with $c$ independent of $k$. Integrating \eqref{x1} over $(0,T)$ and using the convergence properties  in \eqref{pr}, one can show that $\bfv^{(k)}:=\bfu^{(k)}+\bfcalv$ is a solution to \eqref{ana1k} satisfying the requirements of \defref{1} with $\Omega\equiv\Omega_{R_k}$; see \cite[Section 3]{GS1} for details. The final step is to let $k\to\infty$ in \eqref{ana1k} and show that $\bfv^{(k)}$ converges, in a suitable sense, to a weak solution to problem \eqref{SE}. To this end, one extends every  element of the sequence $\{\bfv^{(k)}\}$ to be identically 0 outside $\Omega_{R_k}$. In this way the extended fields, again denoted by $\bfv^{(k)}$, satisfy  \eqref{ana1}, for any arbitrarily fixed $\bfpsi\in \calc_\sharp(\Omega)$, provided $k$ is chosen sufficiently large. From \eqref{perche} it immediately follows the existence of $\bfu\in L^2(0,T;\calh(\Omega))$ such that
\be\ba{ll}\medskip
\bfu^{(k)}\to\bfu\,, \ \mbox{weakly in $L^2(0,T;\calh(\Omega))$}\,,
\\
\|\bfgamma_\sub{\bfu}\|_{L^2(0,T)}+\|\mathbb D(\bfu)\|_{L^2(L^2)}\le c\,\delta\,,
\ea
\eeq{xf}
which, in particular, by the properties of $\bfcalv$, gives
\be
\|\bfgamma_\sub{\bfu}\|_{L^2(0,T)}+\|\nabla\bfv\|_{L^2(L^2)}\le c\,\delta\,.
\eeq{xf1}
However, as is well known, the  type of convergence in \eqref{xf} is not enough to ensure the convergence of the nonlinear term. Nevertheless, this property can be established by means of Lions-Aubin lemma, provided one shows appropriate uniform bounds (in $k$) for $\partial_t\bfv^{(k)}$. In turn, the latter follows by combining  \lemmref{Hprime} and \eqref{perche}, which delivers 
$$
\|\partial_t\bfv^{(k)}\|_{L^1(\calh^{-1}(\Omega_\rho))}+\|\bfv^{(k)}\|_{L^2(W^{1,2}(\Omega_\rho))}\le c(\rho)\,\delta\,,\ \ \mbox{for each {\em fixed} $\rho\in (R_*,R_k)$.}
$$
Since the embedding $W^{1,2}(\Omega_\rho)\subset L^q(\Omega_\rho)$, $q\in [1,6)$, is compact and $L^q(\Omega_\rho)\subset \calh^{-1}(\Omega_\rho)$, by Aubin-Lions lemma we deduce, in particular, that there exists a subsequence, again denoted by $\{\bfv^{(k)}\}$ such that   
\be
\bfv^{(k)}\to \bfv\,,\ \ \mbox{strongly in $L^2(0,T;L^2(\Omega_\rho))$.}
\eeq{stc}
This  subsequence may depend on $\rho$.
However, covering $\Omega$ with an increasing sequence of bounded domains and using
Cantor diagonalization method, we may select a subsequence for which \eqref{stc} holds for all $\rho$. With \eqref{xf}$_1$  and \eqref{xf} in hand, it is then routine to show that the limiting field $\bfv$ is a weak solution in the sense of \defref{1}, which, by \eqref{xf1}  satisfies the estimate stated in the theorem; see \cite[Section 3]{GS1} or \cite[Section 3]{GaSi} for details.\par\hfill$\square$\medskip\par   

\setcounter{equation}{0}
\section{Sufficient Conditions for Self-Propulsion when $\bar{\bfv_*} \neq\0$} The results proved in the previous section allow us 
to find conditions on the boundary velocity distribution  $\bfv_*$ guaranteeing a non-zero net motion of the body, at least when   $\bfv_*$ has a non-zero average over a period. More precisely, provided $\bar{\bfv_*}\neq\0$,  we show that, in the class of weak solutions,  $\mathscr B$ can propel itself whenever $\bar{\bfv_*}$ has a non-zero projection on a suitable three-dimensional subspace of $L^2(\partial\Omega)$. The more involved case $\bar{\bfv_*}=\0$, will be instead addressed in the following sections. \par
To this end, we begin to introduce the following auxiliary fields $(\bfh^{(i)},p^{(i)})$, $i=1,2,3$, 
\be\ba{cc}\medskip\left.\ba{ll}\medskip
\nu\Delta\bfh^{(i)}=\nabla p^{(i)}\\
\Div\bfh^{(i)}=0\ea\right\}\ \ \mbox{in $\Omega$}\\ 
\bfh^{(i)}=\bfe_{i}\ \ \mbox{at $\partial\Omega $}\,.
\ea
\eeq{6.9}
From \cite[Lemma V.4.4]{Gab} it follows that there exists one and only one solution to \eqref{6.9} such that
\be
(\bfh^{(i)},p^{(i)})\in [D^{1,r}(\Omega)\cap L^{q}(\Omega)]\times L^r(\Omega)\,,\ \ \mbox{all $r\in (3/2,\infty)$ and $q\in (3,\infty)$}.
\eeq{6.9_0}
We  define the quantities
\be
\bfg_i:=\bfn\cdot\mathbb T(\bfh^{(i)},p^{(i)})\left.\right|_{\partial\Omega}\,,\ \ i=1,2,3\,,
\eeq{gi}
which represent stress vectors associated to the flows $(\bfh^{(i)},p^{(i)})$,  evaluated at the boundary, and introduce
the matrix $\mathbb A$ with components
\be
\mathbb A_{ki}=\int_{\partial\Omega} (g_i)_k\,,\ \ k,i=1,2,3\,.
\eeq{A}
It is well known that $\mathbb A$ is symmetric and non-singular \cite[Sections 5.2--5.4]{HB}.
The proof of the following result is given in \cite[Theorem 2.1]{GaSP}.
\Bl Let $\bfV_*$ be as in \theoref{Exi}. Then, the problem
\be\ba{cc}\medskip
[\bfu_0,\bfpsi]= 0\,,\, \ \mbox{for all $\bfpsi\in\calh(\Omega)$}\\ \medskip
\Div\bfu_0=0\,,\\
\bfu_0=\bar{\bfV_*}+\bfxi_0\,,\, \ \mbox{at $\partial\Omega$}\,,
\ea
\eeq{Stsp}
has one and only one solution $(\bfu_0,\bfxi_0)\in D^{1,2}(\Omega)\times \real^3$. Moreover, 
\be
\bfxi_0=\mathbb A^{-1}\cdot\bfpzc{G}\,,
\eeq{speed}
where
\be
\bfpzc G:=\sum_{i=1}^3\left(\int_\Omega\bar{\bfV_*}\cdot\bfg_i\right)\bfe_i\,.
\eeq{pzc}
\EL{stoSP}
\par
Next, let $(\bfv,\bfgamma)$ be the weak solution determined in \theoref{Exi} and define the scaled quantities $(\bfu,\bfxi)$ and $(\bfw,\bfchi)$ as follows
\be
\bar{\bfv}=\delta\,\bfu\,,\ \ \bar{\bfgamma}=\delta\,\bfxi\,;\ \ \delta\,\bfw=\bfv-\bar{\bfv}\,,\ \ \delta\,\bfchi=\bfgamma-\bar{\bfgamma}.
\eeq{scaling}
Then, by taking, in particular, in \eqref{ana1} $\bfpsi\in \calc(\Omega)$, 
we get that $\bfu$ satisfies
\be
2  \,\nu\,[\bfu,\bfpsi]
= \delta\left[ \big((\bfu-\bfxi)\cdot\nabla\bfpsi,\bfu)+\bar{ \big((\bfw-\bfchi)\cdot\nabla\bfpsi,\bfw)}\right]\,,\ \ \mbox{for all $\bfpsi\in\calc(\Omega)$}\,,
\eeq{ana1_0}
and, in addition,
$$
\bfu=\bar{\bfV_*}+\bfxi\ \ \mbox{at $\partial\Omega$}\,.
$$
\Bl Let $(\bfv,\bfgamma)$ be the weak solution constructed in \theoref{Exi}. Then, as $\delta\to 0$, 
$$
\bfu-\bfu_0\to \0\,\  \mbox{weakly in $\calh(\Omega)$}\,;\ \ \bfxi\to\bfxi_0\, \ \mbox{in $\real^3$}\,,
$$
where $(\bfu_0,\bfxi_0)$ is the solution to the problem \eqref{Stsp} given in \lemmref{stoSP}. 
\EL{3.2}
{\em Proof.} Let $\{\delta_n\}\subset (0,\delta_0)$ with $\delta_n\to 0$ as $n\to\infty$, and let $(\bfv_n,\bfgamma_n)$ be the corresponding weak solutions of \theoref{Exi}. From \eqref{stima},  \eqref{scaling}, observing that $(\bfu_n-\bar{\bfV})\in\calh(\Omega)$, and also recalling that $D^{1,2}(\Omega)\subset W^{1,2}(\Omega_R)$, we readily deduce that there is $(\hat{\bfu},\hat{\bfxi})\in D^{1,2}(\Omega)\times \real^3$ such that (possibly, along a subsequence)
\be\ba{cc}\medskip
\nabla\bfu_n\to\nabla\hat{\bfu}\,\  \mbox{weakly in $L^2(\Omega)$}\,;\ \ \bfxi_n\to\hat{\bfxi}\, \ \mbox{in $\real^3$}\,,\\
\hat{\bfu}=\bar{\bfV_*}+\hat{\bfxi}\ \ \mbox{at $\partial\Omega$}\,.
\ea
\eeq{conve}
In view of \eqref{ana1_0}, we at once deduce that $(\bfu_n,\bfxi_n,\bfw_n,\bfchi_n)$ satisfies
\be
2  \,\nu\,[\bfu_n,\bfpsi]
= \delta_n\left[ \big((\bfu_n-\bfxi_n)\cdot\nabla\bfpsi,\bfu_n)+\bar{ \big((\bfw_n-\bfchi_n)\cdot\nabla\bfpsi,\bfw_n)}\right]\,,\ \ \mbox{for all $\bfpsi\in\calc(\Omega)$}\,,
\eeq{ana1_n}
Denote by $\cali_n$ the quantity in brackets on the right-hand side of \eqref{ana1_n}.
By applying Schwarz inequality, we show 
\be
|\cali_n|\le \|\nabla\bfpsi\|_\infty\left[\|\bfu_n\|_{2,\Omega_\rho}(|\bfxi_n|+\|\bfu_n\|_{2,\Omega_\rho})+\|\bfw_n\|_{L^2(L^2(\Omega_\rho)}(\|\bfchi_n\|_{L^2(0,T)}+
\|\bfw_n\|_{L^2(L^2(\Omega_\rho)}\right]\,,
\eeq{sn}
where $\Omega_\rho\supset\supp(\bfpsi)$. From \eqref{stima} and the scaling \eqref{scaling} we get
$$
|\bfxi_n|+\|\bfchi_n\|_{L^2(0,T)}\le c\,.
$$
with $c$ independent of $n$.
Furthermore, we observe that, for a.a. $t\in [0,T]$, $\bfv_n(\cdot,t)-\delta\,\bfV(\cdot,t)$ is in $\calh(\Omega)$, and so, by \eqref{SoB}, \lemmref{xten} and \eqref{stima} we infer
$$
\|\bfv_n\|_{L^2(L^2(\Omega_\rho))}\le c\,\delta_n\,,
$$
which implies
$$
\|\bfu_n\|_{2,\Omega_\rho}+\|\bfw_n\|_{L^2(L^2(\Omega_\rho))}\le c\,,
$$
with the constant $c$ independent of $n$. Replacing the above information back in \eqref{sn} and using again \eqref{stima}, we conclude that $\cali_n$ is uniformly bounded in $n$. Therefore, passing to the limit $n\to\infty$ in \eqref{ana1_n} and employing \eqref{conve} we conclude
$$
[\hat{\bfu},\bfpsi]=0\,,\, \ \mbox{for all $\bfpsi\in \calc(\Omega)$}\,;\ \ \Div\hat{\bfu}=0\,;\ \ \hat{\bfu}=\bar{\bfV_*}+\hat{\bfxi}\, \ \mbox{at $\partial\Omega$}\,,
$$
that is, $(\hat{\bfu},\hat{\bfxi})\in D^{1,2}(\Omega)\times\real^3$ is a solution to \eqref{Stsp}. However, by \lemmref{stoSP}, the solution is unique, which gives, on the one hand, $(\hat{\bfu},\hat{\bfxi})\equiv(\bfu_0,\bfxi_0)$ and, on the other hand, that the convergence in \eqref{conve} holds not only along a
subsequence, but as long as $\delta\to 0$.
Finally, since $\bfu-\bfu_0\in \calh$, \eqref{conve}$_1$ gives $\bfu-\bfu_0\to\0$, weakly in $\calh(\Omega)$, as claimed. The proof of the lemma is completed.
\par\hfill$\square$\medskip\par
We are now in a position to show the main result of this section.
\Bt Let $(\bfv,\bfgamma)$ be the weak solution to \eqref{SE} given in \theoref{Exi}. If $\bfV_*$ is such that the corresponding vector $\bfpzc G$ in \eqref{pzc} does not vanish, then also $\bar{\bfgamma}\neq\0$ and we have
\be 
\bar{\bfgamma}=\delta\,\mathbb A^{-1}\cdot\bfpzc G+o(\delta)\,, \,\ \mbox{as $\delta\to 0$.}
\eeq{delta}
In other words, under the above assumption, the body $\mathscr B$ can self-propel.
\ET{NZA}
{\em Proof.} Recalling  that $\bar{\bfgamma}=\delta \,\bfxi$, the theorem is an immediate corollary of the previous lemma and \eqref{speed}.\par 
\par\hfill$\square$\medskip\par

\Br \theoref{NZA} shows, in particular, that, if $\bar{\bfV_*}=\0$,  self-propulsion is a phenomenon to be searched at an order in $\delta$ higher than 1. It also shows that, since the linearized approximation possesses in this case only the identically vanishing solution (see \lemmref{stoSP}), the solution  to the nonlinear problem that would ensure self-propulsion {\em cannot} be obtained by a perturbation argument around its linear counterpart. In other words, when $\bar{\bfV_*}=\0$ self-propulsion  is a {\em strictly nonlinear} phenomenon. As we will show later on, its resolution will require a suitable contradiction argument directly applied on 
the full  set of nonlinear equations. In order to reach this goal, however, we need to prove  well-posednes of the problem \eqref{SE} in a class of  solutions more regular than weak solutions. For such a purpose, we will establish some key results on certain steady-state and time-periodic problems, which will be the main object of the following two sections.  
\ER{3.1}
\setcounter{equation}{0}
\section{On the Resolution of  a Nonlinear Steady-State Problem}
The main objective of this section is to show existence, uniqueness and corresponding estimates for solutions to  the following boundary-value problem in the unknowns $\bfu$, $p$ and $\bfxi_\sub{\bfu}$: 
\be\ba{cc}\medskip\left.\ba{ll}\medskip
\nu\Delta\bfu+\bfxi_\sub{\bfu}\cdot\nabla\bfu-\nabla p=\bff\\ \Div\bfu=0\ea\right\}\ \ \mbox{in $\Omega$}\,,
\\ 
\medskip
\bfu=\bfxi_\sub{\bfu}\ \ \mbox{on $\partial\Omega$\,,}\ \ \Lim{|x|\to \infty}\bfu(x)=\0\,,\\
\Int{\partial\Omega}{}[\mathbb T(\bfu,p)-{\bfscr{F}}]\cdot\bfn=\0\,,
\ea
\eeq{St1}
where  $\bff=\Div\bfscr F$. Notice that \eqref{St1} is a {\em nonlinear} problem. 
Precisely, we have the following.
\Bl Let $\bff\in\sfL^q$, $q\in(1,2]$. Then, there is  $(\bfu,p)\in \sfD^{2,q}\times \sfD^{1,q}$ solving \eqref{St1} that, in addition, satisfies the estimate
\be\ba{rl}\medskip
\|\bfu\|_{{\sf D}^{2,q}}+\|p\|_{{\sf D}^{1,q}}&\!\!\!\!\le C_1\left(\|\bff\|_2+\|\bfscr{F}\|_{q}+\|\bfscr{F}\|_{2}+\|\bfscr{F}\|_{2}^2\right)\\
|\bfxi_\sub{\bfu}|&\!\!\!\!\le c  \,\|\bfscr F\|_2
\,.\ea
\eeq{K0}
Suppose, next, $q\in (1,\frac65]$ and let $(\bfu_1,p_1)\in \sfD^{2,q}\times \sfD^{1,q}$ be another solution to \eqref{St1} corresponding to $\bff_1=\Div\bfscr F_1\in {\sf L}^q$ with $(\bfscr F_1-\bfscr F)\in L^{\frac{3q}{3-q}}(\Omega)$. Then,  there exists $c_0>0$ such that if 
\be
\|\bff\|_{\sfL^q}<c_0\,,\eeq{hyp} 
it follows that $\bfU:=\bfu_1-\bfu_2\in \mathsf D^{2,\frac{3q}{3-q}}$, $P:=p_1-p\in \mathsf D^{1,\frac{3q}{3-q}}$ and, in addition, (with  $\bfscr H:=\bfscr F_1-\bfscr F$)
\be \ba{ll}\medskip
\|\bfU\|_{\frac{3q}{3-2q}}+|\bfU|_{1,\frac{3q}{3-q}}+\|P\|_{\frac{3q}{3-q}}\le c_1(\|\bfscr H\|_2+\|\bfscr H\|_{\frac{3q}{3-q}})\,,\\
|\bfU|_{2,2}+|P|_{1,2}\le c\,\left(\|\Div\bfscr H\|_2+\|\bfscr H\|_2+\|\bfscr H\|_2^2+\|\bfscr H\|_{\frac{3q}{3-q}}\right)
\ea
\eeq{FP}
Thus, in particular, if $\bff\equiv\bff_1$, under the assumption \eqref{hyp} the solution $(\bfu,p)$ is unique in the  class $\sfD^{2,q}\times \sfD^{1,q}$, $q\in (1,\frac65]$.
\EL{Fm}
{\em Proof.} Since the actual value of $\nu$ is irrelevant, we set, for simplicity, $\nu=1$. We begin to put \eqref{St1} in a weak form. To this end, we dot-multiply both sides of \eqref{St1}$_1$ by $\bfphi\in\calc$ and integrate by parts as necessary. We get
$$
-[\bfu,\bfphi]+\bfxi_\sub{\bfphi}\cdot\Int{\partial\Omega}{}[\mathbb T(\bfu,p)-{\bfscr{F}}]\cdot\bfn+(\bfxi_\sub{\bfu}\cdot\nabla\bfu,\bfphi)=-(\bfscr{\bfscr F},\nabla\bfphi)\,,
$$
so that,  reinforcing \eqref{St1} in the preceding relation, we find
\be
-[\bfu, \bfphi]+(\bfxi_\sub{\bfu}\cdot\nabla\bfu,\bfphi)=-(\bfscr{\bfscr F},\nabla\bfphi)\,,\ \ \mbox{for all $\bfphi\in\calc$.}
\eeq{k}
It is also easy to see, conversely, that if $(\bfu,p)$ satisfies \eqref{k} and is sufficiently smooth, then it satisfies  \eqref{St1}. Now, by formally replacing $\bfu$ for $\bfphi$ in \eqref{k1}, and then using Schwarz inequality along with \eqref{WH1}, we obtain  
\be
\|\mathbb D(\bfu)\|_2\le \sqrt{2}\,\|\bfscr{F}\|_2\,,
\eeq{k1}
which, in particular, by \eqref{WH2} implies \eqref{K0}$_2$.  
Combining \eqref{k1} with the Galerkin method, one can show that \eqref{k} has at least one solution $\bfu\in\calh$. To such a solution, we can then associate a pressure field $p$ with $p\in L^2_{\rm loc}(\bar\Omega)$ and such that
\be
\|p\|_{2,\Omega_R}\le c_R\left (\|\bfscr{F}\|_2+|\bfxi_\sub{\bfu}|\,\|\bfu\|_{2,\Omega_R}\right)\,;
\eeq{k2}
see \cite[Lemma V.1.1 and Lemma VII.1.1]{Gab}. Since $\bfu\in\calh$, by \eqref{St1} and \eqref{WH2} we infer, in a distributional sense,
\be 
\Div[\mathbb T(\bfu,p)-\bfscr{F}]=-\bfxi_u\cdot\nabla\bfu\in L^2(\Omega)
\eeq{k0}
and, by \eqref{k1}, \eqref{k2} and the assumption, 
$$
\left(\mathbb T(\bfu,p)-\bfscr{F}\right)\in L^2(\Omega_R)\,.
$$
As a consequence, e.g., \cite[Theorem III.3.2]{Gab}, 
\be\ba{ll}\medskip
\left(\mathbb T(\bfu,p)-\bfscr{F}\right)\cdot\bfn|_{\partial\Omega}\in W^{-\frac12,2}(\partial\Omega)\,,\\ 
\|\left(\mathbb T(\bfu,p)-\bfscr{F}\right)\cdot\bfn\|_{W^{-\frac12,2}(\partial\Omega)}\le c\,\left (\|\bfscr{F}\|_2+|\bfxi_\sub{\bfu}|(\,\|\bfu\|_{2,\Omega_R}+\|\nabla\bfu\|_2)\right)\,.
\ea
\eeq{k3}
Following \cite{KoSo}, we shall next  make a suitable extension of problem \eqref{St1} to the whole space $\real^3$. Consider the vector-valued Neumann problem:
\be
\Delta \bfP=\0 \ \ \mbox{in $\Omega_0$}\,;\ \ \pde{\bfP}{\bfn}= \left(\mathbb T(\bfu,p)-\bfscr{\bfscr F}\right)\cdot\bfn\ \ \mbox{at $\partial\Omega$}\,.
\eeq{NP}
In view of \eqref{k3}, and the fact that,  by \eqref{k}, we have
$$
 \int_{\partial\Omega}\left(\mathbb T(\bfu,p)-\bfscr{\bfscr F}\right)\cdot\bfn=0\,,
$$
in distributional sense, we may assert that
\eqref{NP} admits one and only one (up to a constant) solution $\bfP\in W^{1,2}(\Omega_0)$, and that, setting ${\bfscr G}:=\nabla \bfP$, the following estimate holds
\be
\|\bfscr{G}\|_2\le c\,\|\left(\mathbb T(\bfu,p)+\bfxi_\sub{\bfu}\otimes\bfxi_\sub{\bfu}-\bfscr{F}\right)\cdot\bfn\|_{W^{-\frac12,2}(\partial\Omega)}\le c\,\left (\|\bfscr{F}\|_2+|\bfxi_\sub{\bfu}|(\,\|\bfu\|_{2,\Omega_R}+\|\nabla\bfu\|_2)\right)\,,
\eeq{k4}
where we have used \eqref{WH1} and \eqref{WH2}.
Define
$$
\tilde{{\bfscr F}}(x)=\left\{\ba{ll}\medskip\bfscr{F}(x)\ &\mbox{if $x\in\Omega$}\\ 
{\bfscr G}(x)\ &\mbox{if $x\in\Omega_0$}\ea\right.\,.
$$
Taking into account that, by \lemmref{1}, $\|\bfu\|_{2,\Omega_R}\le c\|\bfu\|_6\le c|\bfu|_{1,2}$, we deduce with the help of \eqref{k1}, \eqref{k4} and, again, \lemmref{1} 
\be
\|\tilde{\bfscr{F}}\|_{q}+\|\tilde{\bfscr{\bfscr F}}\|_2\le c\,\left(\|\bfscr{F}\|_{q}+(1+|\bfxi_\sub{\bfu}|)\|\bfscr{\bfscr F}\|_2\right)
\eeq{K01}
Consider now the problem
\be\ba{ll}\medskip
\left.\ba{ll}\medskip
\Delta\bfv+\bfxi_\sub{\bfu}\cdot\nabla\bfv-\nabla {\sf p}=\Div\tilde{\bfscr{F}}\\ \Div\bfv=0\ea\right\}\ \ \mbox{in $\real^3$}\,.
\ea
\eeq{k5}
It is well known that \eqref{k5} has one and only one distributional solution $(\bfv,{\sf p})\in [D^{1,2}_0(\real^3)\cap D_0^{1,q}(\real^3)]\times [L^{2}(\real^3)\cap L^{q}(\real^3)]$ such that
\be
|\bfv|_{1,2}+|\bfv|_{1,q}+\|{\sf p}\|_2+\|{\sf p}\|_{q}\le c\,\left(\|\tilde{\bfscr{\bfscr F}}\|_2+\|\tilde{\bfscr{\bfscr F}}\|_{q}\right)\,;
\eeq{k6}
see \cite[Theorem IV.2.2 and Theorem VII.4.2]{Gab}. We now extend $(\bfu,p)$ to the whole of $\real^3$ by setting
$$
\tilde{\bfu}=\left\{\ba{ll}\medskip\bfu(x)\ &\mbox{if $x\in\Omega$}\\ 
\bfxi_\sub{\bfu}\ &\mbox{if $x\in\Omega_0$}\ea\right.\ \ \tilde{p}=\left\{\ba{ll}\medskip p(x)\ &\mbox{if $x\in\Omega$}\\ 
0\ &\mbox{if $x\in\Omega_0$}\ea\right.\,.
$$
It is easy to see that $\tilde{\bfu},\tilde{p}$ is  a solution to \eqref{k5} in the sense of distribution. In fact, for all $\psi\in C_0^\infty(\real^3)$,
$$\ba{ll}\medskip
\Int{\real^3}{}[\mathbb T(\tilde{\bfu},\tilde{p})+\bfxi_\sub{\bfu}\otimes\bfu-\tilde{\bfscr F}]\cdot\nabla\psi=\Int{\Omega}{}[\mathbb T({\bfu},{p})+\bfxi_\sub{\bfu}\otimes\bfu-\bfscr{\bfscr F}]\cdot\nabla\psi
+\Int{\Omega_0}{}(\bfxi_\sub{\bfu}\otimes\bfxi_\sub{\bfu}-\bfscr G)\cdot\nabla\psi\\
=\Int{\partial\Omega}{}[\mathbb T({\bfu},{p})+\bfxi_\sub{\bfu}\otimes\bfxi_\sub{\bfu}-\bfscr{\bfscr F}]\cdot\bfn\psi-\Int{\partial\Omega}{}[\bfxi_\sub{\bfu}\otimes\bfxi_\sub{\bfu}-\bfscr{\bfscr G}]\cdot\bfn\psi=\0
\ea
$$
We then deduce that $\bfw:=\tilde{\bfu}-\bfv$ satisfies the following equation
 in distributional sense
$$
\left.\ba{ll}\medskip
\Delta\bfw+\bfxi_\sub{\bfu}\cdot\nabla\bfw-\nabla {\sf \phi}=\0\\ \Div\bfw=0\ea\right\}\ \ \mbox{in $\real^3$}\,.
$$
Therefore, since $\bfw\in D^{1,2}_0(\real^3)$, it  must be $\bfw\equiv\0$ \cite[Theorem IV.2.2 and Theorem VII.4.2]{Gab} which implies, in particular, $\bfu\equiv\bfv$ in $\Omega$. Thus, we have $(\bfu,p)\in D^{1,q}(\Omega)\times [L^{q}(\Omega)\cap L^2(\Omega)]$  and also, from \eqref{k6}, \eqref{K01},
\be
|\bfu|_{1,2}+|\bfu|_{1,q}+\|{p}\|_2+\|{ p}\|_{q}\le c\,\left(\|{\bfscr F}\|_{q}+(1+|\bfxi_\sub{\bfu}|)\|\bfscr{\bfscr F}\|_2\right)\,.
\eeq{k7}
We next observe that, from \eqref{St1}, we obtain that $(\bfu,p)$ can be viewed as a (distributional) solution to the following Stokes problem
\be\ba{cc}\medskip\left.\ba{ll}\medskip
\Delta\bfu-\nabla p=\bfF
\\ \Div\bfu=0\ea\right\}\ \ \mbox{in $\Omega$}\,,\\
\bfu(x)=\bfxi_\sub{\bfu}\ \ \mbox{on $\partial\Omega$\,,}
\ea
\eeq{k8}
where $\bfF:=-\bfxi_\sub{\bfu}\cdot\nabla\bfu+\bff$.
By assumption and what we have shown so far, $\bfF$ is in $L^2(\Omega)$, and so by \cite[Theorem V.5.3]{Gab} we deduce, on the one hand, $(\bfu,p)\in D^{2,2}(\Omega^R)\times D^{1,2}(\Omega^R)$ for all $R>R_*$, and, on the other hand, by local estimates for the Stokes problem \cite[Theorem IV.5.1  and Remark IV.5.1]{Gab} 
$$
|\bfu|_{2,2,\Omega_{2R}}\le c_R\,(\|\bfF\|_{2,\Omega_{3R}}+\|\bfu\|_{2,\Omega_{3R}}+|\bfxi_\sub{\bfu}|)\,.
$$
As a consequence, since by \eqref{SoB} and \eqref{k1}
\be
\|\bfu\|_{2,\Omega_R}\le c_R\|\bfu\|_{6,\Omega_R}\le c_R\|\bfscr F\|_2\,,
\eeq{ch}
from the latter, \lemmref{1}  and \eqref{St1}$_1$, we infer
\be 
(\bfu,p)\in D^{2,2}(\Omega)\times D^{1,2}(\Omega)\,.
\eeq{St11}
We may then use \cite[Lemma V.4.3]{Gab}, \eqref{K0}$_2$ and \eqref{ch},  to show the estimate
$$
|\bfu|_{2,2}+|p|_{1,2}\le c\left(\|\bff\|_{2}+(1+|\bfxi_\sub{\bfu}|)\|\bfscr{F}\|_2\right)\,.
$$
The proof of existence is thus completed.
Let us next prove the validity of \eqref{FP}. Setting $\bfzeta:=\bfxi_\sub{\bfu_1}-\bfxi_\sub{\bfu}$, it follows that
\be\ba{cc}\medskip\left.\ba{ll}\medskip
\Delta\bfU+\bfxi_\sub{\bfu_1}\cdot\nabla\bfU=-\bfzeta\cdot\nabla\bfu+\nabla {P}+\Div\bfscr H\\ \Div\bfU=0\ea\right\}\ \ \mbox{in $\Omega$}\,,
\\ \medskip
\bfU=\bfzeta\ \ \mbox{at $\partial\Omega$}\,,\\
\Int{\partial\Omega}{}[\mathbb T(\bfU,{P})-\bfscr H]\cdot\bfn=\0\,.
\ea
\eeq{St0}
We  now dot-multiply both sides of \eqref{St0}$_1$ by $\bfU$ and integrate  by parts over $\Omega$. In view of the summability properties of $\bfU$ and $\bfu$,  and \eqref{St0}$_4$, we thus get
$$
\|\mathbb D(\bfU)\|_2^2=\int_\Omega\left[\bfzeta\cdot\nabla\bfU\cdot\bfu-\bfscr H:\nabla\bfU\right]\,,
$$
which, in turn, with the help of Schwarz inequality, \eqref{WH1} and \eqref{WH2} leads to
\be
|\bfU|_{1,2}\le c_1\,\left(\|\bfu\|_2|\bfU|_{1,2}+\|\bfscr H\|_2\right)\,.
\eeq{ari}
Since $q\in(1,\frac65]$, we may use the embedding in \lemmref{Emb0} along with \eqref{K0} to show
$$
\|\bfu\|_2\le c_2\,\left(\|\bff\|_2+\|\bfscr F\|_q+\|\bfscr F\|_2+\|\bfscr F\|_2^2\right)\le c_2\,\left(\|\bff\|_{\sf L^q}+\|\bff\|_{\sf L^q}^2\right)\,.
$$
From the latter and \eqref{ari} we deduce that there is a constant $c_0>0$ such that if \eqref{hyp} holds, then $\|\bfu\|_2\le 1/2c_1$, which once replaced in  \eqref{ari}, with the help of \eqref{WH2} allows us to conclude
\be
|\bfzeta|\le c\,|\bfU|_{1,2}\le c_3\,\|\bfscr H\|_2\,.
\eeq{ari1}
Set $s:=3q/(3-q)$. Clearly, $s\in (\frac32,3)$ and, as a result, we can apply \cite[Theorem V.5.1 and Theorem  VII.5.2]{Gab} to the boundary-value problem \eqref{St0}$_{1,2,3}$ and obtain the following inequality
\be
\|\bfU\|_{\frac{3s}{3-s}}+|\bfU|_{1,s}+\|P\|_s\le c\,\left(|\bfzeta|\,\|\bfu\|_s+\|\bfscr H\|_s+|\bfzeta|\right)\,.
\eeq{ari2}
Our next step will be to give a suitable estimate of the right-hand side of \eqref{ari2}. From \lemmref{HSS} and \eqref{K0} we find
\be
\|\bfu\|_s\le c\,\left(\|\bff\|_{\sf L^q}+\|\bff\|_{\sf L^q}^2\right)\,.
\eeq{ari3}
Furthermore, setting
$$
\tilde{\bfU}=\left\{\ba{ll}\medskip\bfU(x)\ &\mbox{if $x\in\Omega$}\\ 
\bfzeta\ &\mbox{if $x\in\Omega_0$}\ea\right.\,,
$$
we have
$$
|\bfzeta|\,|\Omega_0|^{\frac{3-s}{3s}}=\left(\int_{\Omega_0}{|\bfzeta|}^{\frac{3s}{3-s}}\right)^{\frac{3-s}{3s}}\le \|\tilde{\bfU}\|_{\frac{3s}{3-s}}\,,
$$
which, by \lemmref{HSS}, entails
\be
|\bfzeta|\le c\,|\tilde{\bfU}|_{1,s}=c\,|{\bfU}|_{1,s}\,.
\eeq{ari4}
Combining \eqref{ari1}--\eqref{ari4}, we thus prove that there exists a constant $c_0$ such that if \eqref{hyp} holds, then
\be
\|\bfU\|_{\frac{3s}{3-s}}+|\bfU|_{1,s}+\|P\|_s\le c\,\left(\|\bfscr H\|_2+\|\bfscr H\|_s\right)\,,
\eeq{ari5}
which proves \eqref{FP}$_1$. 
We now pass to the estimate of the second derivatives of $\bfU$. Applying \cite[Lemma V.4.3]{Gab} to \eqref{St0} we infer
\be
|\bfU|_{2,2}+|P|_{1,2}\le c\,\left(\|\Div\bfscr H\|_2+|\bfxi_\sub{\bfu_1}|\,|\bfU|_{1,2}+|\bfzeta|\,|\bfu|_{1,2}+\|\bfU\|_{2,\Omega_R}+\|P\|_{2,\Omega_R}\right)\,.
\eeq{ari6}
From \eqref{K0},  \eqref{hyp}, \eqref{ari4}, and \eqref{ari5} it follows that
\be
|\bfzeta|\,|\bfu|_{1,2}\le c\,|\bfU|_{1,s}\le c\,\left(\|\bfscr H\|_2+\|\bfscr H\|_s\right)\,,
\eeq{ari7} 
whereas by \eqref{WH1}, \eqref{SoB} and \eqref{ari1}
\be
\|\bfU\|_{2,\Omega_R}\le c\, \|\bfU\|_{6,\Omega_R}\le c\, |\bfU|_{1,2}\le c\,\|\bfscr H\|_2\,.
\eeq{ari8}
Recalling that $q\in(1,\frac65]$, we have $s=2$ if $q=\frac65$, otherwise $s\in (\frac32,2)$. In the first case, from \eqref{ari5} we have
\be
\|P\|_{2,\Omega_R}\le c\,\left(\|\bfscr H\|_2+\|\bfscr H\|_s\right)\,.
\eeq{ari9}
In the second case, we observe that since $P\in D^{1,2}(\Omega)\cap L^s(\Omega)$, by \lemmref{HSS}, we get $P\in L^6(\Omega)$ with
$$
\|P\|_6\le c\,|P|_{1,2}\,.
$$
Thus, if $s<2$, by elementary interpolation,  the latter inequality and \eqref{ari5} we show, for arbitrary $\varepsilon>0$,
\be
\|P\|_{2,\Omega_R}\le c_\varepsilon\|P\|_s+\varepsilon\,\|P\|_6 \le c_\varepsilon\|P\|_s+\varepsilon\,C\,|P|_{1,2}\le c\,\left(\|\bfscr H\|_2+\|\bfscr H\|_s\right)+\varepsilon\,C\, |P|_{1,2}\,.
\eeq{ari10}
From \eqref{ari6}--\eqref{ari10} and \eqref{ari1} we then deduce
\be
|\bfU|_{2,2}+|P|_{1,2}\le c\,\left[\|\Div\bfscr H\|_2+(|\bfxi_\sub{\bfu_1}|+1)\|\bfscr H\|_2+\|\bfscr H\|_s\right]\,.
\eeq{ari11}
Finally, from
$$
|\bfxi_\sub{\bfu_1}|\le |\bfzeta|+|\bfxi_\sub{\bfu}|\,,
$$ 
and \eqref{K0}, \eqref{hyp} and \eqref{ari1}, we show
$$
|\bfxi_\sub{\bfu_1}|\le c\,(1+\|\bfscr H\|_2)\,,
$$
which, once replaced in \eqref{ari11}, proves \eqref{FP}$_2$. 
 The proof of the lemma is thus  completed. 
\par\hfill$\square$\medskip\par

\setcounter{equation}{0}
\section{On the Resolution of a Time-Periodic Linear Problem}
Our next task is to study the well-posedness of the following problem
\be\ba{cc}\medskip\left.\ba{ll}\medskip
\partial_t\bfw-\nu\Delta\bfw=-\nabla {\sfp}+\bff\\ 
\Div\bfw=0\ea\right\}\ \ \mbox{in $\Omega\times\real$}\,,\\ 
\medskip
\bfw=\bfchi_\sub{\bfw}\ \ \mbox{on $\partial\Omega\times\real$}\,,\\ 
M\,\dot{\bfchi}_\sub{\bfw}+\Int{\partial\Omega}{}\bfT(\bfw,{\sf p})\cdot\bfn=\bfF
\ea
\eeq{2.00}
in a suitable function class of $T$-periodic solutions, 
with $\bff=\bff(x,t)$ and $\bfF(t)$  given $T$-periodic functions. In order to reach this goal, we need to premise a preparatory result.
\Bl  Consider the  boundary-value problems, $i=1,2,3$, $k\in\mathbb Z$: 
\be\ba{cc}\medskip\left.\ba{ll}\medskip
\i\,k\,\omega\,\bfH_k^{(i)}=\Delta \bfH_k^{(i)}-\nabla{\sf \gamma}^{(i)}_k\\
\Div\bfH_k^{(i)}=0\ea\right\}\ \ \mbox{in $\Omega$}\,,\\
\bfH_k^{(i)}|_{\partial\Omega}=\bfe_i\,,\ \ \ \bfH_0^{(i)}=\0\,,
\ea
\eeq{0.2}
with $\omega:=2\pi/T$. The following properties hold.
\begin{itemize}
\item[{\rm (i)}] There is one and only one solutions  $(\bfH_k^{(i)},\gamma_k^{(i)})\in W^{2,2}(\Omega)\times W^{1,2}(\Omega)$. This solution satisfies  the estimates
\be
\|\bfH_k^{(i)}\|_2
\le c\,;\ \ 
\|\nabla\bfH_k^{(i)}\|_2
\le c\,(|k|+1)^{\frac12}\,;\ \ 
|\bfH_k^{(i)}|_{2,2}
\le c\,(\,|k|+1)\,;\ \
c\,\le \|\nabla\bfH_k^{(i)}\|_2\,,
\eeq{0.0.0}
where $c$ is a constant independent of $k$.
\item[{\rm (ii)}] The matrix $\mathbb B$ defined by components
$$
(\mathbb B)_{\ell i}=\int_{\partial\Omega}\mathbb T_{\ell j}(\bfH^{(i)}_k,\gamma_k^{(i)})n_j
$$
satisfies the condition (with ${}^*\equiv \textrm{c.c.}$)
\be
\bfzeta^*\cdot\mathbb B\cdot\bfzeta=\i\,k\,\omega\,\|\zeta_i\bfH^{(i)}\|_2^2 +\|\zeta_i\nabla\bfH^{(i)}\|_2^2\,,
\eeq{B}
for all $\bfzeta=(\zeta_1,\zeta_2,\zeta_3)\in\mathbb C^3$, and it is, therefore, invertible.
\end{itemize}
\EL{1.1}
{\em Proof.} We begin to prove (i). Since the proof is the same for $i=1,2,3$, we chose $i=1$ and, for simplicity, omit the superscript. Let $\phi=\phi(|x|)$ be a (smooth) cut-off function such that
$$
\phi(\bfx)=\left\{\begin{array}{ll}\medskip
1\ & \mbox{in }\Omega_{R}\\
0\ & \mbox{in }\bar{\Omega^{2R}}\,,
\end{array}\right.
$$
and set
$$
\label{FG0}\bfPhi(\bfx)={\rm curl}\,\big(x_2\phi(\bfx)\,\bfe_3\big)\,.
$$
Clearly, $\Div\bfPhi=0$ and since $(\partial_i\equiv\partial/\partial x_i$)
$$
\bfPhi(\bfx)=-\bfe_3\times\nabla\big(x_2\phi(\bfx)\big)=x_2(\partial_2\phi(\bfx)\,\bfe_1-\partial_1\phi(\bfx)\,\bfe_2)+\phi(\bfx)\,\bfe_1\,,
$$
by the property of $\phi$ we deduce $\bfPhi(\bfx)=\bfe_1$ for all $\bfx\in \partial\Omega$. Therefore, $\bfPhi$ is a solenoidal extension of $\bfe_1$ with support contained in $\Omega_{2R}$. Setting $\bfv_k:=\bfH_k-\bfPhi$, from \eqref{0.2} we deduce that $\bfv_k$ satisfies the following boundary-value problem, for all $|k|\ge 1$:
\be\ba{cc}\medskip\left.\ba{ll}\medskip
\i\,k\,\omega\,\bfv_k=\Delta \bfv_k-\nabla{\gamma}_k-\i\,\,k\,\omega\,\bfPhi+\Delta\bfPhi\\
\Div\bfv_k=0\ea\right\}\ \ \mbox{in $\Omega$}\,,\\
\bfv_k|_{\partial\Omega}=\0\,.
\ea
\eeq{0.3}
Existence to \eqref{0.3} in the stated function class can be easily obtained by the Galerkin method combined with the estimate that we are about to derive. Let us dot-multiply both sides of \eqref{0.3}$_1$ by $\bfv_k^*$,  and integrate by parts as necessary. We obtain
\be
\i\,k\,\omega\,\|\bfv_k\|_2^2+\|\nabla\bfv_k\|_2^2=-\i\,\,k\,\omega\,(\bfPhi,\bfv_k^*)+(\Delta\bfPhi,\bfv_k^*)\,,
\eeq{0.4}
which, after taking its imaginary part,  and employing Schwarz inequality and the properties of $\bfPhi$, gives
\be
|k|\,\|\bfv_k\|_2\le c\,\left(\|\nabla\bfv_k\|_2+|k|+1\right)\,.
\eeq{0.5_1}
Similarly, by taking the real part of \eqref{0.4} and using Schwarz inequality, the properties of $\bfPhi$ and \eqref{0.5} it follows that
\be
\|\nabla\bfv_k\|_2^2\le c\,(|k|+1)\|\bfv_k\|_2\,.
\eeq{0.6_1}
Replacing \eqref{0.6_1} into \eqref{0.5_1} and then employing Cauchy-Schwarz inequality furnishes
\be
\|\bfv_k\|_2\le c\,,
\eeq{0.5}
with $c$ independent of $k$, which, in turn, by \eqref{0.6_1} implies
\be
\|\nabla\bfv_k\|_2\le c\,(|k|+1)^{\frac12}\,.
\eeq{0.6}
Furthermore, from classical estimates for the Stokes problem \cite[Lemma 1]{Hey}, from \eqref{0.3} and \eqref{0.6} we get
$$
|\bfv_k|_{2,2}\le c\, (|k|\,\|\bfv_k\|_2+\,|k|+1)
$$
which combined with \eqref{0.5} furnishes
\be
|\bfv_k|_{2,2}\le c\,(|k|+1)\,.
\eeq{0.7}
Recalling that $\bfv_k=\bfH_k-\bfPhi$, from \eqref{0.5}--\eqref{0.7} we readily conclude that $\bfH_k$ satisfies all the properties listed in \eqref{0.0}$_{1,2,3}$. 
Finally, from the trace inequality 
\be
\int_{\partial\Omega}|\bfe_1|^2\le c\,\|\nabla\bfH_k|^2_2
\eeq{0.9}
we prove also \eqref{0.0}$_4$, thus completing the proof of (i). In order to show (ii), set $\bfu:=\zeta_i\bfH^{(i)}$, $\phi=\zeta_i\gamma^{(i)}$ where, for simplicity, the subscript $"k"$ is omitted.\footnote{We use summation convention over repeated indices, unless confusion may arise.} From \eqref{0.2} we thus get
\be\ba{cc}\medskip\left.\ba{ll}\medskip
\i\,k\,\omega\,\bfu=\Delta \bfu-\nabla{\phi}\\
\Div\bfu=0\ea\right\}\ \ \mbox{in $\Omega$}\,,\\
\bfu|_{\partial\Omega}=\bfzeta\,.
\ea
\eeq{0.0}
If we dot-multiply both sides of \eqref{0.0}$_1$ by $\bfu^*$ and integrate over $\Omega$ we show
$$
\bfzeta^*\cdot\mathbb B\cdot\bfzeta=\i\,k\,\omega\,\|\bfu\|_2^2+\|\mathbb D(\bfu)\|_2^2\,,
$$
and so if 0 is an eigenvalue of $\mathbb B$, we must have $\bfu\equiv\0$, which, once evaluated at $\partial\Omega$ produces $\bfxi=\0$, and the proof of the lemma is completed.
\par\hfill$\square$\par 
\Bl Let $q\in(1,\infty)$,  $s=2$ if $q\in(1,2]$, and $s=q$ if $q\in(2,\infty)$. Then, for any $(\bff,\bfF)\in \call_\sharp^{2,q}\times L^s_\sharp(0,T)$,   problem \eqref{2.00}
has one and only one solution $\big(\bfw,\sfp,\bfchi_\sub{\bfw}\big)\in \mathcal W_\sharp^{2,q}
\times \mathcal P^{1,q}\times W^{1,s}_\sharp(0,T)$. This solution satisfies the estimate
\be
\|\bfw\|_{\mathcal W_\sharp^{2,q}}+\|{\sf p}\|_{\mathcal P^{1,q}}+\|\bfchi_\sub{\bfw}\|_{W^{1,2}_\sharp(0,T)}\le C_2\,\Big(\|\bff\|_{\call_\sharp^{2,q}}+\|\bfF\|_{L^s_\sharp(0,T)}\Big)\,.
\eeq{n}
\EL{1.6_0}
{\em Proof.} Since the actual value of  $\nu$ and $M$ is irrelevant to the proof, we set, for simplicity, $\nu=M=1$.  We  write $\bfw=\bfz+\bfu$ where $\bfz$ and $\bfu$ satisfy the following set of equations
\be\ba{cc}\medskip\left.\ba{ll}\medskip
\partial_t\bfz-\Delta\bfz=-\nabla {\sf \tau}+\bff\\ 
\Div\bfz=0\ea\right\}\ \ \mbox{in $\Omega\times\real$}\\
\bfz|_{\partial\Omega}=\0
\ea
\eeq{2.3} 
and
\be\ba{cc}\medskip\left.\ba{ll}\medskip
\partial_t\bfu-\Delta\bfu=-\nabla {\sf q}\\ 
\Div\bfu=0\ea\right\}\ \ \mbox{in $\Omega\times\real$}\\ \medskip
\bfu|_{\partial\Omega}=\bfchi_\sub{\bfu}\,;\\ \dot{\bfchi}_\sub{\bfu}+\Int{\partial\Omega}{}\bfT(\bfu,{\sf q})\cdot\bfn=\bfF-\Int{\partial\Omega}{}\bfT(\bfz,{\sf \tau})\cdot\bfn:=\bfcalf\,.
\ea
\eeq{2.4}
From \cite[Theorem 12]{GaMaH}, it follows that, under the stated assumptions, there exists a unique solution $(\bfz,\tau)\in  \mathcal W_\sharp^{2,q}\times\calp^{1,q}$, $q\in(1,\infty)$ that, in addition, obeys the inequality
\be
\|\bfz\|_{\mathcal W_\sharp^{2,q}}+\|\tau\|_{\mathcal P^{1,q}}\le c\,\|\bff\|_{\mathcal L_\sharp^{2,q}}\,,\ \ q\in(1,\infty).
\eeq{2.5}
Furthermore, by trace theorem,\footnote{Possibly, by modifying $\tau$ by adding to it a suitable function of time.} 
\be
\|\Int{\partial\Omega}{}\mathbb T(\bfz,{\sf \tau})\cdot\bfn\|_{L^q(L^q)}\le c\,\left(\|\bfz\|_{\mathcal W^{2,q}}+\|{\sf \tau}\|_{\mathcal P^{1,q}}\right)\,.
\eeq{2.6}
By \eqref{2.6} we  infer that the function $\bfcalf$ in \eqref{2.4} is in $L^2_\sharp(0,T)$. In order to find solutions to \eqref{2.4}, we
formally expand $\bfu$, ${\sf q}$, and  $\bfchi_\sub{\bfu}\,(\equiv\bfchi)$, in Fourier series:
\be
\bfu(x,t)=\Sum{k\in\mathbb Z}{}\bfu_k(x)\,{\rm e}^{\i k\,t}\,,\ \ {\sf q}(x,t)=\Sum{k\in\mathbb Z}{}{\sf q}_k(x)\,{\rm e}^{\i k\,t}\,,\ \
\bfchi(t)=\Sum{k\in\mathbb Z}{}\bfchi_k \,{\rm e}^{\i k\,t}\,,\ \ \bfu_0\equiv\nabla{\sf q}_0\equiv\bfchi_0\equiv\0\,,
\eeq{Fou}
where $(\bfu_k,{\sf q}_k,\bfchi_k)$ solve the problem
$k\neq0$
\be\ba{cc}\medskip\left.\ba{ll}\medskip
\i\,k\,\omega\,\bfu_k=\Delta \bfu_k-\nabla{\sf q}_k\\
\Div\bfu_k=0\ea\right\}\ \ \mbox{in $\Omega$}\\
\bfu_k|_{\partial\Omega}=\bfchi_k\,,
\ea
\eeq{2.8}
with the further condition
\be
\i\,k\,\omega\,{\bfchi_k}+\Int{\partial\Omega}{}\mathbb T(\bfu_k,{\sf q}_k)\cdot\bfn=\bfcalf_k\,,
\eeq{2.8_1}
where $\{\bfcalf_k\}$ are Fourier coefficients of $\bfcalf$ with $\bfcalf_0\equiv\0$. 
For each fixed $k\in\mathbb Z-\{0\}$, a  solution to \eqref{2.8}--\eqref{2.8_1} is given by
\be
\bfu_k=\sum_{i=1}^3\chi_{ki}\bfH_k^{(i)}\,,\ \ {\sf q}_k=\sum_{i=1}^3\chi_{ki}\gamma_k^{(i)}\,,
\eeq{1.20}
with $(\bfH_k^{(i)},\gamma_k^{(i)})$ given in \lemmref{1.1}, and where $\bfchi_k$ solve the equations
\be
\i\,k\,\omega\,{\bfchi_{k}}+\sum_1^3\chi_{ki}\Int{\partial\Omega}{}\mathbb T(\bfH_k^{(i)},{\gamma}^{(i)}_k)\cdot\bfn=\bfcalf_k\,.
\eeq{chi}
The latter, with the notation of \lemmref{1.1}(i), can be equivalently rewritten as
\be
\left(\i\,k\,\omega\,\mathbb I+\mathbb B\right)\cdot\bfchi_k=\bfcalf_k\,,
\eeq{1.21}
with $\mathbb I$ identity matrix.
The matrix $\i\,k\,\omega\,\mathbb I+\mathbb B$ is invertible. In fact, using \eqref{B}, for all $\bfxi\in\mathbb C^3$, we obtain the relation
$$
\bfxi^*\cdot(\i\,k\,\omega\,\mathbb I+\mathbb B)\cdot\bfxi=\i\,k\,\omega\left(|\bfxi|^2+
\|\xi_{i}\bfH_k^{(i)}\|_2^2\right)+\|\xi_{i}\bfH_k^{(i)}\|_2^2\,,
$$
that shows that $0$ is not an eigenvalue. As a result, for the given $\bfcalf_k$,  \eqref{1.21} has one and only one solution $\bfchi_k$.
If we dot-multiply both sides of \eqref{1.21} by $\bfchi_k^*$ and use \eqref{B} we deduce
$$
\i\,k\,\omega\left(M\,|\bfchi_k|^2+
\|\chi_{ki}\bfH_k^{(i)}\|_2^2\right)+\|\chi_{ki}\bfH_k^{(i)}\|_2^2=(\bfcalf_k,\bfchi_k^*)\,,
$$
that entails, in particular, the estimate
$$ 
|\bfchi_k|\le \frac1{|k|\omega}|\bfcalf_k|\,,\ \ |k|\ge1\,, 
$$
from which we conclude
\be
|\bfchi|^2_{W^{1,2}(0,T)}=\sum_{|k|\ge1}(|k|^2+1)|\bfchi_k|^2\le c\,\|\bfcalf\|_{L^2(0,T)}^2\,.
\eeq{1.46}
Combining \eqref{1.20}, \eqref{1.46} and \eqref{0.0.0} we thus infer
\be
\sum_{|k|\ge1}\left((|k|^2+1)\|\bfu_k\|_2^2+\|\nabla\bfu_k\|_2^2+|\bfu_k|_{2,2}^2\right)\le c\sum_{|k|\ge1}(|k|^2+1)|\bfchi_k|^2\le c\,\|\bfcalf\|_{L^2(0,T)}^2\,.
\eeq{1.49}
From \eqref{chi}, \eqref{1.46},  \eqref{1.49} and it follows that the vector functions $\bfu$ and $\bfchi$ defined in \eqref{Fou} satisfy \eqref{2.4} and, in addition 
\be
(\bfu,{\sf q}, \bfchi)\in \calw_\sharp^2\times\calp^{1,2}\times W^{1,2}_\sharp(0,T)\,.
\eeq{bolo} 
We next  extend $\bfchi$   to a solenoidal function  $\bfh\in \hat{\calw}_\sharp^{2,q}$ \cite[Section III.3]{Gab}  and write the solution $\bfu$ to \eqref{2.4}$_{1,2,3}$ as $\bfu=\bfv+\bfh$. From  \cite[Theorem 12]{GaMaH} and \eqref{bolo} it then follows that $\bfu\in \calw_\sharp^{2,q}$ and  ${\sf q}\in \calp^{1,q}$. We may thus conclude that $(\bfw:=\bfz+\bfu,\sfp:=\tau+{\sf q},\bfchi_\sub{\bfw}:=\bfchi)$ is a solution to \eqref{2.00} in the class $\calw_\sharp^{2,q}\times \calp^{1,q}\times W_\sharp^{1,2}(0,T)$, that satisfies in addition \eqref{n}. This completes the proof of the existence property when $q\in (1,2]$. We shall next show that, if $\bfF\in L^q_\sharp(0,T)$, with $q\in (2,\infty)$, the solution $(\bfw,{\sf p},\bfchi_\sub{\bfw})$ just constructed satisfies the other properties stated in the existence part of the lemma. Actually, in view of \eqref{2.5} and \eqref{2.6}, it is enough to show that if $\bfcalf\in L_\sharp^q(0,T)$, the solution $(\bfu,\bfchi)$ to \eqref{2.4} that we proved to be in the class $\calw_\sharp^2\times W^{1,2}_\sharp(0,T)$ is, in fact, in $\calw_\sharp^{2,q}\times W^{1,q}_\sharp(0,T)$ and that $(\bfw,{\sf p},\bfchi_\sub{\bfw})$  satisfies the estimate
given in \eqref{n} for $q\in(2,\infty)$.
In order to reach our goal, we recall the following inequalities:
\be
\|v\|_{r,\partial\Omega}\le c\,\left(\| v\|_{r,\Omega_R}+\| v\|_{r,\Omega_R}^{\frac1{r'}}| v|_{1,r,\Omega_R}^{\frac1r}\right)\,,\ \ r\in (1,\infty)\,, \ \ v\in W^{1,2}(\Omega_R)\,,
\eeq{in1}
and
\be
\|{\sf q}\|_{r,\Omega_R}\le  c\,\left(\|\nabla \bfu\|_{r,\Omega_R}+\|\nabla \bfu\|_{r,\Omega_R}^{\frac1{r'}}| \bfu|_{2,r,\Omega_R}^{\frac1r}\right)\,,\ \ r\in (1,\infty)\,,
\eeq{in2}
where $(\bfu,{\sf q})$ is the solution to \eqref{2.4} constructed previously. The first is a well known trace inequality (e.g. \cite[Theorem II.4.1]{Gab}), whereas the second one is proved in \cite[Lemma 2.5]{GaKy}. We now employ \eqref{in1} with $v\equiv\nabla\bfu$ and $r=2$. Since $\calw_\sharp^2\subset C([0,T];W^{1,2}(\Omega)$ and $\bfu\in \calw_\sharp^2$,  we obtain
\be
\nabla\bfu\in L^4(L^2(\partial\Omega))\,.
\eeq{sam1}
Likewise, from \eqref{in2}, we deduce
$$
{\sf q}\in L^4(L^2(\Omega_R))\,,
$$
and so, applying \eqref{in1} with $v\equiv{\sf q}$ and $r=2$, it follows that 
\be
{\sf q}\in L^4(L^2(\partial\Omega))\,.
\eeq{sam2}
Therefore, from \eqref{sam1}, \eqref{sam2} and \eqref{2.4}$_4$, we find $\bfchi\in W^{1,q}_\sharp(0,T)$, for all $q\in (2,4]$, provided $\bfF\in L^q_\sharp(0,T)$ and $\bff\in \call_\sharp^{2,q}$. Extending the boundary data $\bfchi$ to the smooth function $\bfh$ as done previously, by \cite[Theorem 12]{GaMaH} we get that $(\bfu,{\sf q})\in \calw_\sharp^{2,q}\times\calp^{1,q}$ as well. With this improved regularity, by an argument entirely analogous to that used before we show that 
$$  
\nabla\bfu,{\sf q}\in L^q(L^2(\partial\Omega))\,,\ \ q\in (4,8]
$$
which gives $(\bfu,{\sf q})\in \calw_\sharp^{2,q}\times\calp^{1,q}$, provided $\bfF\in L^q_\sharp(0,T)$ and $\bff\in \call_\sharp^{2,q}$. With the help of a simple boot-strap procedure, we then show that $(\bfu,{\sf q})\in \calw_\sharp^{2,q}\times\calp^{1,q}$, provided $\bfF\in L^q_\sharp(0,T)$ and $\bff\in \call_\sharp^{2,q}$, for all $q\in (2,\infty)$. Therefore, we conclude that for any $(\bff,\bfF)\in \call_\sharp^{2,q}\times L^q_\sharp(0,T)$, $q\in (2,\infty)$ there exists at least one corresponding solution $\big(\bfw,\sfp,\bfchi_\sub{\bfw}\big)\in \mathcal W_\sharp^{2,q}
\times \mathcal P^{1,q}\times W^{1,q}_\sharp(0,T)$ to problem \eqref{2.00}. Therefore, the validity of the estimate \eqref{n} for $q\in(2,\infty)$ is a consequence of the Open Mapping Theorem, provided we prove uniqueness in the class $\mathcal W_\sharp^{2,q}
\times \mathcal P^{1,q}\times W^{1,s}_\sharp(0,T)$. 
The latter amounts to show that the problem
\be\ba{cc}\medskip\left.\ba{ll}\medskip
\partial_t\bfw-\Delta\bfw=-\nabla {\sfp}\\ 
\Div\bfw=0\ea\right\}\ \ \mbox{in $\Omega\times\real$}\\ \medskip
\bfw|_{\partial\Omega}=\bfchi_\sub{\bfw}\,;\\ \dot{\bfchi}_\sub{\bfw}+\Int{\partial\Omega}{}\bfT(\bfw,{\sf p})\cdot\bfn=\0
\ea
\eeq{1.51}
has only the zero solution in the above function class. This is easily established. In fact, if we dot-multiply \eqref{1.51}$_1$ by $\bfw$, integrate by parts over $\Omega$ and use \eqref{1.51}$_3$, we get
$$
\half\ode{}t(\|\bfw(t)\|_2^2+|\bfchi_\sub{\bfw}(t)|^2)+\|\mathbb D(\bfw(t))\|_2^2=0\,.
$$
Integrating both sides of this equation from $0$ to $T$ and employing the $T$-periodicity leads us to  $\|\mathbb D(\bfw(t))\|_2\equiv 0$ which, in turn, by the characterization of the space $\calh$ given in \lemmref{1}, immediately furnishes $\bfw\equiv\nabla\sfp\equiv\0$. The proof of the lemma is completed.\par\hfill$
\square$

\setcounter{equation}{0}
\section{On the Strong Solvability of the Nonlinear Problem}
The  main objective of this section is  to show that, if the data $\bfv_*$ are sufficiently small (in certain  norms), then the problem \eqref{SE} possesses a unique ``strong" solution in a suitable function class. In order to achieve this goal, we need some further preparatory results. 
\smallskip\par
Set
\be\ba{ll}\medskip
\calv^{q}_\sharp:=\{\bfv\in W^{1,q}(0,T;W^{2-\frac1q,q}(\partial\Omega)):\ \bfv \ \mbox{is $T$-periodic with $\bar{\bfv}=\0$}\}\,,\\
\calv^{2,q}_\sharp:=\calv^{2}_\sharp\cap \calv^{q}_\sharp
\ea
\eeq{OP}
and
$$
\|\bfv\|_{\calv^{2,q}_\sharp}:=
\left(\int_0^T\big[\|\bfv\|_{2,\frac32 (\partial\Omega)}^2+\|\partial_t\bfv\|_{2,\frac32 (\partial\Omega)}^2\big]\right)^{\frac12}
+
\left(\int_0^T\big[\|\bfv\|_{q,2-\frac1q (\partial\Omega)}^q+\|\partial_t\bfv\|_{q,2-\frac1q (\partial\Omega)}^q\big]\right)^{\frac1q}\,.
$$
The following lemma holds.
\Bl Suppose $\bfv_*\in \calv^{2,q}_\sharp$, $q\in(1,\infty)$. Then, there exists a solenoidal extension $\bfV$ of $\bfv_*$ to $\Omega$ such that
$\bfV\in \hat{\calw}_\sharp^{2,q}$,  and 
$$
\|\bfV\|_{ \hat{\calw}_\sharp^{2,q}}\le c\, \|\bfv_*\|_{ \calv_\sharp^{2,q}}\,. 
$$
\EL{4.1}
{\em Proof.} We  construct the extension $\bfV$ in a way similar to that given in \cite[Theorem 1]{GaTP}. For fixed $t\in [0,T]$ consider the following boundary-value problem
\be\ba{cc}\medskip\left.\ba{ll}\medskip
\bfV =\nu\Delta \bfV-\nabla \phi\\
\Div \bfV =0\,
\ea\right\}\ \ \mbox{in $\Omega$\,,}
\\
\qquad\, \bfV=\bfv_*(t)\ \ \mbox{on $ \partial\Omega$}\,.
\ea
\eeq{89}
By the results of \cite[\S VII.5]{Gab}, we know that, for each $t\in [0,T]$ there exists a unique solution $(\bfV(t),\phi(t))$ to \eqref{89} such that
\be 
\bfV(t)\in W^{2,q}(\Omega)\,,\ \ \phi(t)\in D^{1,q}(\Omega)\,,  \ \ \mbox{all $q\in(1,\infty)$}
\eeq{90}
which satisfies, in addition,
\be
\|\bfV(t)\|_{2,q}+\|\nabla\phi(t)\|_q\le c_1\|\bfv_*(t)\|_{q,2-\frac1q (\partial\Omega)}
\,.
\eeq{91}
In view of the uniqueness property and the periodicity of $\bfv_*$ it follows that $\bfV$ is $T$-periodic  with $\bar \bfV\equiv 0$. Furthermore, again by uniqueness,  the regularity assumptions on $\bfv_*$ and \eqref{90} and \eqref{91}, we easily show that $\bfV$ is time-differentiable in the sense of distribution with $ \partial_t\bfV$ satisfying
$$
\ba{cc}\medskip\left.\ba{ll}\medskip
\partial_t\bfV =\nu\Delta \partial_t\bfV-\nabla \partial_t\phi\\
\Div \partial_t\bfV =0\,
\ea\right\}\ \ \mbox{in $\Omega$\,,}
\\
\qquad\, \partial_t\bfV=\partial_t\bfv_*\ \ \mbox{on $ \partial\Omega$}\,,
\ea
$$
from which it follows  $\partial_t\bfV\in L^q(\Omega)$ and
$$
\|\partial_t\bfV\|_q\le c\,\|\partial_t\bfv_*\|_{{q,1-\frac1q}(\partial\Omega)}\,.
$$
The proof of the lemma is thus completed.
\par\hfill$\square$\par\smallskip
Next,  consider the following linear problem in the unknowns  $(\bfV_0,P_0,\bfchi_0)$:
\be\ba{cc}\medskip\left.\ba{ll}\medskip
\partial_t\bfV_0=\nu\Delta\bfV_0-\nabla P_0\\
\Div\bfV_0=0\ea\right\}\ \ \mbox{in $\Omega\times\real$}\\ \medskip
\bfV_0(x,t)=\bfV_*(x,t)+\bfchi_0(t)\ \ \mbox{at $\partial\Omega\times\real$}\\
M\,\dot{\bfchi}_0 +\Int{\partial\Omega}{}\mathbb T(\bfV_0,P_0)\cdot\bfn=\0
\ea
\eeq{6.1}
\Bl Suppose $\bfV_*\in \calv^{2,q}_\sharp$, $q\in(1,\infty)$. Then, problem \eqref{6.1} has one and only one solution $(\bfV_0,P_0,\bfchi_0)\in \hat{\calw}_\sharp^{2,q}\times\calp^{1,q}\times (W^{1,2}_\sharp(0,T)\cap W^{1,q}_\sharp(0,T))$.  This solution satisfies the estimate:
\be
\|\bfV_0\|_{\hat{\calw}_\sharp^{2,q}}+\|P_0\|_{\calp^{1,q}}+\|\bfchi_0\|_{W^{1,s}(0,T)}\le c\,\|\bfV_*\|_{\calv^{2,q}_\sharp}\,,
\eeq{6.2}
where $s=2$ if $q\le2$, and $s=q$ otherwise.
\EL{6.1}
{\em Proof.} We look for a solution of the form $\bfV_0=\bfV_1+\bfV_2$,  where $\bfV_2$ is the extension of $\bfV_*$ constructed in \lemmref{4.1} and
\be\ba{cc}\medskip\left.\ba{ll}\medskip
\partial_t\bfV_1=\nu\Delta\bfV_1-\nabla P_0+\bfh\\
\Div\bfV_1=0\ea\right\}\ \ \mbox{in $\Omega\times\real$}\\ \medskip
\bfV_1(x,t)=\bfchi_0(t)\ \ \mbox{at $\partial\Omega\times\real$}\\
M\,\dot{\bfchi}_0 +\Int{\partial\Omega}{}\mathbb T(\bfV_1,P)\cdot\bfn=-\nu\Int{\partial\Omega}{}\mathbb D(\bfV_2)\cdot\bfn\,.
\ea
\eeq{6.3}
where
$$
\bfh=-\partial_t\bfV_2+\nu\Delta\bfV_2\,.
$$
From \lemmref{4.1} we have
\be
\|\bfh\|_{\call^{2,q}}\le c\,\|\bfV_*\|_{\calv^{2,q}_\sharp}\,,
\eeq{6.4}
and by the trace theorem and \lemmref{4.1},
\be
\|\Int{\partial\Omega}{}\mathbb D(\bfV_2)\cdot\bfn\|_{L^q(L^q)}\le c\,\|\bfV_2\|_{\hat{\mathcal W}_\sharp^{2,q}}\le c\,\|\bfV_*\|_{\calv^{2,q}_\sharp}\,.
\eeq{6.5}
Therefore, the result follows by applying \lemmref{1.6_0} to \eqref{6.3}, and using \eqref{6.4}--\eqref{6.5}.
\par\hfill$\square$\par\medskip
\par
For $\delta>0$, we set 
\be\hat{\bfV}=\delta\bfV_0\,,\ \ \hat{P}=\delta P_0\,,\ \ \hat{\bfchi}=\delta\bfchi_0\,,\ \  \bfv_*=\delta\bfV_*
\,,
\eeq{6.6}
and  look for a solution to \eqref{SE} --with $\bfv_*$ as in \eqref{6.6}--  
``around" $(\hat{\bfV},\hat{P}, \hat{\bfchi})$, namely
\be\ba{ll}\medskip
\bfv=(\bfv-\bar{\bfv})+\bar{\bfv}:=\bfw+\hat{\bfV}+\bfu\,;\ \ \ p=(p-\bar{p})+\bar{p}:=\tau+\hat{P}+\sfp\,;\\
\bfgamma=(\bfgamma-\bar{\bfgamma})+\bar{\bfgamma}:=\bfchi+\hat{\bfchi}+\bfxi\,;\ \ \bar{\bfw}\equiv\bar{\bfchi}\equiv\0\,,\ \ \bar{\tau}\equiv\0.
\ea
\eeq{ABP}
From \eqref{SE}  it immediately follows that $(\bfu,\sfp)$ and $(\bfw,\tau)$ solve the following coupled  system of nonlinear equations:
\be\ba{cc}\medskip\left.\ba{ll}\medskip
\nu\Delta\bfu+\bfxi_\sub{\bfu}\cdot\nabla\bfu-\nabla \sfp=\Div\bfscr{F}(\bfu,\bfw,\bfchi_\sub{\bfw})\\ \Div\bfu=0\ea\right\}\ \ \mbox{in $\Omega$}\,,
\\ \medskip
\bfu=\bfxi_\sub{\bfu}\ \ \mbox{at $\partial\Omega$}\\

\Int{\partial\Omega}{}[\mathbb T(\bfu,\sfp)-{\bfscr{F}(\bfu,\bfw,\bfchi_\sub{\bfw})}]\cdot\bfn=\0\,,
\ea
\eeq{6.10}
and
\be\ba{cc}\medskip\left.\ba{ll}\medskip
\partial_t\bfw-\nu\Delta\bfw=-\nabla {\tau}+\bff(\bfu,\bfw,\bfchi_\sub{\bfw})\\ 
\Div\bfw=0\ea\right\}\ \ \mbox{in $\Omega\times\real$}\,,\\ 
\medskip
\bfw=\bfchi_\sub{\bfw}\ \ \mbox{at $\partial\Omega\times\real$}\\
M\,\dot{\bfchi}_\sub{\bfw}+\Int{\partial\Omega}{}\bfT(\bfw,{\tau})\cdot\bfn=\bfF(\bfchi_\sub{\bfw},\bfxi_\sub{\bfu})
\ea
\eeq{6.11}
with
\be
\bfscr{\bfscr F}(\bfu,\bfw,\bfchi_\sub{\bfw}):=\bfu\otimes\bfu+\bar{(\bfw-\bfchi_\sub{\bfw})\otimes\bfw}+\bar{(\hat{\bfV}-\hat{\bfchi})\otimes\bfw}+\bar{\bfw\otimes\hat{\bfV}}+\bar{(\hat{\bfV}-\bfchi_\sub{\bfw}-\hat{\bfchi})\otimes\hat{\bfV}}\,,
\eeq{6.12}
and
\be\ba{ll}\medskip
\bff(\bfu,\bfw,\bfchi_\sub{\bfw}):=
-(\bfw-\bfchi_\sub{\bfw})\cdot\nabla\bfw-(\hat{\bfV}-\hat{\bfchi})\cdot\nabla\bfw-\bfw\cdot\nabla\hat{\bfV}+\bar{(\bfw-\bfchi_\sub{\bfw})\cdot\nabla\bfw}+\bar{(\hat{\bfV}-\hat{\bfchi})\cdot\nabla\bfw}\\\medskip\qquad+\bar{\bfw\cdot\nabla\hat{\bfV}}
+\bfchi_\sub{\bfw}\cdot\nabla\hat{\bfV}-\bar{\bfchi_\sub{\bfw}\cdot\nabla\hat{\bfV}}
- \bfu\cdot\nabla(\bfw+\hat{\bfV})
-(\bfw-\bfchi_\sub{\bfw})\cdot\nabla\bfu-(\hat{\bfV}-\hat{\bfchi})\cdot\nabla\bfu\\ \medskip\hspace*{4cm}-(\hat{\bfV}-\hat{\bfchi})\cdot\nabla\hat{\bfV}+\bar{(\hat{\bfV}-\hat{\bfchi})\cdot\nabla\hat{\bfV}}\\
\bfF(\bfchi_\sub{\bfw},\bfxi_\sub{\bfu}):= \Int{\partial\Omega}{}[(\bfv_*+\bfxi_\sub{\bfu}+\bfchi_\sub{\bfw}+\hat{\bfchi})\otimes\bfv_*-\bar{(\bfv_*+\bfchi_\sub{\bfw}+\hat{\bfchi})\otimes\bfv_*}]\cdot\bfn
\ea
\eeq{6.13}
Existence and uniqueness of a solution to \eqref{6.10}--\eqref{6.13} for sufficiently small $\delta>0$  will be shown  by a suitable successive approximation scheme. To this end, it is convenient to define first another function class. For $q\in (1,\infty)$, set 
$$
\mathfrak B^q:=\left\{(\bfphi,\bfpsi,\bfchi_\sub{\bfpsi})\in{\sf D}^{2,q}\times \calw_\sharp^{2,q}\times [W^{1,2}_\sharp(0,T)\cap W^{1,q}_\sharp(0,T)]\right\}\,. 
$$
Plainly, $\mathfrak B^q$ becomes a Banach space endowed with the norm
$$
\|(\bfphi,\bfpsi,\bfchi_\sub{\bfpsi})\|_{\mathfrak B^q}:=\|\bfphi\|_{\sfD^{2,q}}+\|\bfpsi\|_{\calw^{2,q}}+\|\bfchi_\sub{\bfpsi}\|_{W^{1,2}(0,T)\cap W^{1,q}_\sharp(0,T)}\,.
$$
Moreover, also with the help of \cite[Exercise II.6.2]{Gab}, we show that $\mathfrak B^q$ is a reflexive space. 

Next, define
\be\ba{ll}\medskip
\bfscr F_0:=\bar{(\bfV_0-\bfchi_0)\otimes\bfV_0}\,,\ \ \ \bff_0:=({\bfV_0-\bfchi_0)\otimes\bfV_0} -\bar{(\bfV_0-\bfchi_0)\otimes\bfV_0}\\
\bfF_0:= \Int{\partial\Omega}{}{[(\bfV_*+{\bfchi}_0)\otimes\bfV_*}-\bar{(\bfV_*+{\bfchi}_0)\otimes\bfV_*}]\cdot\bfn\,.
\ea
\eeq{F_0}
and
\be\ba{ll}\medskip
\tilde{\bfscr F}(\bfu,\bfw,\bfchi_{\sub{\bfw}}):= {\bfscr F}(\bfu,\bfw,\bfchi_{\sub{\bfw}})-\delta^2{\bfscr F}_0\,;\,\ \ \tilde{\bff}(\bfu,\bfw,\bfchi_{\sub{\bfw}}):={\bff}(\bfu,\bfw,\bfchi_\sub{\bfw})-\delta^2\bff_0\,,\\ 
\tilde{\bfF}(\bfchi_\sub{\bfw},\bfxi_{\sub{\bfu}}):=\bfF(\bfchi_\sub{\bfw},\bfxi_{\sub{\bfu}})-\delta^2\bfF_0\,.
\ea
\eeq{Ftilde}
\smallskip\par\noindent
The following lemma holds.
\Bl Assume $q\in(1,\mbox{$\frac32$}]$, and let $\bfU:=(\bfu,\bfw,\bfchi_\sub{\bfw})\in \mathfrak B^q$, and $\hat{\bfV}$ as in \eqref{6.6}. Then, there is a positive constant $c=c(\Omega,q)$ such that
\be\ba{ll}\medskip
\|\Div\tilde{\bfscr F}(\bfU)\|_{2}+\|\tilde{\bfscr F}(\bfU)\|_{2}+\|\tilde{\bfscr F}(\bfU)\|_{q}+\|\tilde{\bff}(\bfU)\|_{\call^{2,q}}\le c\Big(\|\bfU\|_{\mathfrak B^q}^2+\delta\|\bfU\|_{\mathfrak B^q}\Big)\\\medskip 
\|\tilde{\bfF}(\bfchi_\sub{\bfw},\bfxi_{\sub{\bfu}})\|_{L^2(0,T)\cap L^q(0,T)}\le
c\,\delta\,\|\bfU\|_{\mathfrak B^q}\,,
\\
\|\Div{\bfscr F}_0(\bfU)\|_{2}+\|{\bfscr F}_0(\bfU)\|_{2}+\|{\bfscr F}_0(\bfU)\|_{q}\le c\,.
\\
\ea
\eeq{4.11}
\EL{4.2}
{\em Proof.} By assumption, $\bfU\in \mathfrak B^q$, and, by  \lemmref{4.1}, $\bfV_0\in \calw_\sharp^{2,q}$. Thus, using  \eqref{WH2}, it follows at once $\tilde{\bfF}\in L^2_\sharp(0,T)\cap L^2_\sharp(0,T)$ and the validity of \eqref{4.11}$_2$. Moreover, with the help of \lemmref{Emb0} and H\"older inequality, we get
$$
\|\bfu\otimes\bfu\|_2+\|\bfu\otimes\bfu\|_q\le c\,\|\bfu\|_4^2+\|\bfu\|_{\frac{3q}{3-q}}\|\bfu\|_3\le c\,\|\bfu\|_{{\sf D}^{2,q}}^2\le c\,\|\bfU\|_{\mathfrak B^q}^2\,,
$$
and
$$
\|\Div(\bfu\otimes\bfu)\|_2=\|\bfu\cdot\nabla\bfu\|_2\le \|\bfu\|_6\|\nabla\bfu\|_3\le c\,\|\bfu\|_{{\sf D}^{2,q}}^2\le c\,\|\bfU\|_{\mathfrak B^q}^2\,.
$$
In analogous fashion, using this time \lemmref{Mallo}, we get
$$\ba{ll}\medskip
\|\bar{(\bfw-\bfchi_\sub{\bfw})\otimes\bfw}\|_2+\|\bar{(\bfw-\bfchi_\sub{\bfw})\otimes\bfw}\|_q\le c\,\Big[\|\bfchi_\sub{\bfw}\|_{L^\infty(0,T)}\left(\|\bfw\|_{L^2(L^2)}+\|\bfw\|_{L^q(L^q)}\right)\\
\hspace*{2.7cm}+\|\bfw\|_{L^2(L^\infty)}\|\bfw\|_{L^\infty(L^2)}+\|\bfw\|_{L^q(L^\infty)}\|\bfw\|_{L^\infty(L^q)}\Big]\le c\,\|\bfU\|_{\mathfrak B^q}^2\,,
\ea
$$
and also
$$\ba{rl}\medskip
\|\Div\bar{(\bfw-\bfchi_\sub{\bfw})\otimes\bfw}\|_2=\|\bar{(\bfw-\bfchi_\sub{\bfw})\cdot\nabla\bfw}\|_2&\!\!\!\le c\,
\Big[\|\bfchi_\sub{\bfw}\|_{L^\infty(0,T)}\|\nabla\bfw\|_{L^2(L^2)}
+\|\bfw\|_{L^6(L^6)}\|\nabla\bfw\|_{L^3(L^3)}\Big]\\
&\!\!\!\le c\,\|\bfU\|_{\mathfrak B^q}^2\,.
\ea
$$
By an entirely similar procedure, we show
$$
\|\bar{(\hat{\bfV}-\hat{\bfchi})\otimes\bfw+ \bfw\otimes\hat{\bfV}-\bfchi_\sub{\bfw}\otimes\hat{\bfV}}\|_2+\|\bar{(\hat{\bfV}-\hat{\bfchi})\otimes\bfw+ \bfw\otimes\hat{\bfV}-\bfchi_\sub{\bfw}\otimes\hat{\bfV}}\|_q\le c\,\delta\,\|\bfU\|_{\mathfrak B^q}
$$
and
$$
\|\Div\bar{(\hat{\bfV}-\hat{\bfchi})\otimes\bfw+ \bfw\otimes\hat{\bfV}-\bfchi_\sub{\bfw}\otimes\hat{\bfV}}\|_2\le c\,\delta\,\|\bfU\|_{\mathfrak B^q}\,.
$$
The proof of the stated estimate for $\bfscr F_0$ and $\bff$ is very much alike, and we leave it to the reader.
\par\hfill$\square$\par
\ms\par We also have the following result.
\Bl Assume $q\in (1,\frac32)$, and let $$\bfU_i:=(\bfu_i,\bfw_i,\bfchi_\sub{\bfw_i})\in \mathfrak B^q\cap \mathfrak B^{\frac{3q}{3-q}}\,,\ \ i=1,2\,,$$   
with $\|\bfU_i\|_{\mathfrak B^q}\le\delta$, $i=1,2$,  $\delta\in (0,\delta_0]$, suitable $\delta_0>0$. Moreover, let $\hat{\bfV}\in \calw_\sharp^{2,q}\cap \calw_\sharp^{2,\frac{3q}{3-q}}$. Then, there exists a positive constant $c(\Omega,q)$ such that
\be\ba{ll}\medskip
\|\Div[\tilde{\bfscr F}(\bfU_1)-\tilde{\bfscr F}(\bfU_2)]\|_{2}+\|\tilde{\bfscr F}(\bfU_1)-\tilde{\bfscr F}(\bfU_2)\|_{2}+\|\tilde{\bfscr F}(\bfU_1)-\tilde{\bfscr F}(\bfU_2)\|_{\frac{3q}{3-q}}+\|\tilde{\bff}(\bfU_1)-\tilde{\bff}(\bfU_2)\|_{\call^{2,\frac{3q}{3-q}}}\\ \medskip
\hspace*{3.5cm}\le c\,\delta\,\|\bfU_1-\bfU_2\|_{\mathfrak B^{\frac{3q}{3-q}}}\\
\|\tilde{\bfF}(\bfchi_\sub{\bfw_1},\bfxi_{\sub{\bfu_1}})-\tilde{\bfF}(\bfchi_\sub{\bfw_2},\bfxi_{\sub{\bfu_2}})\|_{L^2(0,T)\cap L^{\frac{3q}{3-q}}(0,T)}\le
c\,\delta\,\|\bfU_1-\bfU_2\|_{\mathfrak B^{\frac{3q}{3-q}}}\,.
\ea
\eeq{4.111}
\EL{Cauc}
{\em Proof.} Again, under the given assumption, the proof of \eqref{4.111}$_2$ is straightforward. Next, setting for simplicity $\bfu:=\bfu_1-\bfu_2$, we get
$$
\|\bfu_1\otimes\bfu_1-\bfu_2\otimes\bfu_2\|_2+\|\bfu_1\otimes\bfu_1-\bfu_2\otimes\bfu_2\|_{\frac{2q}{3-q}}\le \|\bfu_1\otimes\bfu\|_2+\|\bfu\otimes \bfu_2\|_2+\|\bfu_1\otimes\bfu\|_{\frac{3q}{3-q}}+\|\bfu\otimes \bfu_2\|_{\frac{3q}{3-q}}\,.
$$
Observing that, under the given assumption,  we have $3q/(3-q)<3$, by the H\"older inequality and  \lemmref{Emb0} we show
$$
\|\bfu_1\otimes\bfu\|_2+\|\bfu_1\otimes\bfu\|_{\frac{3q}{3-q}}\le c\,(\|\bfu_1\|_3\|\bfu\|_6+\|\bfu_1\|_3\|\bfu\|_{\frac{3q}{3-2q}})\le c\,\delta\,\|\bfu\|_{{\sf D}^{2,\frac{3q}{3-q}}}\le c\,\delta\,\|\bfU_1-\bfU_2\|_{\mathfrak B^{\frac{3q}{3-q}}}\,,
$$
and, likewise,
$$
\|\bfu\otimes\bfu_2\|_2+\|\bfu\otimes\bfu_2\|_{\frac{3q}{3-q}}
\le c\,\delta\,\|\bfU_1-\bfU_2\|_{\mathfrak B^{\frac{3q}{3-q}}}\,,
$$
so that
$$
\|\bfu_1\otimes\bfu_1-\bfu_2\otimes\bfu_2\|_2+\|\bfu_1\otimes\bfu_1-\bfu_2\otimes\bfu_2\|_{\frac{2q}{3-q}}\le c\,\delta\,\|\bfU_1-\bfU_2\|_{\mathfrak B^{\frac{3q}{3-q}}}. 
$$
Moreover,
$$
\|\Div(\bfu_1\otimes\bfu_1-\bfu_2\otimes\bfu_2)\|_2\le\|\bfu_1\cdot\nabla \bfu\|_2+\|\bfu\cdot\nabla \bfu_2\|_2\le c\,\left(\|\bfu_1\|_\infty\|\nabla\bfu\|_2+\|\bfu\|_\infty\|\nabla\bfu_2\|_2\right)\,,
$$
which, with the help of \lemmref{Emb0}, furnishes
$$
\|\Div(\bfu_1\otimes\bfu_1-\bfu_2\otimes\bfu_2)\|_2\le c\, \delta\,\|\bfU_1-\bfU_2\|_{\mathfrak B^{\frac{3q}{3-q}}}\,.
$$
Next, setting $\bfw:=\bfw_1-\bfw_2$, $\bfchi:=\bfchi_1-\bfchi_2$, $\bfchi_i:=\bfchi_\sub{\bfw_i}$, $i=1,2$, by H\"older inequality we have
$$\ba{ll}\medskip
\|\bar{(\bfw_1-\bfchi_1)\otimes\bfw_1-(\bfw_2-\bfchi_2)\otimes\bfw_2}\|_2+\|\bar{(\bfw_1-\bfchi_1)\otimes\bfw_1-(\bfw_2-\bfchi_2)\otimes\bfw_2}\|_{\frac{3q}{3-q}}\\ \medskip\quad\quad
\le \|\bar{(\bfw_1-\bfchi_1)\otimes\bfw}\|_2+\|\bar{(\bfw-\bfchi)\otimes\bfw_2}\|_2+\|\bar{(\bfw_1-\bfchi_1)\otimes\bfw}\|_{\frac{3q}{3-q}}+\|\bar{(\bfw-\bfchi)\otimes\bfw_2}\|_{\frac{3q}{3-q}}
\\ \medskip\quad\quad
\le c\,\Big[\|{\bfw_1-\bfchi_1}\|_{L^2(L^\infty)}\|{\bfw}\|_{L^\infty(L^2)}+\|\bfw-\bfchi\|_{L^2(L^\infty)}\|\bfw_2\|_{L^\infty(L^2)}\\ \qquad\quad\qquad
+
\|{\bfw_1-\bfchi_1}\|_{L^{\frac{3q}{3-q}}(L^\infty)}\|{\bfw}\|_{L^{\infty}(L^{\frac{3q}{3-q}})}+\|\bfw-\bfchi\|_{L^{\frac{3q}{3-q}}(L^\infty)}\|\bfw_2\|_{L^{\infty}(L^{\frac{3q}{3-q}})}\,.
\Big]
\ea
$$
Thus, observing that $3q/(3-q)>3/2$, from \lemmref{Mallo} and the assumption we derive 
$$\ba{ll}\medskip
\|\bar{(\bfw_1-\bfchi_1)\otimes\bfw_1-(\bfw_2-\bfchi_2)\otimes\bfw_2}\|_2+\|\bar{(\bfw_1-\bfchi_1)\otimes\bfw_1-(\bfw_2-\bfchi_2)\otimes\bfw_2}\|_{\frac{3q}{3-q}}\\
\hspace*{4.5cm}\le c\, \delta\,\|\bfU_1-\bfU_2\|_{\mathfrak B^{\frac{3q}{3-q}}}\,.
\ea$$
Furthermore, with $r\in (2,4)$, we have
$$\ba{ll}\medskip
\|\Div(\bar{(\bfw_1-\bfchi_1)\otimes\bfw_1-(\bfw_2-\bfchi_2)\otimes\bfw_2})\|_2\le \|\bar{(\bfw_1-\bfchi_1)\cdot\nabla\bfw}\|_2+\|\bar{(\bfw-\bfchi)\cdot\nabla\bfw_2})\|_2\\
\hspace*{1.5cm}\le \|\bfw_1-\bfchi_1\|_{L^r(L^\infty)}\|\nabla\bfw\|_{L^{\frac{2r}{r-2}}(L^2)}+\|\bfw_1-\bfchi\|_{L^r(L^\infty)}\|\nabla\bfw_2\|_{L^{\frac{2r}{r-2}}(L^2)}\,,
\ea
$$
which, again by \lemmref{Mallo} and the assumption, furnishes 
$$\ba{rl}\medskip
\|\Div(\bar{(\bfw_1-\bfchi_1)\otimes\bfw_1-(\bfw_2-\bfchi_2)\otimes\bfw_2})\|_2&\!\!\!\le \|\bar{(\bfw_1-\bfchi_1)\cdot\nabla\bfw}\|_2+\|\bar{(\bfw-\bfchi)\cdot\nabla\bfw_2})\|_2\\&\!\!\!\le c\, \delta\,\|\bfU_1-\bfU_2\|_{\mathfrak B^{\frac{3q}{3-q}}}\,.
\ea
$$
In an entirely (and simpler) fashion one shows
$$
\|\bar{(\hat{\bfV}-\hat{\bfchi})\otimes\bfw+\bfw\otimes\hat{\bfV}-\bfchi\otimes\hat{\bfV}}\|_2+\|\bar{(\hat{\bfV}-\hat{\bfchi})\otimes\bfw+\bfw\otimes\hat{\bfV}-\bfchi\otimes\hat{\bfV}}\|_{\frac{3q}{3-q}}\le c\,\delta\,\|\bfU\|_{\mathfrak B^{\frac{3q}{3-q}}}
$$
and
$$
\|\Div\bar{(\hat{\bfV}-\hat{\bfchi})\otimes\bfw+ \bfw\otimes\hat{\bfV}-\bfchi_\sub{\bfw}\otimes\hat{\bfV}}\|_2\le c\,\delta\,\|\bfU\|_{\mathfrak B^q}\,.
$$
Finally, the proof of the stated estimate of $\bff$ is performed by exactly the same token, and we leave it to the reader.
\par\hfill$\square$\ms\par\par
We are now in a position to show the following existence and uniqueness result.
\Bt Let $\bfv_*$ be as in \eqref{6.6}, with $\bfV_*\in\calv_\sharp^{2,\frac{3q}{3-q}}$,  $q\in (1,\frac65]$. Then, there is $\delta_0>0$ such that if
$
\delta\le \delta_0,
$
the problem  \eqref{6.10}--\eqref{6.13} has one and only one solution $(\bfu,\bfw,\bfchi_\sub{\bfw})\in \mathfrak B^q\cap\mathfrak B^{\frac{3q}{3-q}}$ satisfying
\be
\|(\bfu,\bfw,\bfchi_\sub{\bfw})\|_{\mathfrak B^q}+\|(\bfu,\bfw,\bfchi_\sub{\bfw})\|_{\mathfrak B^{\frac{3q}{3-q}}}\le c\,\delta^2\,.
\eeq{gatta0}
\ET{6.1}
{\em Proof.} We shall employ an  approximating procedure in suitable spaces. Precisely, let's define a sequence $\{\bfu_{n},\bfxi_n,\bfw_{n},\bfchi_n\}$\footnote{For simplicity, we set $\bfxi_\sub{\bfu_n}\equiv\bfxi_n$, $\bfchi_\sub{\bfw_n}\equiv\bfchi_n$.} by recurrence, as follows:
\be\ba{cc}\medskip\left.\ba{ll}\medskip
\nu\Delta\bfu_n+\bfxi_n\cdot\nabla\bfu_n-\nabla \sfp_n=\Div\bfscr{F}(\bfu_{n-1},\bfw_{n-1},\bfchi_{n-1})\\ \Div\bfu_n=0\ea\right\}\ \ \mbox{in $\Omega$}\,,
\\ \medskip
\bfu_n=\bfxi_{n}\ \ \mbox{at $\partial\Omega$}\\

\Int{\partial\Omega}{}[\mathbb T(\bfu_n,\sfp_n)-{\bfscr{F}(\bfu_{n-1},\bfw_{n-1},\bfchi_{n-1})}]\cdot\bfn=\0\,,
\ea
\eeq{6.10_0}
and
\be\ba{cc}\medskip\left.\ba{ll}\medskip
\partial_t\bfw_n-\nu\Delta\bfw_n=-\nabla {\tau}_n+\bff(\bfu_{n-1},\bfw_{n-1},\bfchi_{n-1})\\ 
\Div\bfw_n=0\ea\right\}\ \ \mbox{in $\Omega\times\real$}\,,\\ 
\medskip
\bfw_n=\bfchi_n\ \ \mbox{at $\partial\Omega\times\real$}\\
M\,\dot{\bfchi}_n+\Int{\partial\Omega}{}\bfT(\bfw_n,{\tau}_n)\cdot\bfn=\bfF(\bfchi_{n-1},\bfxi_{n-1})\,,
\ea
\eeq{6.11_0}
where $n\ge1$ and $\bfu_0\equiv\bfw_0\equiv\0$. Our first objective is to show that the sequence is bounded in $\mathfrak B^q\cap \mathfrak B^{\frac{3q}{3-q}}$, provided $\delta$ is sufficiently small. From \eqref{6.6}, \eqref{6.12} and \eqref{6.13} we get
\be\ba{cc}\medskip\left.\ba{ll}\medskip
\nu\Delta\bfu_1+\bfxi_1\cdot\nabla\bfu_1-\nabla \sfp_1=\delta^2\Div(\bar{\bfV_0-\bfchi_0)\otimes\bfV_0}\\ \Div\bfu_1=0\ea\right\}\ \ \mbox{in $\Omega$}\,,
\\ \medskip
\bfu_1=\bfxi_{1}\ \ \mbox{at $\partial\Omega$}\\

\Int{\partial\Omega}{}[\mathbb T(\bfu_1,\sfp_1)-\delta^2\bar{(\bfV_*+{\bfchi}_0)\otimes\bfV_*}]\cdot\bfn=\0\,,
\ea
\eeq{6.10_1}
and
\be\ba{cc}\medskip\left.\ba{ll}\medskip
\partial_t\bfw_1-\nu\Delta\bfw_1=-\nabla {\tau}_1-\delta^2\left({(\bfV_0-\bfchi_0)\cdot\bfV_0}-\bar{(\bfV_0-\bfchi_0)\cdot\bfV_0}\right)\\ 
\Div\bfw_1=0\ea\right\}\ \ \mbox{in $\Omega\times\real$}\,,\\ 
\medskip
\bfw_1=\bfchi_1\ \ \mbox{at $\partial\Omega\times\real$}\\
M\,\dot{\bfchi}_1+\Int{\partial\Omega}{}\bfT(\bfw_1,{\tau}_1)\cdot\bfn=\delta^2\Int{\partial\Omega}{}{[(\bfV_*+{\bfchi}_0)\otimes\bfV_*}-\bar{(\bfV_*+{\bfchi}_0)\otimes\bfV_*}]\cdot\bfn\,.
\ea
\eeq{6.11_1}
In view of the regularity assumption made on $\bfV_*$,\footnote{Notice that $\calv_\sharp^{2,\frac{3q}{3-q}}\subset \calv_\sharp^{2,q}$.} by \lemmref{6.1} we see that all the terms involving $ \bfV_0$, and $\bfchi_0$ meet the hypotheses of \lemmref{Fm} and \lemmref{1.6_0}. Therefore, we conclude that for sufficiently small $\delta>0$ there exists a unique solution  such that
$$
(\bfu_1,\bfw_1,\bfchi_1)\in \mathfrak B^q\cap \mathfrak B^{\frac{3q}{3-q}}\,.
$$
For fixed $q$, let $M$ be any positive constant satisfying 
$$
M\ge \|\Div\bfscr F_0\|_q+\|\bfscr F_0\|_2+\|\bfscr F_0\|_2^2+\|\bfscr F_0\|_q+\|\bff_0\|_{\call^{2,q}_\sharp}+\|\bfF_0\|_{L^2_\sharp(0,T)\cap L^q_\sharp(0,T)}
$$
with $\bfscr F_0,\bff_0$, and $\bfF_0$ given in \eqref{F_0}, 
and let  $C_0=\max\{C_1,C_2\}$, with $C_1,C_2$ the constants defined in \eqref{K0} and \eqref{n}. Then, again from  \lemmref{Fm} and \lemmref{1.6_0} we deduce
\be
\|(\bfu_1,\bfw_1,\bfchi_1)\|_{\mathfrak B^{\frac{3q}{3-q}}}+\|(\bfu_1,\bfw_1,\bfchi_1)\|_\sub{\mathfrak B^q}\le \kappa_0\, \delta^2\,,
\eeq{gatta}
with $\kappa_0:=M\,C_0$, 
provided $\delta$ is small enough. In view of this result, we can employ a classical induction argument to show that 
$$
(\bfu_n,\bfw_n,\bfchi_n)\in \mathfrak B^q\,, \ \ \mbox{for all $n\ge1$},
$$
and that
\be
\|(\bfu_n,\bfw_n,\bfchi_n)\|_{\mathfrak B^q}\le 2\kappa_0\, \delta^2\,, \ \ \mbox{for all $n\ge1$}\,.
\eeq{step1}
Thus, let us assume that \eqref{step1} holds for $n-1$, and show that it is true also for $n$. To this end, set
$$
\bfU_n:=(\bfu_n,\bfw_n,\bfchi_n)\,,\ \ \ n\ge1\,,
$$
and
$$\ba{ll}\medskip
{\bfscr F}_{n-1}:= {\bfscr F}(\bfU_{n-1})-{\bfscr F}_0\,;\,\ \ \bff_{n-1}:={\bff}(\bfU_{n-1})-\bff_0\,,\\ 
\bfF_{n-1}:=\bfF(\bfchi_{n-1},\bfxi_{n-1})-\bfF_0\,.
\ea
$$
From \lemmref{4.2} and \eqref{step1} evaluated at $n-1$ it follows that
$$
\ba{ll}\medskip
\|\Div{\bfscr F}_{n-1}\|_{2}+\|{\bfscr F}_{n-1}\|_{2}+\|\bfscr F_{n-1}\|_{q}+\|{\bff}_{n-1}\|_{\call^{2,q}}\le c\Big(\|\bfU_{n-1}\|_{\mathfrak B^q}^2+\delta\|\bfU_{n-1}\|_{\mathfrak B^q}\Big)\le c_1\,\delta^3\,,\\
\|{\bfF}_{n-1}\|_{L^2(0,T)\cap L^q(0,T)}\le
c\,\delta\,\|\bfU_{n-1}\|_{\mathfrak B^q}\le c_2\,\delta^3\,,
\ea
$$
for $\delta \le1$. Thus, applying \lemmref{Fm} and \lemmref{1.6_0} to \eqref{6.10_0}, \eqref{6.11_0}, with the help of the above estimates we deduce for small $\delta$
\be
\|\bfU_n\|_{\mathfrak B^q}\le c_0\delta^3+\kappa_0\,\delta^2\,,
\eeq{Un}
which, in turn, again by taking $\delta$ below a certain positive constant, implies \eqref{gatta}. We shall next show that the sequence $\{\bfU_n\}$ is Cauchy in $\mathfrak B^{\frac{3q}{3-q}}$ by proving that it satisfies the following inequality
\be
\|\bfU_{n+1}-\bfU_n\|_{\mathfrak B^{\frac{3q}{3-q}}}\le C\,\delta \|\bfU_{n}-\bfU_{n-1}\|_{\mathfrak B^{\frac{3q}{3-q}}}\,,\ \ \ \mbox{all $n\ge 1$}\,,
\eeq{gatta1}
and then choosing $\delta<C^{-1}$. However, \eqref{gatta1} is an immediate consequence of the second part of \lemmref{Fm},  \lemmref{1.6_0}, and \lemmref{Cauc}. Thus, denoting by $\bfU:=(\bfu,\bfw,\bfchi_\sub{\bfw})\in \mathfrak B^{\frac{3q}{3-q}}$ the limiting point of the sequence, we obtain that $\bfU$ is the solution to the  problem \eqref{6.10}--\eqref{6.12} in the above class, for sufficiently small $\delta>0$. Furthermore, from \eqref{gatta1} and recalling that $\bfU_0\equiv\0$, we infer, again by standard arguments, that
$$
\|\bfU_{n+k}-\bfU_n\|_{\mathfrak B^{\frac{3q}{3-q}}}\le \frac{(C\delta)^n}{1-C\delta} \|\bfU_{n}\|_{\mathfrak B^{\frac{3q}{3-q}}}\,,\ \ \ \mbox{all $n\ge 1$, $k\ge 0$}.
$$
Thus, letting $k\to\infty$ in this relation and choosing $n=1$, with the help of \eqref{gatta} we arrive at \eqref{gatta0}. Finally, since $\mathfrak B^q$ is reflexive, from \eqref{Un} we deduce that there is $\bfU^*\in \mathfrak B^q$ such that $\bfU_n\to\bfU^*$ weakly in $\mathfrak B^q$, along a subsequence, with $\bfU^*$ satisfying \eqref{Un}.  As a result, we show that it must be  $\bfU^*=\bfU$, which completes the proof of the theorem.
\par\hfill$\square$\setcounter{equation}{0}
\section{Sufficient Conditions for Self-Propulsion in the case $\bar{\bfv_*}=\0$}
As \theoref{Exi}, also \theoref{6.1} furnishes a generic existence result of $T$-periodic solutions and, as such, it does not ensure $\bfxi_\sub{\bfu}\neq\0$, namely, that the body will perform a non-zero net motion. Our next task is to give sufficient conditions on $\bfv_*$ that, in fact, ensure this property, when $\bfv_*$ has zero average. In other words, we want 
to provide a sufficient condition on the   boundary data that --under fixed geometrical and physical properties of the body $\mathscr B$ and the liquid--  ensure that $\mathscr B$ is able to self-propel; see \eqref{SPC}.  As will be shown in Section 9, this condition is also necessary, in the sense that we will prove, by means of an explicit example, that if the  boundary data are not such as to satisfy that condition, the body $\mathscr B$  may just ``oscillate" without performing any non-zero net motion,  regardless of  its shape and physical properties.

To accomplish all the above, we  begin to define  the vector $\bfG$ characterized as follows: 
\be
\bfG=\bfG(\mathscr B,\bfV_*,M,\nu):=\sum_{i-1}^3\left(\int_{\Omega}\bar{(\bfV_0-\bfchi_0)\cdot\nabla\bfh^{({i})}\cdot \bfV_0}\right)\bfe_i\,,
\eeq{G}
where  $\bfV_0$ is the velocity field of the solution given in \lemmref{6.1}, and  the fields $(\bfh^{(i)},p^{(i)})$ are given in \eqref{6.9}. Observe that from \lemmref{6.1} and \eqref{6.9_0},  we deduce that $\bfG$ is well defined. Notice also that, for a given body and liquid,  $\bfG$ depends only on the boundary velocity distribution.   
\par
In order to show our  main finding, we need the following lemma.
\Bl Let the assumption of \theoref{6.1} be satisfied, and let $\bfU:=(\bfu,\bfw,\bfchi_\sub{\bfw})\in \mathfrak B^{q}\cap \mathfrak B^{\frac{3q}{3-q}}$ be the solution to  \eqref{6.10}--\eqref{6.13} there constructed. Further,  let $\bfscr G$ be a tensor field with  $\bfscr G\in L^3(\Omega)$. Then, the following inequality holds
$$
\left|\int_\Omega\tilde{\bfscr F}(\bfU):\bfscr G\right|\le c\,\delta^3\|\bfscr G\|_3\,,
$$
where $\tilde{\bfscr F}$ is defined in \eqref{Ftilde}.
\EL{FM}
{\em Proof.} From H\"older inequality, we have
\be
\left|\int_\Omega\tilde{\bfscr F}(\bfU):\bfscr G\right|\le \|\tilde{\bfscr F}\|_{\frac32}\|{\bfscr G}\|_3\,.
\eeq{a}
Next, we observe that, since $q\in (1,\frac65]$,  by simple interpolation we get
$$
\|\tilde{\bfscr F}\|_{\frac32}\le \|\tilde{\bfscr F}\|_{q}^\theta\|\tilde{\bfscr F}\|_{2}^{1-\theta}\,,
$$
with $\theta=\frac13\frac q{2-q}$.
From the latter and \lemmref{4.2}, we thus deduce
$$
\|\tilde{\bfscr F}\|_{\frac32}\le   c\,\left(\|\bfU\|_{\mathfrak B^{q}}^2 +\delta\,\|\bfU\|_{\mathfrak B^{q}}\right)\,,
$$
which, in view of \eqref{gatta0} and \eqref{a}, completes the proof of the lemma.
\par\hfill$\square$\ms\par
Our next  result shows that the vector $\bfG$  acts as ``thrust" for self-propulsion. More precisely, we have the following theorem.
\Bt Let $\bfv_*=\delta\,\bfV_*$ be as in \theoref{6.1}. Suppose that
\be
 \bfG\neq \0\,.
\eeq{SPC}
Then, there exists $\delta_0>0$ such that for all $\delta\in(0,\delta_0)$ the unique solution to \eqref{SE} constructed in \theoref{6.1} and corresponding to $\bfv_*$, must have $\bfxi_\sub{\bfu}\neq\0$.
\ET{6.2}
{\em Proof.} Let $(\bfu,\bfw,\bfchi_\sub{\bfw})\in \mathfrak B^q\cap\mathfrak B^{\frac{3q}{3-q}}$ be the solution given in  \theoref{6.1}  corresponding to $\bfv_*$.
and assume, ad absurdum, that $\bfxi_\sub{\bfu}=\0$. If we then  dot-multiply both sides of \eqref{6.10}$_1$ (with $\bfxi_\sub{\bfu}=\0$) by $\bfh^{(i)}$,  integrate by parts over $\Omega$ and use \eqref{6.10}$_4$, we get
\be
\nu\int_\Omega\mathbb D(\bfu):\mathbb D(\bfh^{(i)})=\int_\Omega\bfscr F:\nabla\bfh^{(i)}\,.
\eeq{GG1}
Likewise, if we dot-multiply both sides of \eqref{6.9}$_1$ by $\bfu$, integrate by parts over $\Omega$ and assume $\bfxi_\sub{\bfu}=\0$, we show
$$
\int_\Omega\mathbb D(\bfu):\mathbb D(\bfh^{(i)})=0\,,
$$
which, combined with \eqref{GG1}, furnishes
\be 
\int_{\Omega}\bfscr F:\nabla\bfh^{(i)}=0\,,\ \ \mbox{for {\em all} $i=1,2,3$}\,.
\eeq{6.20}
Let $G\equiv G_j\neq 0$, for some $j\in\{1,2,3\}$. From \eqref{6.20} we get 
\be
\left|\Int{\Omega}{}(\bfscr F:\nabla\bfh^{(j)}-\tilde{\bfscr F}:\nabla\bfh^{(j)})\right|=
\delta^2|G|
\eeq{6.21}
where
$\tilde{\bfscr F}$ is given in \eqref{Ftilde}.
Since, in particular, by  \eqref{6.9_0} we have  $\bfh^{(i)}\in D^{1,r}(\Omega)$, for all $r\in(\frac32,\infty)$, thanks to \lemmref{FM} we  show that
$$  
\left|\Int{\Omega}{}\tilde{\bfscr F}:\nabla\bfh^{(j)}\right|\le c_0\,\delta^3\,,
$$
which, once replaced in  \eqref{6.21} allows us to infer, in particular,
\be
\left|\Int{\Omega}{}\bfscr F:\nabla\bfh^{(j)}\right|\ge\delta^2(|G|-c_0\,\delta)\,.
\eeq{pot}
However, if we choose $\delta<|G|/c_0$ in \eqref{pot}, we contradict \eqref{6.20}  and thus, in turn,  the assumption $\bfxi_\sub{\bfu}=\0$. The proof of the theorem is therefore completed. \par\hfill$\square$\par
\Br It is clear, from \theoref{6.1} and \theoref{6.2}, that the set of boundary distributions assuring the validity of the self-propulsion condition \eqref{SPC} is open in the space $\calv^{2,q}_\sharp$. As a result, also the set of ``thrusts" $\bfG$ is open in $\real^3$.
\ER{6.3}
\setcounter{equation}{0}
\section{On the Velocity of Self-Propulsion in the case $\bar{\bfv_*}=\0$} In the previous section we have furnished sufficient conditions ensuring self-propulsion of the body $\mathscr B$. Our objective now is to give an estimate of the velocity of self-propulsion, $\bfxi_\sub{\bfu}$, and, in particular, its relation to the ``thrust" $\bfG$ defined in \eqref{G}.\par
To this end, we begin to show the following result.
\Bl Let the assumptions of \theoref{6.2} be satisfied,  and consider the  boundary-value problem:
\be\ba{cc}\medskip\left.\ba{ll}\medskip
\nu\Delta\tilde{\bfh}^{(i)}-\bfxi_\sub{\bfu}\cdot\nabla\tilde{\bfh}^{(i)}=\nabla \tilde{p}^{(i)}\\
\Div\tilde{\bfh}^{(i)}=0\ea\right\}\ \ \mbox{in $\Omega$}\\ 
\tilde{\bfh}^{(i)}=\bfe_i\ \ \mbox{at $\partial\Omega$}\,,\ \ \Lim{|x|\to\infty}\tilde{\bfh}^{(i)}(x)=\0\,.
\ea
\eeq{7.1_0}
Then, there exists a unique solution 
\be
\tilde{\bfh}^{(i)}\in L^{\frac{2\sigma}{2-\sigma}}(\Omega)\cap D^{1,\frac{4\sigma}{4-\sigma}}(\Omega)\cap D^{2,\sigma}(\Omega)\,,\ \ p\in D^{1,\sigma}(\Omega)\,,\ \ \ \mbox{all $\sigma\in(1,2)$}
\eeq{7.2_0}
that, in addition, satisfies the following properties ($i=1,2,3$):
\be\ba{ll}\medskip
\Int{\partial\Omega}{}\mathbb T(\tilde{\bfh}^{(i)},\tilde{p}^{(i)})\cdot\bfn=\Int{\partial\Omega}{}\mathbb T({\bfh}^{(i)},{p}^{(i)})\cdot\bfn+g_i(\delta)\,,\ \ |g_i(\delta)|\le c_0\,\delta^{\frac34}\,,\\
\tilde{\bfh}^{(i)}={\bfh}^{(i)}+\bfsf h^{(i)}(\delta)\,,\ \ |\bfsf h^{(i)}|_{1,3}\le c_0\,\delta^{\frac34}\,.
\ea
\eeq{7.2_00}
\EL{7.1_0}
{\em Proof.} Existence (and uniqueness) of the pair $(\tilde{\bfh}^{(i)},\tilde{p}^{(i)})$ in the class \eqref{7.2_0} is well known \cite[Theorem VII.7.1]{Gab}. In particular, such a solution satisfies the estimate
\be
|\bfxi_\sub{\bfu}|^{\frac14}|\tilde{\bfh}^{(i)}|_{1,r}\le C\,,\ \ \mbox{all  $r\in (\frac43,4)$}\,,
\eeq{7.3_0}
with $C=C(\Omega,r)$. Setting $\bfsf h^{(i)}:=\tilde{\bfh}^{(i)}-\bfh^{(i)}$, ${\sf p}^{(i)}=\tilde{p}^{(i)}-{p}^{(i)}$, from \eqref{7.1_0} and \eqref{6.9}, we deduce that 
\be\ba{cc}\medskip\left.\ba{ll}\medskip
\nu\Delta\bfsf h^{(i)}=\nabla {\sf p}^{(i)}+\bfxi_\sub{\bfu}\cdot\nabla\tilde{\bfh}^{(i)}\\
\Div\bfsf h^{(i)}=0\ea\right\}\ \ \mbox{in $\Omega$}\\ 
\bfsf h^{(i)}=\0\ \ \mbox{at $\partial\Omega $}\,.
\ea
\eeq{7.4_0}
Therefore,
from \eqref{7.3_0}, \eqref{7.4_0} and classical estimates on the Stokes problem, we get \cite[Theorem V.4.8]{Gab}
\be\ba{rl}\medskip
\|\bfsf h^{(i)}\|_{\frac{3r}{3-2r}}+|\bfsf h^{(i)}|_{1,\frac{3r}{3-r}}+|\bfsf h^{(i)}|_{2,r}+|{\sf p}^{(i)}|_{1,r}&\!\!\!\le c\,|\bfxi_\sub{\bfu}|\,|\tilde{\bfh}^{(i)}|_{1,r}\le c\,\delta^{\frac34}\,,\ \ \mbox{all $r\in(\frac43,\frac32)$}\,.
,\\
|\bfsf h^{(i)}|_{2,s}&\!\!\!\le  c\,|\bfxi_\sub{\bfu}|\,|\tilde{\bfh}^{(i)}|_{1,s} \le c\,\delta^{\frac34}\,,\ \ \mbox{all $s\in(3,4)$}\,.\ea
\eeq{7.5_0}
The first  inequality in \eqref{7.5_0}along with trace theorems, proves  \eqref{7.2_00}$_1$. Moreover, by classical embedding results for homogeneous Sobolev spaces \cite[Theorem II.9.1]{Gab}, we have
$$
|\bfsf{h}^{(i)}|_{1,3}\le c\,\left(|\bfsf{h}^{(i)}|_{1,\frac{3r}{3-r}}+|\bfsf{h}^{(i)}|_{2,s}\right)\,,
$$
so that \eqref{7.2_00}$_2$ is a consequence of the latter and \eqref{7.5_0}. 
\par\hfill$\square$\par 
\smallskip\par
We are now in a position to show the main result of this section.
\Bt Under the assumptions of  \theoref{6.2}, the velocity of self-propulsion $\bfxi_\sub{\bfu}$ is (non-zero and) given by the following formula
\be
\bfxi_\sub{\bfu}=\delta^2\mathbb A^{-1}\cdot\bfG+\bfsigma(\delta)\,,
\eeq{7.6_0}
where $\mathbb A$ is the symmetric, nonsingular matrix defined in \eqref{A}, 
and $\bfsigma$ is vector function depending also on $\delta$ such that 
$$
\bfsigma(\delta)=O(\delta^{\frac{11}4}).
$$
\ET{7.1_0}
{\em Proof.}  We  dot-multiply both sides of \eqref{6.10}$_1$ by $\tilde{\bfh}^{(i)}$,  integrate by parts over $\Omega$ and employ \eqref{6.10}$_4$ to show
\be 
\Int{\Omega}{}\left[\nu\,\mathbb D(\bfu):\mathbb D(\tilde{\bfh}^{(i)})+\bfxi_\sub{\bfu}\cdot\nabla \tilde{\bfh}^{(i)}\cdot\bfu-\bfscr F:\nabla \tilde{\bfh}^{(i)}\right]=0\,.
\eeq{7.7_0}
In a similar way, dot-multiplying both sides of \eqref{7.1_0}$_1$ by $\bfu$ and proceeding as before, we get
\be
\bfxi_\sub{\bfu}\cdot\Int{\partial\Omega}{}\mathbb T(\tilde{\bfh}^{(i)},\tilde{p}^{(i)})\cdot\bfn-
\Int{\Omega}{}\left[\nu\,\mathbb D(\bfu):\mathbb D(\tilde{\bfh}^{(i)})+\bfxi_\sub{\bfu}\cdot\nabla \tilde{\bfh}^{(i)}\cdot\bfu\right]\,,
\eeq{7.8_0}
so that, summing \eqref{7.7_0} and \eqref{7.8_0} side-by-side, we conclude
\be
\bfxi_\sub{\bfu}\cdot\Int{\partial\Omega}{}\mathbb T(\tilde{\bfh}^{(i)},\tilde{p}^{(i)})\cdot\bfn=\int_\Omega \bfscr F:\nabla \tilde{\bfh}^{(i)}\,,\ \ i=1,2,3\,.
\eeq{7.9_0}
From \lemmref{7.1_0},  \eqref{WH2} and \eqref{gatta0}  it immediately follows  that 
\be
\bfxi_\sub{\bfu}\cdot\Int{\partial\Omega}{}\mathbb T(\tilde{\bfh}^{(i)},\tilde{p}^{(i)})\cdot\bfn=\xi_\sub{\bfu k}\mathbb A_{ki} + O(\delta^{\frac{11}4})\,.
\eeq{7.11_0}
Further, we have
\be
\int_\Omega \bfscr F:\nabla \tilde{\bfh}^{(i)}=\delta^2 G_i+\int_\Omega \tilde{\bfscr F}:\nabla {\bfh}^{(i)}+\int_\Omega \tilde{\bfscr F}:\nabla {\bfsf h}^{(i)}+\delta^2\Int{\Omega}{}\bfscr F_0:\nabla \bfsf{h}^{(i)}\,,
\eeq{clo}
where, we recall, the quantities $\tilde{\bfscr F}$ and $\bfscr F_0$ are defined in \eqref{F_0} and \eqref{Ftilde}.
From \lemmref{FM} and  \lemmref{7.1_0},
\be
\int_\Omega \tilde{\bfscr F}:\nabla {\bfh}^{(i)}+\delta^2\Int{\Omega}{}\bfscr F_0:\nabla \bfsf{h}^{(i)}=O(\delta^{\frac{11}4})\,,
\eeq{7.12_0}
and, 
\be
\int_\Omega \tilde{\bfscr F}:\nabla {\bfsf h}^{(i)}=O(\delta^3)\,.
\eeq{7.13_0}
The theorem then follows from \eqref{7.9_0}--\eqref{7.13_0}.
\par\hfill$\square$\par
\setcounter{equation}{0}
\section{Some Comments on the Self-Propelling Condition \eqref{SPC}}
From \theoref{7.1_0}, it turns out that if $\bar{\bfv_*}=\0$,  self-propulsion  manifests itself at the second order in $\delta$, {\em provided $\bfG\neq\0$}. The natural question to ask is then the  following one. {\em Suppose $\bfG=\0$. Can then self-propulsion  occur at an order in $\delta$  higher than 2}?
 The aim of this section is to prove that, in general, the answer is {\em negative}. In fact, we shall provide an example of boundary data,$\bfv_*$, of arbitrary magnitude $\delta>0$, for which the corresponding solution to \eqref{6.1} does {\em not} satisfy \eqref{SPC} and the averaged field $\bfu$ in problem \eqref{6.10} is {\em identically zero}. More precisely, our example shows that given a body $\mathscr B$ of {\em any shape and mass}, and an arbitrary period $T>0$, we can always find a $T$-periodic boundary velocity $\bfv_*$ such that \eqref{SPC} is violated {\em and} the net motion of $\mathscr B$ is zero, that is, $\mathscr B$ can only ``oscillate." In order to show all the above, 
we premise the following result.
\Bl Let $\psi\in D^{1,2}(\Omega)$ be the solution to the  Neumann problem:
\be
\Delta\psi=0\ \ \mbox{in $\Omega$}\,;\ \ \pde{\psi}{\bfn}=f\ \ \mbox{at $\partial\Omega$}\,;\ \ \lim_{|x|\to\infty}\psi(x)=0
\eeq{7.1}
where $f$ is a given smooth function satisfying 
\be
\int_{\partial\Omega}f=0\,.
\eeq{7.2}
Moreover, let $a\in W^{1,q}(0,T)$, $q\in(1,\infty)$,  be  $T$-periodic  with $\bar{a}=0$. Then, the triple
\be
\bfsf V=a(t)\,\nabla\psi(x)\,,\ \ {\sf P}=-\dot{a}\,\psi\,,\ \ \bfsf z=\frac{a(t)}M\int_{\partial\Omega}\psi\,\bfn
\eeq{7.3}
is a solution to \eqref{6.1} in the class $\hat{\calw}_\sharp^{2,q}\times\calp^{1,q}\times W^{1,q}_\sharp(0,T)$, corresponding to the boundary data 
\be
\bfsf V_*:=a(t)\left(\nabla\psi|_{\partial\Omega}-\frac1M\int_{\partial\Omega}\psi\,\bfn\right).
\eeq{BoDa} 
Finally, for this solution we have
\be
\bfG=\0\,.
\eeq{7.4}
\EL{7.1}
{\em Proof.} We begin to observe that,  from classical results on the   exterior Neumann problem, we have $\psi\in C^\infty(\Omega)\cap W^{2,q}(\Omega_R)$, for all $q\in (1,\infty)$ and $R>R_*$. Furthermore, in view of \eqref{7.2}, it follows that \be D^\alpha\psi=O(|x|^{-2-\alpha})\,, \ \ |\alpha|=0,1,\ldots\,;
\eeq{p6} 
see \cite[Exercise V.3.6]{Gab}.  As a consequence, the fields \eqref{7.3} are in the class $ \hat{\calw}_\sharp^{2,q}\times\calp^{1,q}\times W^{1,q}_\sharp(0,T)$. Moreover, since $\psi$ is harmonic, it immediately follows that $(\bfsf V,{\sf P})$ is a solution to  \eqref{6.1}. Next, we have (in the trace sense)
$$
\int_{\partial\Omega}\mathbb T(\bfsf V,{\sf P})\cdot\bfn=\int_{\partial\Omega}[2a\,\nabla(\nabla\psi)\cdot\bfn+\dot{a}\psi\,\bfn]:= 2a\bfI+\int_{\partial\Omega}\dot{a}\,\psi\,\bfn\,.
$$
By integrating by parts and using the fact that $\psi$ is harmonic, we get, for $i=1,2,3$,
$$ 
I_{i}=\int_{\partial\Omega}\partial_i\partial_j\psi\,n_j=\int_{\Omega_R}\partial_i(\partial_j\partial_j\psi)-\int_{\partial B_R}\partial_i\partial_j\psi\,n_j=-\int_{\partial B_R}\partial_i\partial_j\psi\,n_j\,,
$$
so that, letting $R\to\infty$ into this relation and using the asymptotic properties of $\psi$, we conclude
$$
\bfI=\0\,.
$$
Therefore, \eqref{6.1}$_{3,4}$ are also satisfied if we choose $\bfsf z$ as in \eqref{7.3}$_3$, and $\bfsf V_* $ as stated. Next, we have
$$
G^{(i)}=\frac1T\int_0^Ta(t)\left(\int_{\Omega}({a(t)}\partial_k\psi-{\sf z}_{k}(t))\partial_k h_\ell^{(i)}\partial_\ell\psi\right)\,{\rm d}t:=\frac1T\int_0^Ta(t)I(t)\,{\rm d}t
.
$$
Thus,  integrating  by parts over $\Omega$, we get
$$\ba{rl}\medskip
I(t)=&\!\!\!\Int{\Omega}{}\partial_k[(a\,\partial_k\psi-{\sf z}_k)h_\ell\partial_\ell\psi]-\Int{\Omega}{}[(a\,\partial_k\psi-{\sf z}_k)\partial_\ell\partial_k\psi]h_\ell\\
=&\!\!\!\Int{\partial\Omega}{}n_k(a\,\partial_k\psi-{\sf z}_k)\bfe_i\cdot\nabla\psi-\Int{\partial\Omega}{}\bfe_i\cdot\bfn(\half \,a\,|\nabla\psi|^2-\bfsf z\cdot\nabla\psi):=I_1(t)-I_2(t)\,.
\ea
$$
However, using the fact that $\psi$ is harmonic in conjunction with Gauss theorem, we deduce
$$
I_1(t)=\Int{\Omega}{}\partial_k\big[(a\,\partial_k\psi-{\sf z}_k)\bfe_i\cdot\nabla\psi\big]=\Int{\Omega}{}\partial_i\big[\half a\,|\nabla\psi|^2-{\sf z}\cdot\nabla\psi\big]=I_2(t)\,,
$$
which thus proves \eqref{7.4}
\par\hfill$\square$\medskip\par
Before proceeding further, we believe it is important to emphasize that, by \lemmref{1.6_0},  for the indicated  boundary data, the solution provided in \eqref{7.3} is the {\em only one} in the relevant function class. 
We shall now show,
 with the help of \lemmref{7.1},  that the functions
\be
\bfv:= \bfsf V \,,\ \ p:=a(t)\bfgamma(t)\cdot\nabla\psi-\half\,a^2(t)|\nabla\psi|^2+{\sf P}\,,\ \ \bfgamma:= \bfsf z
\eeq{VV}
satisfy \eqref{SE} and $\bfv_*:=\bfsf V_*$ defined in \eqref{BoDa}.
In fact, observing that
$$
(\bfv-\bfgamma)\cdot\bfv=\nabla\Big[a(t)\bfgamma(t)\cdot\nabla\psi-\half\,a^2(t)|\nabla\psi|^2\Big]\,,
$$
it is immediately checked that the fields $\bfv$ and $p$ thus defined satisfy \eqref{SE}$_{1,2,4}$.
Let us now turn to the surface integral in \eqref{SE}$_5$  that, with $\bfv$ and $p$ defined above,  becomes
\be\ba{ll}\smallskip
\Int{\partial\Omega}{}\Big[2a\,\nu\,\nabla(\nabla\psi)\cdot\bfn-\dot{a}\psi\,\bfn-a\,(\bfgamma\cdot\nabla\psi\,\bfn-\nabla\psi\,\bfgamma\cdot\bfn)\\
\hspace*{2.7cm}+a^2(\half|\nabla\psi|^2\bfn-\nabla\psi\nabla\psi\cdot\bfn)\Big]
:=2a\,\nu\,\bfI_1-\dot{a} \bfI_2+a\,\bfI_3+a^2\bfI_4\,.
\ea
\eeq{p7}
As in the proof of \lemmref{7.1},   we obtain
\be
\bfI_1=\0\,.
\eeq{p8}
Next, we have
$$
\bfgamma\cdot\nabla\psi\,\bfn-\nabla\psi\,\bfgamma\cdot\bfn=(\bfn\times\nabla\psi)\times \bfgamma\,,
$$
and 
$$
\int_{\partial\Omega}\bfn\times\nabla\psi=-\int_{\partial B_R}\bfn\times\nabla\psi+\int_{\Omega_R}\nabla\times(\nabla\psi)=-\int_{\partial B_R}\bfn\times\nabla\psi\,.
$$
As a result, letting $R\to\infty$ in the latter and employing \eqref{p6} we deduce
\be
\bfI_3=\0\,.
\eeq{p9}
Finally, again by integration by parts, $i=1,2,3$,
$$\ba{rl}\medskip
I_{4i}=&\!\!\!\!\Int{\Omega_R}{}\big[\half\partial_i|\nabla\psi|^2-\partial_\ell(\partial_i\psi\,\partial_\ell\psi)\big]-\Int{\partial B_R}{}\big(\half|\nabla\psi|^2n_i-\partial_i\psi\partial_\ell\psi\,n_\ell\big)\\
=&\!\!\!\!
-\Int{\partial B_R}{}\big(\half|\nabla\psi|^2n_i-\partial_i\psi\partial_\ell\psi\,n_\ell\big)
\ea
$$
and again by \eqref{p6}, we may let $R\to\infty$ in this relation and infer
\be
\bfI_4=\0\,.
\eeq{p10}
Collecting \eqref{p7}--\eqref{p10}, we then infer that the functions $\bfv$, $p$, and $\bfchi$ defined by \eqref{VV} solve \eqref{SE}. 
However, for such a solution we have $\bfu:=\bar{\bfv}\equiv\0$, that is,  {\em self-propulsion cannot occur}. We also notice that the solution \eqref{VV} can be written in such a form  as to satisfy  the hypotheses of \theoref{6.1}, thus, in this case, providing {\em the only solution} corresponding to the given $\bfv_*$, for ``small" $\delta$. To show this, it suffices to define the fields
$$
\bfV_0:=\frac1\delta\,\bfsf V\,,\ \ P_0:=\frac1\delta\,{\sf P}\,,\ \ \bfchi_0:=\frac1\delta\,\bfsf z\,,\ \ \bfV_*:=\frac1\delta\,\bfsf V_*\,,
$$
and rewrite \eqref{VV} as follows
$$
\bfv:=\delta\,\bfV_0\,,\ \ p:=a(t)\bfgamma(t)\cdot\nabla\psi-\half\,a^2(t)|\nabla\psi|^2+\delta\,P_0\,,\ \ \bfgamma:=\delta\bfchi_0\,,
$$
which, of course, is in the form \eqref{ABP} for arbitrary $\delta>0$, with $\bfw\equiv\bfchi\equiv\bfu\equiv{\sf p}\equiv 0$ and $\tau:=a(t)\bfgamma(t)\cdot\nabla\psi-\half\,a^2(t)|\nabla\psi|^2$.
\Br
The example furnished in the previous remark also shows that, in general, if the boundary data have zero average, then a non-zero thrust $(\bfG\neq\0)$ is expected to  be produced by a boundary velocity distribution possessing a non-vanishing mass flow-rate through the body-liquid interface.  
\ER{6.2}
\medskip\par
{\bf Acknowledgment.} {I would like to thank Mr. Jan A. Wein for sharing several conversations on the topic of self-propulsion.}

\ed

Notice that any solution in the class $\mathfrak B$ to \eqref{4.1}--\eqref{4.2} with $\bfxi_\sub{\bfu}=\0$ is necessarily a solution to \eqref{5.4}--\eqref{5.5} in the class $\mathfrak C_0$. We shall show, however, that for ``almost all" $\calf\in \mathfrak L$ problem \eqref{5.4}--\eqref{5.5} is {\em not} well-posed in the class $\mathfrak C_0$, and this will imply that for ``almost all" $\calf\in \mathfrak L$ of ``small" norm, the corresponding solution to \eqref{4.1}--\eqref{4.2}  must have $\bfxi_\sub{\bfu}\neq\0$. In turn, the above lack of well-posedness  is secured if we show that the operator $\calm$ is Fredholm of {\em negative} index. Since from \lemmref{Op} and \lemmref{1.6} it follows that the index of $\cala$, ${\rm ind}\,\cala$, is $-3$. in order to show the same property for $\calm$ it suffices to prove that, at any given $\bfU_0\in \mathfrak C_0$,  the Fr\'echet derivative of $\caln+\call$ is compact.

Before performing this study, however, we would like to emphasize that whatever the shape of $\mathscr B$, and physical properties of $\mathscr B$ and $\mathscr L$, we can always find a boundary distribution that {\em does not} produce a non-zero net motion. This implies that, for a given body $\mathscr B$, some restrictions on the data are indeed necessary to guarantee  self-propulsion.  
Actually, let $\psi=\psi(x)$,  $x\in\Omega$, be a scalar function satisfying the following Neumann problem
\be
\Delta\psi=0\ \ \mbox{in $\Omega$}\,,\ \ \pde{\psi}{n}=f\ \ \mbox{at $\partial\Omega$}\ \ \ \mbox{with}\ \,\int_{\partial\Omega} f=0\,.
\eeq{p5}
Under condition \eqref{p5}$_3$ on the data, it is then well known that (for example, \cite[Exercise V.3.6]{Gab})
\be
D^\alpha \psi(x)=O(|x|^{-2-|\alpha|})\ \ \mbox{as $|x|\to\infty$}\,,\ \ |\alpha|\ge 0.
\eeq{p6}
Next, let $a=a(t)$ be a smooth $T$-periodic function with $\bar{a}=0$, and set
\be
\bfv(x,t):=a(t)\,\nabla\psi(x)\,.
\eeq{nb}
It is immediately checked that, in view of \eqref{p5}$_1$ and \eqref{p6}, the field $\bfv$ thus defined satisfies \eqref{SE}$_{1,2,4}$ with $ {\bfscr H}\equiv\0$, provided we choose 
\be 
p= -\dot{a}(t)\psi+a(t)\bfgamma(t)\cdot\nabla\psi-\half a^2(t)|\nabla\psi|^2\,.
\eeq{nb1}
Let us now turn to the surface integral in \eqref{SE}$_5$ (with $\bfscr{\bfscr H}\equiv\0$) that, with $\bfv$ and $p$ defined above becomes
\be\ba{ll}\smallskip
\Int{\partial\Omega}{}\Big[2a\,\nu\,\nabla(\nabla\psi)\cdot\bfn-\dot{a}\psi\,\bfn-a\,(\bfgamma\cdot\nabla\psi\,\bfn-\nabla\psi\,\bfgamma\cdot\bfn)\\
\hspace*{2.7cm}+a^2(\half|\nabla\psi|^2\bfn-\nabla\psi\nabla\psi\cdot\bfn)\Big]
:=2a\,\nu\,\bfI_1-\dot{a} \bfI_2+a\,\bfI_3+a^2\bfI_4\,.
\ea
\eeq{p7}
By a simple integration by parts and \eqref{p5}$_1$, we find, for $i=1,2,3$,
$$ 
I_{1i}=\int_{\partial\Omega}\partial_i\partial_j\psi\,n_j=\int_{\Omega_R}\partial_i(\partial_j\partial_j\psi)-\int_{\partial B_R}\partial_i\partial_j\psi\,n_j=-\int_{\partial B_R}\partial_i\partial_j\psi\,n_j\,,
$$
so that, letting $R\to\infty$ into this relation and using \eqref{p6}, we find
\be
\bfI_1=\0\,.
\eeq{p8}
We next observe that
$$
\bfgamma\cdot\nabla\psi\,\bfn-\nabla\psi\,\bfgamma\cdot\bfn=(\bfn\times\nabla\psi)\times \bfgamma\,,
$$
and that
$$
\int_{\partial\Omega}\bfn\times\nabla\psi=-\int_{\partial B_R}\bfn\times\nabla\psi+\int_{\Omega_R}\nabla\times(\nabla\psi)=-\int_{\partial B_R}\bfn\times\nabla\psi\,.
$$
As a result, letting $R\to\infty$ in the latter and employing \eqref{p6} we deduce
\be
\bfI_3=\0\,.
\eeq{p9}
Finally, again by integration by parts, $i=1,2,3$,
$$\ba{rl}\medskip
I_{4i}=&\!\!\!\!\Int{\Omega_R}{}\big[\half\partial_i|\nabla\psi|^2-\partial_\ell(\partial_i\psi\,\partial_\ell\psi)\big]-\Int{\partial B_R}{}\big(\half|\nabla\psi|^2n_i-\partial_i\psi\partial_\ell\psi\,n_\ell\big)\\
=&\!\!\!\!
-\Int{\partial B_R}{}\big(\half|\nabla\psi|^2n_i-\partial_i\psi\partial_\ell\psi\,n_\ell\big)
\ea
$$
and again by \eqref{p6}, we may let $R\to\infty$ in this relation and infer
\be
\bfI_4=\0\,.
\eeq{p10}
Collecting \eqref{p7}--\eqref{p10}, we then conclude that  a solution to \eqref{SE}  (with $\bfscr H\equiv\0$) is given by the fields $\bfv,p$ defined in \eqref{nb}--\eqref{nb1} together with 
$$
\bfgamma(t):=\frac{a(t)}M\int_{\partial\Omega}\psi\, \bfn\,,\ \ \bfv_*:=a(t)\,\nabla\psi-\bfgamma\,.
$$
Since, of course, $\bar{\bfgamma}=0$, no net motion occurs. Notice that $(\bfv,p,\bfgamma)$ is in the functional class $\mathcal W^{2,\frac32}_\sharp\times\calp^{1,\frac32}\times W^{1,2}_\sharp(0,T)$. 
 

Set 

Then, with the help of \lemmref{7.1}, we shall  show that the functions
\be
\bfv:=\delta\,\bfV_0\,,\ \ p:=a(t)\bfgamma(t)\cdot\nabla\psi-\half\,a^2(t)|\nabla\psi|^2+\delta\,P_0\,,\ \ \bfgamma:=\delta\bfchi_0
\eeq{VV}
satisfy \eqref{SE} with  ${\bfscr H}\equiv\0$ and $\bfv_*:=\delta\,\bfV_*\equiv \bfsf V_*$ defined in \lemmref{7.1}.

to the whole of $\real^3$ by a procedure entirely analogous to that used earlier on. Precisely, let $\bfscr S:=\nabla\bfP$, with $\bfP$ solving \eqref{NP} with $\bfscr F\equiv\bfscr H$ and set
$$
\tilde{\bfU}=\left\{\ba{ll}\medskip\bfU(x)\ &\mbox{if $x\in\Omega$}\\ 
\bfzeta\ &\mbox{if $x\in\Omega_0$}\ea\right.\,,\  \ \tilde{P}=\left\{\ba{ll}\medskip {P}(x)\ &\mbox{if $x\in\Omega$}\\ 
0\ &\mbox{if $x\in\Omega_0$}\ea\right.\,,\ \ \tilde{{\bfu}}(x)=\left\{\ba{ll}\medskip\bfu(x)\ &\mbox{if $x\in\Omega$}\\ 
{\bfxi_\sub{\bfu}}\ &\mbox{if $x\in\Omega_0$}\ea\right.\,,
$$ 
and
$$
\tilde{{\bfscr H}}(x)=\left\{\ba{ll}\medskip\bfscr{H}(x)\ &\mbox{if $x\in\Omega$}\\ 
{\bfscr S}(x)\ &\mbox{if $x\in\Omega_0$}\ea\right.\,.
$$
Then, arguing exactly as we did before, we show that $(\tilde{\bfU},\tilde{P})$ is a distributional solution to the problem
\be\left.\ba{ll}\medskip
\Delta\tilde{\bfU}+\bfxi_\sub{\bfu_1}\cdot\nabla\tilde{\bfU}=-\bfzeta\cdot\nabla\tilde{\bfu}+\nabla \tilde{P}+\Div\tilde{\bfscr H}\\ \Div\bfU=0
\ea\right\}\ \ \mbox{in $\real^3$}\,.
\eeq{St00}
\Bl
Assume   
\be
\Div  \bfscr{\bfscr G}\in \sfL^2\,, \ \ \bfb\in \call_\sharp^{2,\frac32}\,,\ \ \bfv_*\in W^{1,2}_\sharp\,,\ \ \bfg\in L^2_\sharp(0,T)\,.
\eeq{4.10}
Then, the following inequalities hold with $\bfU:=(\bfphi,\bfpsi,\bfchi_\sub{\bfpsi})$, 
\be\ba{ll}\medskip
\|\Div\bfscr{\bfscr F}(\bfU)\|_{\sfL^2}\le c\Big(\|\bfU\|_{\mathfrak B}^2+\|\bfU\|_{\mathfrak B}\|\bfv_*\|_{W^{1,2}(W^{\frac32,2}(\partial\Omega))}+\|\bfv_*\|_{W^{1,2}(W^{\frac32,2}(\partial\Omega))}^2+\|\Div\bfscr{\bfscr G}\|_{\sfL^2}\Big)\\ \medskip
\|\bff(\bfU)\|_{\call^{2,\frac32}}\le c\Big(\|\bfU\|_{\mathfrak B}^2\!\!+\!\|\bfv_*\|_{W^{1,2}(W^{\frac32,2}(\partial\Omega))}(1+\|\bfU\|_{\mathfrak B})+\|\bfv_*\|_{W^{1,2}(W^{\frac32,2}(\partial\Omega))}^2\!+\!\|\bfb\|_{\call^{2,\frac32}}\Big)\\
\|\bfF(\bfchi_\sub{\bfpsi},\bfxi_\sub{\bfphi})\|_{L^2(0,T)}\le
c\Big(\|\bfv_*\|_{W^{1,2}(W^{\frac32,2}(\partial\Omega))}(\|\bfU\|_{\mathfrak B}+1)+\|\bfv_*\|_{W^{1,2}(W^{\frac32,2}(\partial\Omega))}^2+\|\bfg\|_{L^2(0,T)}\Big)\,.
\ea
\eeq{4.11}
In particular, the range of the map $M$ is contained in $\mathfrak B$. 
\EL{4.2}
{\em Proof.} By assumption, $(\bfphi,\bfpsi,\bfchi_\sub{\bfpsi})\in \mathfrak B$, and, by  \lemmref{4.1}, $\bfV\in \calw_\sharp^{2,\frac32}$. Thus, using  \eqref{WH2}, it follows at once $\bfF(\bfchi_\sub{\bfpsi},\bfxi_\sub{\bfphi} )\in L^2_\sharp(0,T)$ and that it satisfies \eqref{4.11}$_3$. The proof of \eqref{4.11}$_{1,2}$ becomes also quite straightforward, if we make use of the following continuous embedding properties
\be 
{\sf D}^{2,\frac32}\subset L^r(\Omega)\,,\ \ \mathcal W_\sharp^{2,\frac32}\subset L^q(L^\infty)\,,\  \ \mbox{for all $r\in [3,\infty]$ and $q\in [1,4)$\,.}
\eeq{4.11_1}
The second embedding is a particular case of \cite[Theorem 2.1]{Mallo}. To show the first one, it is enough to prove it for $r=\infty$. However, the latter follows from  the well-known inequality
\be
\|\bfphi\|_\infty\le c\,\big(\|\bfphi\|_3+|\bfphi|_{2,2}\big)\,;
\eeq{4.11_2}
see \cite[Theorem II.6.1 and Theorem II.9.1]{Gab}. Once the inequalities in \eqref{4.11} have been established, from \lemmref{Fm} and \lemmref{1.6_0} applied to \eqref{4.3} and \eqref{4.4}, respectively, we deduce that $(\bfu,\bfw,\bfchi_\sub{\bfu}):=M(\bfU)$ is in $\mathfrak B$, which completes the proof of the lemma. \par\hfill$\square$\par
\begin{thebibliography}{99}
\bibitem{SC}\v{C}ani\'c, S., Moving Boundary Problems, {\em Bull. Amer. Math. Soc.} 
doi.org/10.1090/bull/1703 (2020) 
\bibitem{Court}Court, S., Existence of 3D strong solutions for a system modeling a deformable solid inside a viscous incompressible fluid, {\em J. Dynam. Differential Equations} {\bf 29}  (2017),  737--782
\bibitem{GaSP}Galdi, G.P., On the steady self-propelled motion of a body in a viscous incompressible fluid. {\em Arch. Ration. Mech. Anal}. {\bf 148} (1999),  53--88
\bibitem{GaRev}Galdi, G.P., An introduction to the Navier-Stokes initial-boundary value problem. {\em Fundamental directions in mathematical fluid mechanics}, 1--70, {\em Adv. Math. Fluid Mech.}\,, Birkh\"auser, Basel, (2000)
\bibitem{Garev}Galdi, G.P., On the motion of a rigid body in a viscous liquid: a mathematical analysis with applications. {\em Handbook of mathematical fluid dynamics}, Vol. I, 653--791, North-Holland, Amsterdam, 2002
\bibitem{GaNWP}Galdi, G.P., A steady-state exterior Navier--Stokes problem that is not well-posed. {\em Proc. Amer. Math. Soc.} {\bf 137} (2009) 679--684
\bibitem{Gab}Galdi, G.P., {\em An introduction to the mathematical theory of the Navier-Stokes equations.
Steady-state problems}, Second edition. Springer Monographs in Mathematics,
Springer, New York (2011)
\bibitem{GaTP} Galdi, G.P., On time-periodic flow of a viscous liquid past a moving cylinder. {\em Arch. Ration. Mech. Anal.} {\bf 210} (2013), 451--498
\bibitem{GaKy}Galdi, G.P., Kyed, M., Time--periodic flow of a viscous liquid past a body. {\em Partial differential equations in fluid mechanics}, 20--49, London Math. Soc. Lecture Note Ser., 452, Cambridge Univ. Press, Cambridge, 2018
\bibitem{GaMaH}
Galdi, G.P.,  Kyed, M., 
Time-periodic solutions to the Navier--Stokes equations.  {\em Handbook of mathematical analysis in mechanics of viscous fluids}, 509--578, Springer, Cham, 2018
\bibitem{GS1}Galdi, G.P.,  Silvestre A.L., Existence of time-periodic solutions to the Navier-Stokes equations around a moving body. {\em Pacific J. Math.}  {\bf 223}  (2006) 251--267   
\bibitem{GaSi} Galdi, G.P.,  Silvestre, A.L., On the motion of a rigid body in a Navier--Stokes liquid under the action of a time-periodic force. {\em Indiana Univ. Math. J.} {\bf 58} (2009), 2805--2842
\bibitem{HB}Happel, J., Brenner, H., {\em Low Reynolds number hydrodynamics with special applications to particulate media}. Prentice--Hall, Inc., Englewood Cliffs, N.J. 1965
\bibitem{Hey} Heywood, J.G.,  
The Navier-Stokes equations: on the existence, regularity and decay of solutions, {\it Indiana Univ. Math. J.}, {\bf 29}, 639--681 (1980)
\bibitem{Hishida} Hishida, T., Silvestre, A. L., and Takahashi, T., A boundary control problem for the steady self-propelled motion of a rigid body in a Navier-Stokes fluid. {\em  Ann. Inst. H. Poincar\'{e} Anal. Non Lin\'{e}aire} {\bf 34} (2017) 1507--1541
\bibitem{KoSo} Kozono, H., and Sohr, H., On stationary Navier--Stokes equations in unbounded domains. {\em Ricerche Mat.} {\bf 42} (1993), 69--86
\bibitem{MS}M\'{a}cha, V.,  Ne\v{c}asov\'{a}, \v{S}.,  Self-propelled motion in a viscous compressible fluid.{\em  Proc. Roy. Soc. Edinburgh} Sect. A {\bf 146} (2016) 415--433
\bibitem{Macha}M\'{a}cha, V.,  Ne\v{c}asov\'{a}, \v{S}., Self-propelled motion in a viscous compressible fluid--unbounded domains. {\em Math. Models Methods Appl. Sci}. {\em 26} (2016) 627--643
\bibitem{Necasova} Ne\v{c}asov\'{a}, \v{S}., Takahashi, T., and Tucsnak, M., Weak solutions for the motion of a self-propelled deformable structure in a viscous incompressible fluid. {\em Acta Appl. Math.} {\bf 116} (2011) 329--352
\bibitem{Raymond} Raymond, J.-P., Vanninathan, M., A fluid-structure model coupling the Navier-Stokes equations and the Lam\'{e} system. {\em J. Math. Pures Appl.} {\bf 102} (2014) 546--596 
\bibitem{MT}San Martín, J., Scheid, J--F., Takahashi, T., and  Tucsnak, M., An initial and boundary value problem modeling of fish-like swimming. {\em Arch. Ration. Mech. Anal.} {\bf 188} (2008),  429--455
\bibitem{Mallo}Solonnikov, V.A., Estimates of the solutions of the nonstationary Navier--Stokes system.  Boundary value problems of mathematical physics and related questions in the theory of functions,  {\em Zap. Naucn. Sem. Leningrad. Otdel. Mat. Inst. Steklov}. (LOMI) {\bf 38} (1973), 153--231
\bibitem{Staro}Starovoitov, V.N., Solvability of the problem of the self-propelled motion of several rigid bodies in a viscous incompressible fluid.{\em  Comput. Math. Appl.} {\bf 53} (2007)  413--435
\end{thebibliography}
